\documentclass{amsart}
\usepackage{amssymb}
\usepackage{amsfonts}

\newtheorem{theorem}{Theorem}
\theoremstyle{plain}
\newtheorem{acknowledgement}{Acknowledgement}

\newtheorem{corollary}{Corollary}

\numberwithin{equation}{section}
\input{tcilatex}

\begin{document}
\title[Quantum nonadapted calculus]{A nonadapted stochastic calculus and non
stationary evolution in Fock scale}
\author{V.P. Belavkin}
\address{M.I.E.M., B. Vusovski Street 3/12 Moscow\\
109028, USSR. M.I.E.M., 109028, USSR.}
\urladdr{}
\date{December 26, 1989}
\dedicatory{}
\thanks{This paper is published in: Quantum Probability and Related Topics 6
137--179 World Scientific, Singapore 1991.}

\begin{abstract}
A generalized definition of quantum stochastic (QS) integrals and
differentials is given in the free of adaptiveness and dimensionality form
in terms of Malliavin derivative on a projective Fock space, and their
uniform continuity with respect to the inductive limite convergence is
proved. A new form of QS calculus based on an inductive $\star $-algebraic
structure in an indefinite space is developed and a nonadaptive
generalization of the QS It\^{o} formula for its representation in Fock
space is derived. The problem of solution of general QS evolution equations
in a Hilbert space is solved in terms of the constructed operator
representation of chronological products, defined in the indefinite space,
and isometry and $\star $-homomorphism property respectively for operators
and maps of these solutions, corresponding to the peseudounitary and $\star $%
-homomorphism property of the QS integrable generators is proved.
\end{abstract}

\maketitle

\section{Introduction}

The noncommutative generalization of It\^{o} stochastic calculus developed
in [1--6], gives an adequate instrument of studying of the behavior of open
quantum dynamical systems in a singular coupling with Bose stochastic
fields. The quantum stochastic (QS) calculus enables us to solve the old
problem of the stochastic description of continuous collapse of the quantum
system under a continuous observation by using the stochastic theory of
quantum nondemolition measurements and filtering theory [7--9]. This gives
the examples of the stochastic nonunitary, nonstationary and even nonadapted
evolution equations in Hilbert space, the solution of which requires one to
define the chronologically ordered stochastic exponents of operators and
maps in an appropriate way.

Here we solve this general problem in the framework of a new QS calculus in
Fock space, based on the explicit definition of the QS integrals free of the
adaptedness restriction in a uniform inductive topology, given in [10]. We
derive the general (nonadapted) It\^{o} formula as a differential of the
Wick formula for the normal ordered products, represented in an inductive $%
\star $-algebra with respect to an indefinite metric structure. The QS
generalization of It\^{o} formula for adapted processes was obtained by
Hudson and Parhasarathy in [1], where the unitary QS evolution was
constructed for the case of time and field independent QS generators $L$.
They used the QS integral for an adapted operator-valued function $D_{t}$ as
the limit of It\^{o} integral sums in the weak operator topology, defined as
in classical case due to commutativity of forward QS differentials $d\Lambda
(t)=\Lambda (t+dt)-\Lambda (t)$ with $D_{t}$. In this approach the QS
evolution for nonstationary generating operators of QS differential
equations was obtained for some finite dimensional cases by Holevo [11].

An another definition of QS integrals, based on the Berezin-Bargman calculus
in terms of kernels of operators in Fock space was proposed by Maassen [3].
One can show that Maassen kernel calculus corresponds to the particular
cases of our QS calculus, which is given directly in terms of the Fock
representation of integrated operators, instead of kernels [3,4]. Using this
new calculus we construct also the explicit solution of the nonstationary,
non Markovian, even nonadapted QS Langevin equations for a QS differentiable
stochastic process in the sense [12,13] over a unital $\star $-algebra $%
\mathcal{A}\subseteq \mathcal{B}(\mathcal{H})$ as the Fock representation of
an recursively defined operator-valued process in a pseudo Hilbert space
with noninner QS-integrable generators. Such QS evolution in a Markovian
stationary case was constructed recently by Evans and Hudson [14] and in a
nonstationary case by Lindsay and Parthasarathy [15]. We shall obtain the
existence and uniqueness of Evans-Hudson flow in free dimensional Markovian
case by the estimating of the explicit solution in the introduced inductive
uniform topology under the natural integrability conditions of time
dependent structural coefficients.

\section{Nonadapted QS integrals and differentials}

Let $\mathcal{H}$ be a Hilbert space with probability vectors $h\in \mathcal{%
H}\ ,\ \Vert h\Vert =1$, of a quantum dynamical object, described at any
instant $t\in \mathbb{R}_{+}$ by the algebra $\mathcal{B}(\mathcal{H})$ of
all linear bounded operators $L$ in $\mathcal{H}$ with Hermitian involution $%
L\rightarrow L^{\ast }$ and identity operator $I$. Let $X$ be a Borel space
with a positive measure $dx$, say $X=\mathbb{R}_{+}\times \mathbb{R}^{d}$,
and let $\{\mathcal{E}(x),x\in X\}$ be a family of complex Euclidean
(Hilbert) spaces $\mathcal{E}(x)$, describing the quantum field (noise) at a
point $x\in X$ of a dimensionality $\dim \mathcal{E}(x)\leq \infty $,
usually identified with a space $\mathcal{E}$. We denote by $\mathcal{K}%
=\int^{\oplus }\mathcal{E}(x)dx$ the Hilbert integral of the field state
spaces $\mathcal{E}(x)$, that is the space of all square integrable
vector-functions $k\colon x\rightarrow k(x)\in \mathcal{E}(x)$, 
\begin{equation*}
\langle k|k\rangle =\int \Vert k(x)\Vert ^{2}dx<\infty \ ,\ \Vert k(x)\Vert
^{2}=\langle k|k\rangle (x),
\end{equation*}%
and by $\Gamma (\mathcal{K})$ the Fock space of symmetrical tensor-functions 
$k(x_{1},\dots ,x_{n})$, $n=0,1,\dots ,$ with values in $\mathcal{E}%
(x_{1})\otimes \dots \otimes \mathcal{E}(x_{n})$. Let us assume the absolute
continuity $\mathrm{d}x=\lambda (t,\mathrm{d}x)\mathrm{d}t$ with respect to
a measurable map $t\colon X\rightarrow \mathbb{R}_{+}$, say $t(x)=t$, $%
\lambda (t,\mathrm{d}x)=\mathrm{d}\mathbf{x}$ for $x=(t,\mathbf{x})\in 
\mathbb{R}_{+}\times \mathbb{R}^{d}$, such that 
\begin{equation*}
\int_{\Delta }f(t(x))\mathrm{d}x=\int_{0}^{\infty }f(t)\lambda (t,\Delta )%
\mathrm{d}t
\end{equation*}%
for any integrable $\Delta \subseteq X$ and essentially bounded function $%
f\colon \mathbb{R}_{+}\rightarrow \mathbb{C}$. Then one can represent the
Fock space $\Gamma (\mathcal{K})$ as the Hilbert integral $\mathcal{F}=\int_{%
\mathcal{X}}^{\oplus }\mathcal{E}^{\otimes }(\varkappa )\mathrm{d}\varkappa $
of the functions 
\begin{equation*}
k\colon \varkappa \rightarrow k(\varkappa )\in \mathcal{E}^{\otimes
}(\varkappa ),\mathcal{E}^{\otimes }(\varkappa )=\otimes _{x\in \varkappa }%
\mathcal{E}(x)
\end{equation*}%
over the set $\mathcal{X}$ of all finite chains $\varkappa =(x_{1},\dots
,x_{n})$, identified with the indexed subsets $\{x_{1},\dots ,x_{n}\}\subset
X$ of cardinality $|\varkappa |=n<\infty $ and $\mathrm{d}\varkappa
=\dprod_{x\in \varkappa }\mathrm{d}x$ under the order $t(x_{1})<\dots
<t(x_{n})$. We shall denote by $t(\varkappa )$ the chains (subsets) $%
\{t(x)|x\in \varkappa \}$, $\emptyset \in \mathcal{X}$ denotes the empty
chain and $1_{\emptyset }\in \mathcal{F}$ denotes the vacuum function: $%
1_{\emptyset }(\varkappa )=0$, if $\varkappa \not=\emptyset $; $1_{\emptyset
}(\emptyset )=1$.

This can be done as in the case $X=\mathbb{R}_{+},t(x)=x$ by the isometry 
\begin{equation*}
\sum_{n=0}^{\infty }{\frac{1}{n!}}\int_{X^{n}}\Vert k(\varkappa )\Vert ^{2}%
\mathrm{d}\varkappa =\sum_{n=0}^{\infty }\;\idotsint\nolimits_{t_{1}<\dots
<t_{n}}\Vert k(x_{1},\dots ,x_{n})\Vert ^{2}\mathrm{d}x_{1}\dots \mathrm{d}%
x_{n},
\end{equation*}%
where the integrals in right hand side is taken over all $\varkappa
=\{x_{1}<\dots <x_{n}\}$ with different $t_{i}=t(x_{i})$ due to 
\begin{equation*}
{\frac{1}{n!}}\int_{0}^{\infty }\dots \int_{0}^{\infty }f(t_{1},\ldots
,t_{n})\mathrm{d}t_{1}\cdots \mathrm{d}t_{n}=\int_{0}^{\infty }\mathrm{d}%
t_{1}\int_{t_{1}}^{\infty }\mathrm{d}t_{2}\dots \int_{t_{n-1}}^{\infty }%
\mathrm{d}t_{n}f(t_{1},t_{2},\dots ,t_{n})
\end{equation*}%
for the symmetrical function%
\begin{equation*}
f(t_{1},\dots ,t_{n})=\int_{X^{n}}\Vert k(\varkappa )\Vert
^{2}\dprod_{i=1}^{n}\lambda (t_{i},\mathrm{d}x_{i}).
\end{equation*}

One can consider the set $X$ as the space with a casual preorder $\lesssim $
[12], and the increasing map $t\colon x\lesssim x^{\prime }\Rightarrow
t(x)\leq t(x^{\prime })$ as the local time, if for any $x\in X$ and $%
t^{\prime }>t(x)$ there exists $x^{\prime }\in X$ such that $t(x^{\prime
})=t^{\prime }$ (As it is for the map $t(x)=t$ with respect to the Galilean
or Einsteinian order in space-time $X=\mathbb{R}_{+}\times \mathbb{R}^{%
\mathrm{d}}$).

Let us denote by $\mathcal{F}(\xi )=\int_{\mathcal{X}}^{\oplus }\xi
^{|\varkappa |}\mathcal{E}^{\otimes }(\varkappa )\mathrm{d}\varkappa $ for
all $\xi >0$ the Hilbert scale of Fock spaces $\mathcal{F}(\xi )\subseteq 
\mathcal{F}(\zeta )$, $\xi \geq \zeta $, defined by the scalar products 
\begin{equation*}
\Vert k\Vert ^{2}(\xi )=\sum_{n=0}^{\infty }\xi
^{n}\idotsint\nolimits_{0\leq t_{1}<\dots <t_{n}<\infty }\Vert k(x_{1},\dots
,x_{n})\Vert _{\xi }^{2}\mathrm{d}x_{1}\cdots \mathrm{d}x_{n}\ ,
\end{equation*}%
by $\mathcal{G}(\xi )=\mathcal{H}\otimes \mathcal{F}(\xi )$ the Hilbert
tensor products, by $\mathcal{G}^{+}=\mathcal{G}(\xi ^{+})$, $\mathcal{G}=%
\mathcal{G}(1)$, $\mathcal{G}_{-}=\mathcal{G}(\xi _{-})$ the Hilbert
subspaces $\mathcal{G}^{+}\subseteq \mathcal{G}\subseteq \mathcal{G}_{-}$
for some $\xi ^{+}\geq 1\geq \xi _{-}$, and let us note that any linear
operator $L\in \mathcal{B}(\mathcal{H})$ can be considered as $(\xi ^{+},\xi
_{-})$-continuous (bounded) operator $B:\mathcal{G}^{+}\rightarrow \mathcal{G%
}_{-}$ of the form $B=L\otimes \hat{1}$, where $\hat{1}$ means the identity
operator $\hat{1}=\int_{\mathcal{X}}^{\oplus }I^{\otimes }(\varkappa )%
\mathrm{d}\varkappa \equiv I^{\otimes }$ in $\mathcal{F}=\mathcal{F}(1)$, $%
I^{\otimes }(\varkappa )=\otimes _{x\in \varkappa }I(x)$, considered as the
identical map $\mathcal{F}(\xi ^{+})\rightarrow \mathcal{F}(\xi _{-})$.
Following [2,8] we define the QS integral $\Lambda ^{t}(\mathbf{D}%
)=\int_{0}^{t}\mathrm{d}\Lambda ^{s}(\mathbf{D})$ for a table $\mathbf{D}%
=(D_{\nu }^{\mu })_{\nu =0,+}^{\mu =-,0}$ of functions $\{D_{\nu }^{\mu
}(x),x\in X\}$ with values in continuous operators 
\begin{equation*}
D_{0}^{0}(x):\mathcal{G}^{+}\otimes \mathcal{E}(x)\rightarrow \mathcal{G}%
_{-}\otimes \mathcal{E}(x)\ ,\ D_{+}^{-}(x):\mathcal{G}^{+}\rightarrow 
\mathcal{G}_{-},\eqno(1.1a)
\end{equation*}%
\begin{equation*}
D_{+}^{0}(x):\mathcal{G}^{+}\rightarrow \mathcal{G}_{-}\otimes \mathcal{E}%
(x)\ ,\ \;D_{0}^{-}(x):\mathcal{G}^{+}\otimes \mathcal{E}(x)\rightarrow 
\mathcal{G}_{-}\ ,\eqno(1.1.b)
\end{equation*}%
as the sum $\Lambda ^{t}(\mathbf{D})=\sum_{\mu ,\nu }\Lambda _{\mu }^{\nu
}(t,D_{\nu }^{\mu })$ of the operators $\Lambda _{\mu }^{\nu }(t,D):a\in 
\mathcal{G}\mapsto \Lambda _{\mu }^{\nu }(t,D)a$, acting as 
\begin{eqnarray*}
&[\Lambda _{0}^{0}(t,D_{0}^{0})a](\varkappa )=\sum_{x\in \varkappa
^{t}}[D_{0}^{0}(x)\dot{a}(x)](\varkappa \backslash x)&(1.2a) \\
&[\Lambda _{0}^{+}(t,D_{+}^{0})a](\varkappa )=\sum_{x\in \varkappa
^{t}}[D_{+}^{0}(x)a](\varkappa \backslash x)&(1.2b) \\
&[\Lambda _{-}^{0}(t,D_{0}^{-})a](\varkappa )=\int_{X^{t}}[D_{0}^{-}(x)\dot{a%
}(x)](\varkappa )\mathrm{d}x&(1.2c) \\
&[\Lambda _{-}^{+}(t,D_{+}^{-}a](\varkappa
)=\int_{X^{t}}[D_{+}^{-}(x)a](\varkappa )\mathrm{d}x\ ,&(1.2\mathrm{d}) \\
&&
\end{eqnarray*}%
Here $\varkappa ^{t}=\varkappa \cap X^{t}\ ,\ X^{t}=\{x\in X:t(x)<t\}$,$\
\varkappa \backslash x=\{x^{\prime }\in \varkappa :x^{\prime }\not=x\}$, and 
$a\mapsto \dot{a}(x)$ is the point (Malliavin [16,17]) derivative $\mathcal{G%
}^{+}\rightarrow \mathcal{G}^{+}\otimes \mathcal{E}(x)$, evaluated in the
Fock representation almost everywhere as $[\dot{a}(x)](\varkappa )=a(x\sqcup
\varkappa )$, where $x\sqcup \varkappa =\{x,\varkappa :x\notin \varkappa \}$
is the disjoint union of the chains $x,\,\varkappa \in \mathcal{X}$. The
operator-functions (1.2) were defined in [1] as the limits of the QS It\^{o}
integral sums with respect to the gage, creation, annihilation, and time
processes respectively for the bounded adapted operator valued functions $%
D(x)=A(x)\otimes \hat{1}_{[t}$, where $t=t(x),\hat{1}_{[t}=I_{[t}^{\otimes }$
is the identity operator in $\mathcal{F}_{[t}=\int_{\mathcal{X}%
_{[t}}^{\oplus }\mathcal{E}^{\otimes }(\varkappa )\mathrm{d}\varkappa $,$\ 
\mathcal{X}_{[t}=\{\varkappa \in \mathcal{X}|t(\varkappa )\geq t\}$. As it
follows from theorem 1 in [9,10] the operators (1.2) are densely defined in $%
\mathcal{G}$ as $(\zeta ^{+},\zeta _{-})$-continuous operators $\mathcal{G}%
(\zeta ^{+})\rightarrow \mathcal{G}(\zeta _{-})$ for any $\zeta ^{+}>\xi
^{+}\ ,\ \zeta _{-}<\xi _{-}$ even for the nonadapted and unbounded $D$,
satisfying local QS-integrability conditions 
\begin{equation*}
\Vert D_{0}^{0}\Vert _{\xi ^{+},\infty }^{\xi _{-},t}<\infty ,\Vert
D_{+}^{0}\Vert _{\xi ^{+},2}^{\xi _{-},t}<\infty ,\Vert D_{0}^{-}\Vert _{\xi
^{+},2}^{\xi _{-},t}<\infty ,\Vert D_{+}^{-}\Vert _{\xi ^{+},1}^{\xi
_{-},t}<\infty \ ,\eqno(1.3)
\end{equation*}%
for all $t\in \mathbb{R}_{+}$ and some $\xi _{-}\ ,\ \xi ^{+}>0$, where 
\begin{equation*}
\Vert D\Vert _{\xi ^{+},p}^{\xi _{-},t}=\left( \int_{X^{t}}\left( \Vert
D(x)\Vert _{\xi ^{+}}^{\xi _{-}}\right) ^{p}\mathrm{d}x\right)
^{1/p},\;\;\Vert D\Vert _{\xi ^{+}}^{\xi _{-}}=\sup \{\Vert D\mathbf{a}\Vert
(\xi _{-})/\Vert \mathbf{a}\Vert (\xi ^{+})\}.
\end{equation*}

Let us now define the multiple QS integral 
\begin{equation*}
\Lambda _{\lbrack 0,t)}(B)=\sum_{n=0}^{\infty }\;\idotsint_{0\leq
t_{1}<\dots <t_{n}<t}\mathrm{d}\Lambda ^{t_{1},\dots ,t_{n}}(B)\equiv
\int_{0\leq \tau <t}\mathrm{d}\Lambda ^{\tau }(B)
\end{equation*}%
for the operator-valued function $B(\pmb\vartheta )$ on the table $\pmb%
\vartheta =(\vartheta _{\nu }^{\mu })_{\nu =0,+}^{\mu =-,0}$ of four subsets 
$\vartheta _{\nu }^{\mu }\in \mathcal{X}$ with values 
\begin{equation*}
B\left( 
\begin{matrix}
\vartheta _{0}^{-} & \vartheta _{+}^{-} \\ 
\vartheta _{0}^{0} & \vartheta _{+}^{0}%
\end{matrix}%
\right) :\mathcal{G}(\eta ^{+})\otimes \mathcal{E}^{\otimes }(\vartheta
_{0}^{-})\otimes \mathcal{E}^{\otimes }(\vartheta _{0}^{0})\rightarrow 
\mathcal{G}(\eta _{-})\otimes \mathcal{E}^{\otimes }(\vartheta
_{0}^{0})\otimes \mathcal{E}^{\otimes }(\vartheta _{+}^{0})\eqno(1.4)
\end{equation*}%
as the operators in $\mathcal{G}$ with the action 
\begin{equation*}
\lbrack \Lambda _{\lbrack 0,t)}(B)a](\varkappa )=\sum_{\vartheta
_{0}^{0}\sqcup \vartheta _{+}^{0}\subseteq \varkappa ^{t}}\int_{\mathcal{X}%
^{t}}\int_{\mathcal{X}^{t}}[B(\pmb\vartheta )\dot{a}(\vartheta
_{0}^{-}\sqcup \vartheta _{0}^{0})](\vartheta _{-}^{0})\mathrm{d}\vartheta
_{0}^{-}\mathrm{d}\vartheta _{+}^{-}\ .\eqno(1.5)
\end{equation*}

Here $\vartheta _{-}^{0}=\varkappa \cap \overline{(\vartheta _{0}^{0}\sqcup
\vartheta _{+}^{0})}=\varkappa \backslash \vartheta _{0}^{0}\backslash
\vartheta _{+}^{0}$ is the difference of a subset $\varkappa \subset X$ and
the partition $\vartheta _{0}^{0}\sqcup \vartheta _{+}^{0}$ as the disjoint
union $\vartheta _{0}^{0}\bigcup \vartheta _{+}^{0}\subseteq \varkappa $, $%
\vartheta _{0}^{0}\cap \vartheta _{+}^{0}=\emptyset $, $\mathcal{X}%
^{t}=\{\varkappa \in \mathcal{X}|\varkappa \subset X^{t}\}$, and the point
(Malliavin [16]) derivative 
\begin{equation*}
\dot{a}(\vartheta )=\int^{\oplus }a(\varkappa \sqcup \vartheta )\mathrm{d}%
\varkappa \in \mathcal{G}^{+}\otimes \mathcal{E}^{\otimes }(\vartheta )
\end{equation*}%
is defined for almost all $\varkappa \in \mathcal{X}$, $\varkappa \cap
\vartheta =\emptyset $ as $\dot{a}(\varkappa ,\vartheta )=a(\varkappa \sqcup
\vartheta )$ by a vector-function $a\in \mathcal{G}^{+}$. We shall say that
the function $B$ is locally QS integrable (in a uniform inductive limit), if
for any $t\in \mathbb{R}_{+}$ there exists a pair $(\eta ^{\bullet },\eta
_{\bullet })$ of triples $\eta ^{\bullet }=(\eta ^{-},\eta ^{0},\eta ^{+})$,$%
\;\eta _{\bullet }=(\eta _{-},\eta _{0},\eta _{+})$ of numbers $\eta ^{\mu
}>0$, $\eta _{\nu }>0$, for which $\Vert B\Vert _{\eta ^{\bullet }}^{\eta
_{\bullet }}(t)<\infty $, where 
\begin{equation*}
\Vert B\Vert _{\eta ^{\bullet }}^{\eta _{\bullet }}(t)=\int_{\mathcal{X}^{t}}%
\mathrm{d}\vartheta _{+}^{-}\left( \int_{\mathcal{X}^{t}}\int_{\mathcal{X}%
^{t}}{\frac{(\eta _{+})^{|\vartheta _{+}^{0}|}}{(\eta ^{-})^{|\vartheta
_{0}^{-}|}}}\sup_{\mathcal{X}^{t}}{\frac{(\eta _{0})^{|\vartheta _{0}^{0}|}}{%
(\eta ^{0})^{|\vartheta _{0}^{0}|}}}\left( \Vert B(\pmb\vartheta )\Vert
_{\eta ^{+}}^{\eta _{-}}\right) ^{2}\mathrm{d}\vartheta _{+}^{0}\mathrm{d}%
\vartheta _{0}^{-}\right) ^{1/2}\eqno(1.6)
\end{equation*}%
(sup is taken as essential supremum over $\vartheta _{0}^{0}\in \mathcal{X}%
^{t}$). As it follows from the next theorem, the function $B(\pmb\vartheta )$
in QS integral (1.5) can be defined up to the equivalence having the kernel $%
B\approx 0\Leftrightarrow \Vert B\Vert _{\eta ^{\bullet }}^{\eta _{\bullet
}}(t)=0$ for all $\eta _{\bullet },\eta ^{\bullet }$ and $t$. In particular,
one can define it only for the tables $\pmb\vartheta =(\vartheta _{\nu
}^{\mu })$, which are partitions $\varkappa =\sqcup \vartheta _{\nu }^{\mu }$
of the chains $\varkappa \in \mathcal{X}$, i.e. for $\pmb\vartheta =\sqcup
_{x\in \varkappa }\mathbf{x}$, where $\mathbf{x}$ means one from the four
single point (elementary) tables 
\begin{equation*}
\mathbf{x}_{0}^{0}=\left( 
\begin{matrix}
\emptyset  & \emptyset \cr x & \emptyset 
\end{matrix}%
\right) ,\ \mathbf{x}_{+}^{0}=\left( 
\begin{matrix}
\emptyset  & \emptyset \cr\emptyset  & x%
\end{matrix}%
\right) ,\ \mathbf{x}_{0}^{-}=\left( 
\begin{matrix}
x, & \emptyset \cr\emptyset  & \emptyset 
\end{matrix}%
\right) ,\ \mathbf{x}_{+}^{-}=\left( 
\begin{matrix}
\emptyset  & x\cr\emptyset  & \emptyset 
\end{matrix}%
\right) \ .
\end{equation*}

\begin{theorem}
If $B$ is a locally QS integrable function (1.4), then the multiple integral
(1.5) is a $(\xi ^{+},\xi _{-})$ continuous operator $U^{t}:\mathcal{G}(\xi
^{+})\rightarrow \mathcal{G}(\xi _{-})$ for $\xi ^{+}\geq \sum_{\mu }\eta
^{\mu },\xi _{-}^{-1}\geq \sum_{\nu }\eta _{\nu }^{-1}$, having the estimate 
\begin{equation*}
\Vert \Lambda _{\lbrack 0,t)}(B)a\Vert (\xi _{-})\leq \ \Vert B\Vert _{\eta
^{\bullet }}^{\eta _{\bullet }}\Vert a\Vert (\xi ^{+})\ ,\;\forall a\in 
\mathcal{G}(\xi ^{+})\ .
\end{equation*}%
The formally conjugated in $\mathcal{G}$ operator is defined as QS integral 
\begin{equation*}
\Lambda _{\lbrack 0,t)}(B)^{\ast }=\Lambda _{\lbrack 0,t)}(B^{{\star }}),B^{{%
\star }}(\pmb\vartheta )=B(\pmb\vartheta ^{{\star }})^{\ast },\pmb\vartheta
^{{\star }}=\left( 
\begin{matrix}
\vartheta _{+}^{0} & \vartheta _{+}^{-}\cr\vartheta _{0}^{0} & \vartheta
_{0}^{-}%
\end{matrix}%
\right) \ ,\eqno(1.7)
\end{equation*}%
which is the continuous operator $\mathcal{G}(1/\xi _{-})\rightarrow 
\mathcal{G}(1/\xi ^{+})$ with 
\begin{equation*}
\Vert \Lambda _{\lbrack 0,t)}(B^{{\star }})\Vert _{1/\xi _{-}}^{1/\xi
^{+}}=\Vert \Lambda _{\lbrack 0,t)}(B)\Vert _{\xi ^{+}}^{\xi _{-}}\equiv
\sup \{\Vert U^{t}a\Vert (\xi _{-})/\Vert a\Vert (\xi ^{+})\}\ .
\end{equation*}%
The QS process $U^{t}=\Lambda _{\lbrack 0,t)}(B)$ has a QS differential $%
\mathrm{d}U^{t}=\mathrm{d}\Lambda ^{t}(\mathbf{D})$ in the sense 
\begin{equation*}
\Lambda _{\lbrack 0,t)}(B)=B(\pmb\emptyset )+\Lambda ^{t}(\mathbf{D})\ ,\
D_{\nu }^{\mu }(x)=\Lambda _{\lbrack 0,t(x))}(\dot{B}(\mathbf{x}_{\nu }^{\mu
}))
\end{equation*}%
with $(\xi ^{+},\xi _{-})$-continuous QS derivatives $\mathbf{D}=(D_{\nu
}^{\mu })$, densely defined as in (1.5) by $\dot{B}(\mathbf{x},\pmb\vartheta
)=B(\mathbf{x}\sqcup \pmb\vartheta )$ for almost all $\pmb\vartheta
=(\vartheta _{\nu }^{\mu })$, where $\mathbf{x}=(\varkappa _{\nu }^{\mu })$
is one from the elementary tables $\mathbf{x}_{\nu }^{\mu }$, $\mu \not=+$, $%
\nu \not=-$ with $\varkappa _{\nu }^{\mu }=x$, and $\vartheta _{\nu }^{\mu
}\in \mathcal{X}^{t(x)}$.

Let $B(\pmb\vartheta )$ be defined for any partition $\varkappa =\sqcup
\vartheta _{\nu }^{\mu }\in \mathcal{X}$ as the solution $B(\pmb\vartheta )=%
\mathbf{L}^{\triangleleft }(\pmb\vartheta )\odot B(\pmb\emptyset )$ of the
recurrency 
\begin{equation*}
B(\mathbf{x}\sqcup \pmb\vartheta )=L(\mathbf{x})\odot B(\pmb\vartheta
),\vartheta _{\nu }^{\mu }\in \mathcal{X}^{t(x)}
\end{equation*}%
with $B(\pmb\emptyset )=T^{0}\otimes \hat{1}$, i.e. 
\begin{equation*}
\dot{B}(\mathbf{x},\pmb\vartheta )=(L(\mathbf{x})\otimes \hat{1})\cdot B(\pmb%
\vartheta )\ ,\ B(\pmb\emptyset )=T^{0}\otimes \hat{1}\eqno(1.8)
\end{equation*}%
with a table $\mathbf{L}=(L_{\nu }^{\mu })_{\nu =0,+}^{\mu =-,0}$ of
operator-valued functions $L_{\nu }^{\mu }(x)=L(\mathbf{x}_{\nu }^{\mu })$,%
\begin{eqnarray*}
L_{0}^{0}(x) &:&\mathcal{H}\otimes \mathcal{E}(x)\rightarrow \mathcal{H}%
\otimes \mathcal{E}(x)\ ,\;L_{+}^{-}(x):\mathcal{H}\rightarrow \mathcal{H}\
,(1.9a) \\
L_{+}^{0}(x) &:&\mathcal{H}\rightarrow \mathcal{H}\otimes \mathcal{E}(x)\
,\;L_{0}^{-}(x):\mathcal{H}\otimes \mathcal{E}(x)\rightarrow \mathcal{H}\
,(1.9b)
\end{eqnarray*}%
$L\odot B=(L\otimes \hat{1})\cdot B$, and 
\begin{equation*}
B(\pmb\varkappa )\cdot B(\pmb\vartheta )=(B(\pmb\varkappa )\otimes
I^{\otimes }(\vartheta _{0}^{0}\sqcup \vartheta _{+}^{0}))(B(\pmb\vartheta
)\otimes I^{\otimes }(\varkappa _{0}^{-}\sqcup \varkappa _{0}^{0}))\ ,
\end{equation*}%
where $I^{\otimes }(\varkappa )=\otimes _{x\in \varkappa }I(x),I(x)$ is the
identity operator in $\mathcal{E}(x)$. ($(B(\mathbf{x})\cdot B(\pmb\vartheta
)$ in (1.8) means usual product of operators in $\mathcal{G}$, if dim$%
\mathcal{E}=1$).

Then the process $U^{t}=\Lambda _{\lbrack 0,t)}(B)$ satisfies the QS
differential equation $\mathrm{d}U^{t}=\mathrm{d}\Lambda ^{t}(\mathbf{L}%
\odot U^{t})$ in the sense 
\begin{equation*}
U^{t}=U^{0}+\Lambda ^{t}(\mathbf{L}\odot U^{t})\ ,\ (\mathbf{L}\odot U)_{\mu
}^{\mu }(x)=(L(\mathbf{x}_{\nu }^{\mu })\otimes \hat{1})U^{t(x)}\ .\eqno%
(1.10)
\end{equation*}
\end{theorem}

\begin{proof}
$\;$Using~the sum-point integral$\;$property%
\begin{equation*}
\;\int \sum_{\sqcup \vartheta _{\nu }=\vartheta }f(\vartheta _{-},\vartheta
_{0},\vartheta _{+})\mathrm{d}\vartheta =\iiint f(\vartheta _{-},\vartheta
_{0},\vartheta _{+}){\dprod _{\nu }\mathrm{d}\vartheta _{\nu }}
\end{equation*}
of the multiple sum-point integral, we obtain from definition (1.5) for $%
a,c\in \mathcal{G}$: 
\begin{equation*}
\int \langle c(\vartheta )|[U^{t}a](\vartheta )\rangle \mathrm{d}\vartheta
=\int_{\mathcal{X}^{t}}\mathrm{d}\vartheta _{+}^{-}\int_{\mathcal{X}^{t}}%
\mathrm{d}\vartheta _{0}^{-}\int_{\mathcal{X}^{t}}\mathrm{d}\vartheta
_{+}^{0}\int_{\mathcal{X}^{t}}\mathrm{d}\vartheta _{0}^{0}\langle \dot{c}%
(\vartheta _{0}^{0}\sqcup \vartheta _{+}^{0})|B(\pmb\vartheta )\dot{a}%
(\vartheta _{0}^{-}\sqcup \vartheta _{0}^{0})\rangle =
\end{equation*}%
\begin{equation*}
\int_{\mathcal{X}^{t}}\mathrm{d}\vartheta _{+}^{-}\int_{\mathcal{X}^{t}}%
\mathrm{d}\vartheta _{0}^{-}\int_{\mathcal{X}^{t}}\mathrm{d}\vartheta
_{+}^{0}\int_{\mathcal{X}^{t}}\mathrm{d}\vartheta _{0}^{0}\langle B(\pmb%
\vartheta )^{\ast }\dot{c}(\vartheta _{0}^{0}\sqcup \vartheta _{+}^{0}))|%
\dot{a}(\vartheta _{0}^{-}\sqcup \vartheta _{0}^{0})\rangle =\int \langle
\lbrack U^{t{\ast }}c](\vartheta )|a(\vartheta )\rangle \mathrm{d}\vartheta
\ ,
\end{equation*}%
that is $U^{t{\ast }}$ acts as $\Lambda _{\lbrack 0,t)}(B^{{\star }})$ in
(1.5) with $B^{{\star }}(\pmb\vartheta )=B(\pmb\vartheta ^{{\star }})^{\ast }
$. Moreover this equation gives $\Vert \Lambda _{\lbrack 0,t)}(B)\Vert _{\xi
}^{1/\zeta }=\Vert \Lambda _{\lbrack 0,t)}(B^{{\star }})\Vert _{\zeta
}^{1/\xi }$ as 
\begin{equation*}
\Vert U\Vert _{\xi }^{1/\zeta }=\text{sup}|\langle c|Ua\rangle |/\Vert
a\Vert (\xi )\Vert c\Vert (\zeta )=\sup |\langle U^{\ast }c|a\rangle |/\Vert
c\Vert (\zeta )\Vert a\Vert (\xi )=\Vert U^{\ast }\Vert _{\zeta }^{1/\xi }\ .
\end{equation*}%
Let us estimate the integral $\langle c|U^{t}a\rangle $, using the Schwarz
inequality 
\begin{equation*}
\int \Vert \dot{c}(\vartheta )\Vert (\eta _{-}^{-1})\Vert \dot{a}(\vartheta
)\Vert (\eta _{+})(\eta _{0}/\eta ^{0})^{|\vartheta |/2}\mathrm{d}\vartheta
\leq \Vert \dot{c}\Vert (\eta _{-}^{-1},\eta _{0}^{-1})\Vert \dot{a}\Vert
(\eta ^{+},\eta ^{0})
\end{equation*}%
and the following isometricity property of the multiple derivative: 
\begin{equation*}
\Vert \dot{a}\Vert (\xi ,\eta )=\left( \iint \xi ^{|\vartheta |}\eta
^{|\sigma |}\Vert a(\vartheta \sqcup \sigma )\Vert ^{2}\mathrm{d}\vartheta 
\mathrm{d}\sigma \right) ^{1/2}=\Vert a\Vert (\xi +\eta )\ .
\end{equation*}%
This gives $|\langle c|U^{t}a\rangle |=|\int \langle c(\varkappa )|[\imath
_{\lbrack 0,t)}^{\otimes }(B)a](\varkappa )\rangle \mathrm{d}\varkappa \leq 
\newline
$ 
\begin{eqnarray*}
&\leq &\newline
\int_{\mathcal{X}^{t}}\mathrm{d}\vartheta _{+}^{-}\int_{\mathcal{X}^{t}}%
\mathrm{d}\vartheta _{0}^{-}\int_{\mathcal{X}^{t}}\mathrm{d}\vartheta
_{+}^{0}\int_{\mathcal{X}^{t}}\mathrm{d}\vartheta _{0}^{0}\Vert \dot{c}%
(\vartheta _{0}^{0}\sqcup \vartheta _{+}^{0})\Vert (\eta _{-}^{-1})|\,\Vert
B(\pmb\vartheta )\Vert _{\eta ^{+}}^{\eta _{-}}\Vert \dot{a}(\vartheta
_{0}^{-}\sqcup \vartheta _{0}^{0})\Vert (\eta ^{+})\newline
\\
&\leq &\int_{\mathcal{X}^{t}}\mathrm{d}\vartheta _{+}^{-}\int_{\mathcal{X}%
^{t}}\mathrm{d}\vartheta _{0}^{-}\int_{\mathcal{X}^{t}}\mathrm{d}\vartheta
_{+}^{0}\Vert \dot{c}(\vartheta _{+}^{0})\Vert (\eta _{-}^{-1}+\eta
_{0}^{-1})\Vert \dot{a}(\vartheta _{0}^{-})\Vert (\eta ^{+}+\eta ^{0})%
\newline
\Vert B\Vert _{\eta ^{+},\eta ^{0}}^{\eta _{-},\eta _{0}}(t,\pmb\vartheta
\backslash \pmb\vartheta _{0}^{0}) \\
&\leq &\Vert c\Vert (\sum_{\nu =-}^{+}\eta _{\nu }^{-1})\Vert a\Vert
(\sum_{\mu =-}^{+}\eta ^{\mu })\newline
\int_{\mathcal{X}^{t}}\mathrm{d}\vartheta _{+}^{-}\left[ \int_{\mathcal{X}%
^{t}}\mathrm{d}\vartheta _{0}^{-}\int_{\mathcal{X}^{t}}\mathrm{d}\vartheta
_{+}^{0}{\frac{(\eta _{+})^{|\vartheta _{+}^{0}|}}{(\eta ^{-})^{|\vartheta
_{0}^{-}|}}}\Vert B\Vert _{\eta ^{+},\eta ^{0}}^{\eta _{-},\eta _{0}}(t,\pmb%
\vartheta \backslash \pmb\vartheta _{0}^{0})^{2}\right] ^{\frac{1}{2}}
\end{eqnarray*}%
where $\Vert B\Vert _{\eta ^{+},\eta ^{0}}^{\eta _{-},\eta _{0}}(t,\pmb%
\vartheta \backslash \pmb\vartheta _{0}^{0})=\mathrm{esssup}_{\vartheta
_{0}^{0}\in \mathcal{X}^{t}}{\frac{(\eta _{0})^{|\vartheta _{0}^{0}|/2}}{%
(\eta ^{0})^{|\vartheta _{0}^{0}|/2}}}\Vert B(\pmb\vartheta )\Vert _{\eta
^{+}}^{\eta _{-}}$. Hence 
\begin{equation*}
|\langle c|U^{t}a\rangle |\leq \Vert c\Vert (\xi _{-}^{-1})\Vert a\Vert (\xi
^{+})\Vert B\Vert _{\eta ^{\bullet }}^{\eta _{\bullet }}(t)
\end{equation*}%
for $\xi ^{+}\geq \sum_{\mu }\eta ^{\mu },\xi _{-}^{-1}\geq \sum_{\nu }\eta
_{\nu }^{-1}$.

Using the definition (1.5) and the property 
\begin{equation*}
\int_{\mathcal{X}^{t}}f(\vartheta )\mathrm{d}\vartheta =f(\emptyset )+\int_{%
\mathcal{X}^{t}}\mathrm{d}x\int_{\mathcal{X}^{t(x)}}\dot{f}(x,\vartheta )%
\mathrm{d}\vartheta \ ,\ \dot{f}(x,\vartheta )=f(x\sqcup \vartheta )\ ,
\end{equation*}%
one can obtain 
\begin{equation*}
\lbrack (U^{t}-U^{0})a](\varkappa )=[(\imath _{\lbrack 0,t)}^{\otimes }(B)-B(%
\pmb\emptyset ))a](\varkappa )=
\end{equation*}%
\begin{equation*}
\int_{X^{t}}\mathrm{d}x\sum_{\vartheta _{0}^{0}\sqcup \vartheta
_{+}^{0}\subseteq \varkappa }^{t(\vartheta _{\nu }^{0})<t(x)}\int_{\mathcal{X%
}^{t(x)}}\mathrm{d}\vartheta _{+}^{-}\int_{\mathcal{X}^{t(x)}}\mathrm{d}%
\vartheta _{0}^{-}[\dot{B}(\text{$\mathbf{x}$}_{+}^{-},\pmb\vartheta )\dot{a}%
(\vartheta _{0}^{-}\sqcup \vartheta _{0}^{0})+\dot{B}(\text{$\mathbf{x}$}%
_{0}^{-},\pmb\vartheta )\dot{a}(x\sqcup \vartheta _{0}^{-}\sqcup \vartheta
_{0}^{0})](\vartheta _{-}^{0})
\end{equation*}%
\begin{equation*}
+\sum_{x\in \varkappa ^{t}}\sum_{\vartheta _{0}^{0}\sqcup \vartheta
_{+}^{0}\subseteq \varkappa }^{t(\vartheta _{\nu }^{0})<t(x)}\int_{\mathcal{X%
}^{t(x)}}\mathrm{d}\vartheta _{+}^{-}\int_{\mathcal{X}^{t(x)}}\mathrm{d}%
\vartheta _{0}^{-}[\dot{B}(\text{$\mathbf{x}$}_{+}^{0},\pmb\vartheta )\dot{a}%
(\vartheta _{0}^{-}\sqcup \vartheta _{0}^{0})+\dot{B}(\text{$\mathbf{x}$}%
_{0}^{0},\pmb\vartheta )\dot{a}(x\sqcup \vartheta _{0}^{-}\sqcup \vartheta
_{0}^{0})](\vartheta _{-}^{0})
\end{equation*}%
\begin{equation*}
=\int_{X^{t}}\mathrm{d}x[D_{+}^{-}(x)a+D_{0}^{-}(x)\dot{a}(x)](\varkappa
)+\sum_{x\in \varkappa ^{t}}[D_{+}^{0}(x)a+D_{0}^{0}(x)\dot{a}(x)](\varkappa
/x)\ .
\end{equation*}%
Hence 
\begin{equation*}
U^{t}-U^{0}=\Lambda _{-}^{+}(t,D_{+}^{-})+\Lambda
_{-}^{0}(t,D_{0}^{-})+\Lambda _{0}^{+}(t,D_{+}^{0})+\Lambda
_{0}^{0}(t,D_{0}^{0})\ ,
\end{equation*}%
where $\Lambda _{\mu }^{\nu }(t)$ are the QS integrals (1.2) of operators 
\begin{equation*}
\lbrack D_{+}^{\mu }(x)a](\varkappa )=\sum_{\vartheta _{0}^{0}\sqcup
\vartheta _{+}^{0}\subseteq \varkappa }^{t(\vartheta _{\nu }^{0})<t(x)}\int_{%
\mathcal{X}^{t(x)}}\mathrm{d}\vartheta _{+}^{-}\int_{\mathcal{X}^{t(x)}}%
\mathrm{d}\vartheta _{0}^{-}[\dot{B}(\text{$\mathbf{x}$}_{+}^{\mu },\pmb%
\vartheta )\dot{a}(\vartheta _{0}^{-}\sqcup \vartheta _{0}^{0})](\vartheta
_{-}^{0}),
\end{equation*}%
\begin{equation*}
\lbrack D_{0}^{\mu }(x)b](\varkappa )=\sum_{\vartheta _{0}^{0}\sqcup
\vartheta _{+}^{0}\subseteq \varkappa }^{t(\vartheta _{\nu }^{0})<t(x)}\int_{%
\mathcal{X}^{t(x)}}\mathrm{d}\vartheta _{+}^{-}\int_{\mathcal{X}^{t(x)}}%
\mathrm{d}\vartheta _{0}^{-}[\dot{B}(\text{$\mathbf{x}$}_{0}^{\mu },\pmb%
\vartheta )\dot{b}(\vartheta _{0}^{-}\sqcup \vartheta _{0}^{0})](\vartheta
_{-}^{0})\ ,
\end{equation*}%
for $a\in \mathcal{G}^{+}\ ,\ b\in \mathcal{G}^{+}\otimes \mathcal{E}(x)$.
This can be written in terms of (1.5) as $D_{\nu }^{\mu }(x)=\imath
_{\lbrack 0,t(x))}^{\otimes }(\dot{B}(\mathbf{x}_{\nu }^{\mu })).$ Due to
the inequality%
\begin{equation*}
\Vert U^{t}\Vert _{\xi ^{+}}^{\xi _{-}}\leq \Vert B\Vert _{\eta ^{\bullet
}}^{\eta _{\bullet }}(t),\;\;\xi ^{+}\geq \sum \eta ^{\mu },\;\;\;\xi
_{-}^{-1}\geq \sum \eta _{\nu }^{-1}
\end{equation*}%
one obtains $\Vert D_{+}^{-}\Vert _{\xi ^{+},1}^{\xi _{-},t}\leq \Vert
B\Vert _{\eta ^{\bullet }}^{\eta _{\bullet }}(t)$: 
\begin{equation*}
\int_{X^{t}}\Vert D_{+}^{-}(x)\Vert _{\xi ^{+}}^{\xi _{-}}\mathrm{d}x\leq
\int_{X^{t}}\Vert \dot{B}_{+}^{-}(x)\Vert _{\eta ^{\bullet }}^{\eta
_{\bullet }}[t(x)]\mathrm{d}x=\int_{X^{t}}\mathrm{d}x\int_{\mathcal{X}%
^{t(x)}}\Vert B_{+}^{-}(x\sqcup \vartheta )\Vert _{\eta ^{\bullet }}^{\eta
_{\bullet }}[t(x)]\mathrm{d}\pmb\vartheta 
\end{equation*}%
\begin{equation*}
=\int_{\mathcal{X}^{t}}\Vert B_{+}^{-}({\mathcal{X}})\Vert _{\eta ^{\bullet
}}^{\eta _{\bullet }}(t)\mathrm{d}\varkappa -\Vert B_{+}^{-}(\pmb\emptyset
)\Vert _{\eta ^{\bullet }}^{\eta _{\bullet }}(t)=\Vert B\Vert _{\eta
^{\bullet }}^{\eta _{\bullet }}(t)-\Vert B_{+}^{-}(\pmb\emptyset )\Vert
_{\eta ^{\bullet }}^{\eta _{\bullet }}(t),
\end{equation*}%
where $B_{+}^{-}(\varkappa ,\pmb\vartheta )=B(\pmb\varkappa )1_{\emptyset
}(\vartheta _{+}^{-}),\,\pmb\varkappa =%
\begin{pmatrix}
\vartheta _{0}^{-} & \varkappa \cr\vartheta _{0}^{0} & \vartheta _{+}^{0}%
\end{pmatrix}%
$. In the same way we obtain 
\begin{equation*}
\int_{X^{t}}\left( \Vert D_{0}^{-}(x)\Vert _{\xi ^{+}}^{\xi _{-}}\right) ^{2}%
\mathrm{d}x\leq \int_{X^{t}}\left( \Vert \dot{B}_{0}^{-}(x)\Vert _{\eta
^{\bullet }}^{\eta _{\bullet }}[t(x)]\right) ^{2}\mathrm{d}x\leq \eta
^{-}(\Vert B\Vert _{\eta ^{\bullet }}^{\eta _{\bullet }}(t))^{2}\ ,
\end{equation*}%
\begin{equation*}
\int_{X^{t}}\left( \Vert D_{+}^{0}(x)\Vert _{\xi ^{+}}^{\xi _{-}}\right) ^{2}%
\mathrm{d}x\leq \int_{X^{t}}\left( \Vert \dot{B}_{+}^{0}(x)\Vert _{\eta
^{\bullet }}^{\eta _{\bullet }}[t(x)]\right) ^{2}\mathrm{d}x\leq \eta
_{+}^{-1}(\Vert B\Vert _{\eta ^{\bullet }}^{\eta _{\bullet }}(t))^{2}\ ,
\end{equation*}%
and 
\begin{equation*}
\mathrm{ess}\sup {}_{x\in X^{t}}\Vert D_{0}^{0}(x)\Vert _{\xi ^{+}}^{\xi
_{-}}\leq \mathrm{ess}\sup_{x\in X^{t}}\Vert \dot{B}_{0}^{0}(x)\Vert _{\eta
^{\bullet }}^{\eta _{\bullet }}[t(x)]\leq \sqrt{\eta ^{0}/\eta _{0}}\ \Vert
B\Vert _{\eta ^{\bullet }}^{\eta _{\bullet }}(t)\ .
\end{equation*}%
This proves the QS--integrability (1.3) of the derivatives $D_{\nu }^{\mu
}(x)$ with respect to the $(\xi ^{+},\xi _{-})$ norms.

If $B(\pmb\vartheta )$ satisfies the recurrence (1.8), then $\dot{B}(\mathbf{%
x}_{\nu }^{\mu })=L_{\nu }^{\mu }(x)\odot B$, and $D_{\nu }^{\mu }(x)=L_{\nu
}^{\mu }(x)\odot U^{t(x)}$ due to the property $\imath _{\lbrack
0,t)}^{\otimes }(L\odot B)=L\odot \imath _{\lbrack 0,t)}^{\otimes }(B)$ for $%
L\odot B=(L\otimes \hat{1})\cdot B$, following immediately from the
definition (1.5). Hence $U^{t}=\imath _{\lbrack 0,t)}^{\otimes
}(B)=U^{0}+\imath ^{t}(\mathbf{D})$ with $U^{0}=T^{0}\otimes \hat{1}$ and $B(%
\pmb\vartheta )$, defined for the partitions $\varkappa =\sqcup \vartheta
_{\nu }^{\mu }$ of the chains $\varkappa \in \mathcal{X}$ as $\mathbf{L}%
^{\triangleleft }(\pmb\vartheta )\odot T^{0}$, where $\mathbf{L}%
^{\triangleleft }(\pmb\vartheta )=L(\mathbf{x}_{n})\cdot \cdot \cdot L(%
\mathbf{x}_{1})$ for $\pmb\vartheta =\sqcup _{i=1}^{n}\mathbf{x}_{i}$,
satisfies the equation (1.10).\hfill 
\end{proof}

\begin{corollary}
Let $B(\pmb\vartheta )=L(\pmb\vartheta )\otimes \hat{1}$ be defined by the
QS--integrable operator--valued function 
\begin{equation*}
L\left( 
\begin{matrix}
\vartheta _{0}^{-} & \vartheta _{+}^{-}\cr\vartheta _{0}^{0} & \vartheta
_{+}^{0}%
\end{matrix}%
\right) :\mathcal{H}\otimes \mathcal{E}^{\otimes }(\vartheta
_{0}^{-})\otimes \mathcal{E}^{\otimes }(\vartheta _{0}^{0})\rightarrow 
\mathcal{H}\otimes \mathcal{E}^{\otimes }(\vartheta _{0}^{0})\otimes 
\mathcal{E}^{\otimes }(\vartheta _{+}^{0})
\end{equation*}%
with $\Vert L\Vert _{\eta ^{-},\eta ^{0}}^{\eta _{0},\eta _{+}}=\Vert B\Vert
_{\eta ^{\cdot }}^{\eta _{\cdot }}\ <\infty $ for $\eta ^{+},\eta
_{-}^{-1}\geq 1$. Then the QS integral $\imath _{\lbrack 0,t)}^{\otimes
}(B)=U^{t}$ defines an adapted $(\xi ^{+},\xi _{-})$ continuous process $%
U^{t}$ for $\xi ^{+}\geq \eta ^{-}+\eta ^{0}+1\ ,\ \xi _{-}^{-1}\geq \eta
_{+}^{-1}+\eta _{0}^{-1}+1$ in the sense $U^{t}(a^{t}\otimes c)=b^{t}\otimes
c$ for all $a^{t}\in \mathcal{G}^{t}(\xi ^{+}),c\in \mathcal{F}_{[t}$ with $%
b^{t}\in \mathcal{G}^{t}(\xi _{-})$, where $\mathcal{G}^{t}(\xi )=\mathcal{H}%
\otimes \mathcal{F}^{t}(\xi ),$ 
\begin{equation*}
\mathcal{F}^{t}(\xi )=\int_{\mathcal{X}^{t}}^{\oplus }\xi ^{|\varkappa |}%
\mathcal{E}_{\xi }^{\otimes }(\varkappa )\mathrm{d}\varkappa ,\;\;\mathcal{X}%
^{t}=\{\varkappa \in \mathcal{X}|t(\varkappa )\subset \lbrack 0,t)\},
\end{equation*}%
It has the adapted QS derivatives $D_{\nu }^{\mu }(x)=A_{\nu }^{\mu
}(x)\otimes \hat{1}_{[t(x)},$ where $\hat{1}_{[t}$ is the identity operator
in $\mathcal{F}_{[t(x)}$ with $A_{\nu }^{\mu }(x)$, defined in $\mathcal{G}%
^{t(x)}$. If $U^{t}$ is an adapted QS process with $\Vert U\Vert _{\xi
,\infty }^{\xi _{+},t}<\infty $, then $\mathrm{d}\imath ^{t}(\mathbf{B}\odot
U)=\mathrm{d}\imath ^{t}(\mathbf{B})U^{t}$ in the sense%
\begin{equation*}
\imath ^{t}(\mathbf{B}\odot U)a=\int_{0}^{t}\mathrm{d}\imath ^{s}(\mathbf{B}%
)U^{s}a,
\end{equation*}%
where $\mathbf{B}\odot U=\mathbf{B}\cdot (U\otimes \mathbf{1})$ and the left
hand side is defined as the limit of the It\^{o} integral sums $\imath ^{t}(%
\mathbf{B}\odot U)=\lim_{n\rightarrow \infty }\sum_{i=0}^{n}\left[ \imath
^{t_{i+1}}(\mathbf{B})-\imath ^{t_{i}}(\mathbf{B})\right] U^{t_{i}}\ ,\
t_{0}=0,t_{n+1}=t$ in the uniform $(\xi ,\xi _{-})$--topology.
\end{corollary}

Indeed, the QS integral (1.5) for $B=L\otimes \hat{1}$ and $a=a^{t}\otimes c$
with $a^{t}\in \mathcal{G}^{t}\ ,\ c\in \mathcal{F}_{[t}$ can be written as 
\begin{equation*}
\lbrack \imath _{\lbrack 0,t)}^{\otimes }(B)a](\varkappa )=c(\varkappa
_{\lbrack t})\otimes \sum_{\sqcup \vartheta _{\nu }^{0}=\varkappa ^{t}}\int_{%
\mathcal{X}^{t}}\int_{\mathcal{X}^{t}}[L(\pmb\vartheta )\otimes I(\vartheta
_{-}^{0})]a(\vartheta _{0}^{-}\sqcup \vartheta _{0}^{0}\sqcup \vartheta
_{-}^{0})\mathrm{d}\vartheta _{0}^{-}\mathrm{d}\vartheta _{+}^{-}\ .
\end{equation*}%
The norm $\Vert L\otimes \hat{1}\Vert _{\eta ^{\bullet }}^{\eta _{\bullet }}$
for $\eta ^{+},\eta _{-}^{-1}\geq 1$ does not depend on $\eta ^{+},\eta
_{-}^{-1}$, hence%
\begin{equation*}
\Vert \imath _{\lbrack 0,t)}^{\otimes }(L\otimes \hat{1})\Vert _{\xi
^{+}}^{\xi _{-}}\leq \Vert L\Vert _{\eta ^{-},\eta ^{0}}^{\eta _{0},\eta
_{+}},
\end{equation*}%
if $\xi ^{+}\geq \sum \eta ^{\mu }\ ,\ \xi _{-}^{-1}\geq \sum \eta _{\nu
}^{-1}$ with $\eta ^{+}=1=\eta _{-}$.

The derivatives $D_{\nu }^{\mu }$ are adapted as multiple QS integrals $%
\imath _{\lbrack o,t(x))}^{\otimes }(\dot{B}(\mathbf{x}_{\nu }^{\mu }))$ of $%
\dot{B}(\mathbf{x})=\dot{L}(\mathbf{x})\otimes \hat{1}$. If $U^{s}\ ,\ s<t$
is a simple adapted function $U^{s}=\sum_{i=0}^{n}U_{i}1_{[t_{i},t_{i+1})}(s)
$ with $t_{0}=0,t_{n+1}=t\ ,1_{[t,t_{+})}(s)=1$ for $s\in \lbrack t,t_{+})$,
otherwise $1_{[t,t_{+})}(s)=0$, then 
\begin{equation*}
\imath ^{t}(\mathbf{B}\odot U)a=\sum_{i=0}^{n}(\imath ^{t_{i+1}}-\imath
^{t_{i}})(\text{$\mathbf{B}$}\odot U_{i})a=\sum_{i=0}^{n}[\imath ^{t_{i+1}}(%
\text{$\mathbf{B}$})-\imath ^{t_{i}}(\text{$\mathbf{B}$})]b_{i}\ ,
\end{equation*}%
where $b_{i}=U_{i}a$, if $U$ is a constant adapted process on $[r,s)$: 
\begin{equation*}
\lbrack \imath _{\lbrack r,s)}^{\otimes }(BU)a](\varkappa )=\sum_{\vartheta
_{+}^{0}\sqcup \vartheta _{0}^{0}\subseteq \varkappa _{r}^{s}}\int_{\mathcal{%
X}_{r}^{s}}\int_{\mathcal{X}_{r}^{s}}[B(\pmb\vartheta )\dot{b}(\vartheta
_{0}^{-}\sqcup \vartheta _{0}^{0})](\vartheta _{-}^{0})\mathrm{d}\vartheta
_{0}^{-}\mathrm{d}\vartheta _{+}^{-}\ .
\end{equation*}

\section{A nonadapted QS calculus and It\^{o} formula}

Now we shall consider the operators $U=\epsilon (T)$ acting in $\mathcal{G}=%
\mathcal{H}\otimes \mathcal{F}$ as the multiple QS integrals (1.5) with $%
B=L\otimes \hat{1}$, and $t=\infty $ according to the formula 
\begin{equation*}
\lbrack \epsilon (T)a](\varkappa )=\sum_{\varkappa _{0}^{0}\sqcup \varkappa
_{+}^{0}=\varkappa }\iint T(\pmb\varkappa )a(\varkappa _{0}^{0}\sqcup
\varkappa _{0}^{-})\mathrm{d}\varkappa _{0}^{-}\mathrm{d}\varkappa _{+}^{-}%
\eqno(2.1)
\end{equation*}%
Here the sum is taken over all partitions of the chain $\varkappa \in 
\mathcal{X}$, and the operator-valued function $T_{(}\pmb\varkappa )$ is in
one to one correspondence 
\begin{align*}
& T\left( 
\begin{matrix}
\varkappa _{0}^{-} & \varkappa _{+}^{-}\cr\varkappa _{0}^{0} & \varkappa
_{+}^{0}%
\end{matrix}%
\right) =\sum_{\vartheta \subseteq \varkappa _{0}^{0}}L\left( 
\begin{matrix}
\varkappa _{0}^{-} & \varkappa _{+}^{-}\cr\vartheta  & \varkappa _{+}^{0}%
\end{matrix}%
\right) \otimes I^{\otimes }(\varkappa _{0}^{0}\backslash \vartheta ) \\
& L\left( 
\begin{matrix}
\vartheta ^{-} & \vartheta _{+}^{-} \\ 
\vartheta  & \vartheta _{+}%
\end{matrix}%
\right) =\sum_{\varkappa \subseteq \vartheta }(-1)^{|\varkappa |}T\left( 
\begin{matrix}
\vartheta ^{-} & \vartheta _{+}^{-} \\ 
\vartheta \backslash \varkappa  & \vartheta _{+}%
\end{matrix}%
\right) \otimes I^{\otimes }(\varkappa ) \\
&
\end{align*}%
with the operator-valued function $L(\pmb\vartheta )$, defining the integral
representation $U=\Lambda _{\lbrack 0,\infty )}(L\otimes \hat{1})$.

Using the arguments in section 1., one can prove, that the operator $%
\epsilon (T)$ is $(\xi ^{+},\xi _{-})$-continuous, if $T$ is $(\zeta
^{\bullet },\zeta _{\bullet })$-bounded for $\zeta ^{\bullet }=(\zeta
^{-},\zeta ^{0})$ and $\zeta _{\bullet }=(\zeta _{0},\zeta _{+})$,
satisfying the inequalities $\zeta ^{-}+\zeta ^{0}\leq \xi ^{+}$, $\zeta
_{0}^{-1}+\zeta _{+}^{-1}\leq \xi _{-}^{-1}$, because 
\begin{equation*}
\Vert \epsilon (T)\Vert _{\xi ^{+}}^{\xi _{-}}\leq \Vert T\Vert _{\zeta
^{\bullet }}^{\zeta _{\bullet }}\equiv \int \left( \iint \frac{(\zeta
_{+})^{|\varkappa _{+}^{0}|}}{(\zeta ^{-})^{|\varkappa _{0}^{-}|}}\mathrm{es}%
\text{$\mathrm{s\sup }$}_{\varkappa _{0}^{0}\in \mathcal{X}}{\frac{(\zeta
_{0})^{|\varkappa _{0}^{0}|}}{(\zeta ^{0})^{|\varkappa _{0}^{0}|}}}\Vert T(%
\pmb\varkappa )\Vert ^{2}\mathrm{d}\varkappa _{+}^{0}\mathrm{d}\varkappa
_{0}^{-}\right) ^{1/2}\mathrm{d}\varkappa _{+}^{-}\ .
\end{equation*}%
In this case the formally conjugated operator 
\begin{equation*}
U^{\ast }=\epsilon (T^{{\star }})\ ,\ T^{{\star }}(\pmb\varkappa )=T(\pmb%
\varkappa ^{{\star }})^{\ast }\ ,\ \left( 
\begin{matrix}
\varkappa _{0}^{-} & \varkappa _{+}^{-}\cr\varkappa _{0}^{0} & \varkappa
_{+}^{0}%
\end{matrix}%
\right) ^{{\star }}=\left( 
\begin{matrix}
\varkappa _{+}^{0} & \varkappa _{+}^{-}\cr\varkappa _{0}^{0} & \varkappa
_{0}^{-}%
\end{matrix}%
\right) \eqno(2.2)
\end{equation*}%
exists as $(\xi ^{+},\xi _{-})$-continuous operator $\mathcal{G}(\xi
_{+})\rightarrow \mathcal{G}(\xi ^{-})$ with $\Vert U^{\ast }\Vert _{\xi
^{+}}^{\xi _{-}}=\Vert U\Vert _{1/\xi _{-}}^{1/\xi ^{+}}$, if $\xi ^{+}\geq
\zeta _{0}^{-1}+\zeta _{+}^{-1}$, $\xi _{-}^{-1}\geq \zeta ^{-}+\zeta ^{0}$.

As we shall prove now, the map $\epsilon $ is the representation in $%
\mathcal{G}$ of a unital $\star $-algebra of operator-valued functions $T(%
\pmb\varkappa )$, satisfying the relative boundedness condition

\begin{equation*}
\Vert T\Vert (\pmb\zeta )=\text{$\mathrm{ess\sup }$}_{\pmb\varkappa }\{\Vert
T(\pmb\varkappa )\Vert /\dprod_{\mu \leq \nu }\zeta _{\nu }^{\mu }(\varkappa
_{\nu }^{\mu })\}<\infty \ ,\eqno(2.3)
\end{equation*}%
where $\zeta (\varkappa )=\dprod_{x\in \varkappa }\zeta (x)$, with respect
to a triangular matrix-function $\pmb\zeta (x)=[\zeta _{\nu }^{\mu }(x)]$, $%
\mu ,\nu =-,0,+$ $\zeta _{\nu }^{\mu }=0$ for $\mu >\nu $ under the order $%
-<0<+\ $,$\ \zeta _{-}^{-}(x)=1=\zeta _{+}^{+}(x)$ with positive $L^{p}$%
-integrable functions $(\zeta _{\nu }^{\mu })_{\nu =0,+}^{\mu =-,0}$ for
corresponding $p=1,2,\infty $: 
\begin{equation*}
\Vert \zeta _{+}^{-}\Vert _{1}\leq \infty ,\Vert \zeta _{0}^{-}\Vert
_{2}<\infty \ ,\ \Vert \zeta _{+}^{0}\Vert _{2}<\infty \ ,\ \Vert \zeta
_{0}^{0}\Vert _{\infty }<\infty \ ,
\end{equation*}%
where $\Vert \zeta \Vert _{p}=(\int \zeta ^{p}(x)\mathrm{d}x)^{1/p}$. In
this case the operator $U=\epsilon (T)$ is $\zeta $-bounded, as it follows
from the next theorem, for $\zeta >\Vert \zeta _{0}^{0}\Vert =\mathrm{%
ess\sup }_{x\in X}\zeta _{0}^{0}(x)$ in the sense of $(\xi ^{+},\xi _{-})$%
-continuity of $U$ for all $\xi _{-}>0\ ,\ \xi ^{+}\geq \zeta \cdot \xi _{-}$%
. This is due to the estimate 
\begin{equation*}
\Vert T\Vert _{\zeta ^{\bullet }}^{\zeta _{\bullet }}\leq \int \left( \iint 
\frac{(\zeta _{+})^{|\varkappa _{+}^{0}|}}{(\zeta ^{-})^{|\varkappa
_{0}^{-}|}}\text{$\mathrm{ess\sup }$}_{\varkappa _{0}^{0}}{\frac{(\zeta
_{0})^{|\varkappa _{0}^{0}|}}{(\zeta ^{0})^{|\varkappa _{0}^{0}|}}}\left[
\prod_{\mu <\nu }\zeta _{\nu }^{\mu }(\varkappa _{\nu }^{\mu })\right] ^{2}%
\mathrm{d}\varkappa _{+}^{0}\mathrm{d}\varkappa _{0}^{-}\right) ^{1/2}%
\mathrm{d}\varkappa _{+}^{-}\Vert T\Vert (\pmb\zeta )
\end{equation*}%
\begin{equation*}
=\int \dprod_{x\in \varkappa }\zeta _{+}^{-}(x)\mathrm{d}\varkappa \left(
\int \dprod_{x\in \varkappa }{\frac{\zeta _{0}^{-}(x)^{2}}{\zeta ^{-}}}%
\mathrm{d}\varkappa \int \dprod_{x\in \varkappa }{\frac{\zeta _{+}^{0}(x)^{2}%
}{\zeta _{+}^{-1}}}\mathrm{d}\varkappa \text{$\mathrm{ess\sup }$}%
\dprod_{x\in \varkappa }{\frac{\zeta _{0}^{0}(x)^{2}}{\zeta ^{0}\zeta
_{0}^{-1}}}\right) ^{1/2}\Vert T\Vert (\pmb\zeta )
\end{equation*}%
\begin{equation*}
\leq \exp \{\int (\zeta _{+}^{-}(x)+(\zeta _{0}^{-}(x)^{2}+\zeta
_{+}^{0}(x)^{2})/2\varepsilon )\mathrm{d}x\}\Vert T\Vert (\pmb\zeta )\ ,\eqno%
(2.4)
\end{equation*}%
for $\zeta ^{-},\zeta _{+}^{-1}\geq \varepsilon >0$, and $\zeta ^{0}\zeta
_{0}^{-1}\geq \Vert \zeta _{0}^{0}\Vert _{\infty }^{2}$, giving $\epsilon
(T)=0$, if $T(\pmb\varkappa )=0$ for almost all $\pmb\varkappa $. Hence the
operator (2.1) is defined even if $T(\pmb\varkappa )$ is described for
almost all $\pmb\varkappa =(\varkappa _{\nu }^{\mu })$, in particular, only
for the partitions $\varkappa =\sqcup \varkappa _{\nu }^{\mu }$ of the
chains $\varkappa \in \mathcal{X}$.

\begin{theorem}
If the operator-valued function 
\begin{equation*}
T\left( 
\begin{matrix}
\varkappa _{0}^{-} & \varkappa _{+}^{-}\cr\varkappa _{0}^{0} & \varkappa
_{+}^{0}%
\end{matrix}%
\right) :\mathcal{H}\otimes \mathcal{E}^{\otimes }(\varkappa
_{0}^{-})\otimes \mathcal{E}^{\otimes }(\varkappa _{0}^{0})\rightarrow 
\mathcal{H}\otimes \mathcal{E}^{\otimes }(\varkappa _{0}^{0})\otimes 
\mathcal{E}^{\otimes }(\varkappa _{+}^{0})
\end{equation*}%
satisfies the condition (2.3), then the conjugated operators $U=\epsilon
(T),U^{\ast }=\epsilon (T^{{\star }})$ are $\zeta $-bounded in $\mathcal{G}$
for any $\zeta >\zeta _{0}^{0}$, and the operator $U^{\ast }U$ is defined in 
$\mathcal{G}$ as $\zeta ^{2}$-bounded operator 
\begin{equation*}
\epsilon (S\cdot T)=\epsilon (S)\epsilon (T)\;,\quad S=T^{{\star }}
\end{equation*}%
by the following product formula 
\begin{equation*}
(S\cdot T)(\pmb\varkappa )=\sum_{\vartheta _{\nu }^{\mu }\subseteq \varkappa
_{\nu }^{\mu }}^{\mu <\nu }\sum_{\sigma _{+}^{-}\cap \rho _{+}^{-}=\vartheta
_{+}^{-}}^{\sigma _{+}^{-}\cup \rho _{+}^{-}=\varkappa _{+}^{-}}S\left( 
\begin{matrix}
\vartheta _{0}^{-}\sqcup \vartheta _{+}^{-}, & \varkappa _{+}^{-}\backslash
\sigma _{+}^{-}\cr\varkappa _{0}^{0}\sqcup \vartheta _{+}^{0}, & \varkappa
_{+}^{0}\backslash \vartheta _{+}^{0}%
\end{matrix}%
\right) T\left( 
\begin{matrix}
\varkappa _{0}^{-}\backslash \vartheta _{0}^{-}, & \varkappa
_{+}^{-}\backslash \rho _{+}^{-}\cr\varkappa _{0}^{0}\sqcup \vartheta
_{0}^{-}, & \vartheta _{+}^{-}\sqcup \vartheta _{+}^{0}%
\end{matrix}%
\right) \ .\eqno(2.5)
\end{equation*}%
This induces a unital $\star $-algebraic structure on the inductive space $%
\mathcal{U}$ of all relatively bounded functions $T$ with 
\begin{equation*}
\Vert T^{{\star }}\Vert (\pmb\zeta )=\Vert T\Vert (\pmb\zeta ^{{\star }%
}),\;\;\;\Vert T^{{\star }}\cdot T\Vert (\pmb\xi )\leq \lbrack \Vert T\Vert (%
\pmb\zeta )]^{2},
\end{equation*}%
if $\xi _{\nu }^{\mu }\geq (\pmb\zeta ^{{\star }}\pmb\zeta )_{\nu }^{\mu }$,
where $\pmb\zeta ^{{\star }}(x)=\mathbf{g}\pmb\zeta (x)^{\ast }\mathbf{g}$
and $(\pmb\zeta ^{{\star }}\pmb\zeta )(x)=\pmb\zeta ^{{\star }}(x)\pmb\zeta
(x)$ are defined by usual product of the matrices 
\begin{equation*}
\mathbf{g}=\left[ 
\begin{matrix}
0 & 0 & 1\cr0 & 1 & 0\cr1 & 0 & 0%
\end{matrix}%
\right] ,\ \pmb\zeta (x)=\left[ 
\begin{matrix}
1 & \zeta _{0}^{-} & \zeta _{+}^{-}\cr0 & \zeta _{0}^{0} & \zeta _{+}^{0}\cr0
& 0 & 1%
\end{matrix}%
\right] (x),\ \pmb\zeta ^{{\star }}(x)=\left[ 
\begin{matrix}
1 & \zeta _{+}^{0} & \zeta _{+}^{-}\cr0 & \zeta _{0}^{0} & \zeta _{0}^{-}\cr0
& 0 & 1%
\end{matrix}%
\right] (x)\ .\eqno(2.6)
\end{equation*}%
If the multiple QS integral $U^{t}=\Lambda _{\lbrack 0,t)}(B)$ is defined by 
$B(\pmb\vartheta )=\epsilon (L(\pmb\vartheta ))$ with 
\begin{equation*}
\left\Vert L\left( 
\begin{matrix}
\vartheta ^{-} & \vartheta _{+}^{-}\cr\vartheta  & \vartheta _{+}%
\end{matrix}%
\right) \right\Vert (\pmb\xi )\leq c\lambda _{0}^{0}(\vartheta )\lambda
_{+}^{0}(\vartheta _{+})\lambda _{0}^{-}(\vartheta ^{-})\lambda
_{+}^{-}(\vartheta _{+}),\;\lambda (\vartheta )=\dprod_{x\in \vartheta
}\lambda (x)\geq 0\ ,
\end{equation*}%
then $\Lambda _{\lbrack 0,t)}\circ \epsilon =\epsilon \circ N_{[0,t)}$ i.e. $%
U^{t}=\epsilon (T^{t})$, where 
\begin{equation*}
T^{t}(\pmb\varkappa )=\sum_{\pmb\vartheta \subseteq {\pmb\varkappa }^{t}}L(%
\pmb\vartheta ,\pmb\varkappa \backslash {\pmb\vartheta })\equiv N_{[0,t)}(L)(%
\pmb\varkappa )\eqno(2.7)
\end{equation*}%
with $\Vert T^{t}\Vert (\pmb\zeta )\leq c$, if $\zeta _{\nu }^{\mu }(x)\geq
\xi _{\nu }^{\mu }(x)+\lambda _{\nu }^{\mu }(x)$ for $t(x)<t$, and $\zeta
_{\nu }^{\mu }(x)\geq \xi _{\nu }^{\mu }(x)$ for $t(x)\geq t$. The QS
derivatives $D_{\nu }^{\mu }(x)=\Lambda _{\lbrack 0,t(x))}(\dot{B}(\mathbf{x}%
_{\nu }^{\mu }))$ for the process $U^{t}=\epsilon (T^{t})$ have the natural
difference form $\mathbf{D}=\mathbf{G}-\mathbf{U}$, described by the
representations of 
\begin{equation*}
\dot{T}^{t}(\mathbf{x},\pmb\varkappa )=T^{t}(\pmb\varkappa \sqcup \mathbf{x})
\end{equation*}%
with $\mathbf{x}=\mathbf{x}_{\nu }^{\mu },\mu <+,\nu >-$ at $t\searrow t(x)$ 
\begin{equation*}
U_{\nu }^{\mu }(x)=\epsilon (\dot{T}^{t(x)}(\mathbf{x}_{\nu }^{\mu })),\quad
G_{\nu }^{\mu }(x)=\epsilon (\dot{T}^{t(x)]}(\mathbf{x}_{\nu }^{\mu }))\ ,%
\eqno(2.8)
\end{equation*}%
where $T^{s]}(\pmb\varkappa )=N_{[0,t)}(L)(\pmb\varkappa )$ for any $t>s\ ,\
t\leq t(x)$\ $\forall x\in \varkappa _{s}=\sqcup \varkappa _{\nu }^{\mu
}\cap t^{-1}(s,\infty )$, and $\mathbf{x}_{\nu }^{\mu }$ denotes an
elementary table $\pmb\vartheta =(\vartheta _{\lambda }^{\kappa })$ with $%
\vartheta _{\lambda }^{\kappa }=\emptyset $ except $\vartheta _{\nu }^{\mu
}=x$. The QS differential $\mathrm{d}U^{\ast }=\mathrm{d}\Lambda (\mathbf{D}%
^{{\star }})$ is defined by the derivative $\mathbf{D}^{{\star }}=\mathbf{G}%
^{{\star }}-\mathbf{U}^{{\star }}$, and 
\begin{equation*}
\mathrm{d}(U^{\ast }U)=\mathrm{d}\Lambda (\mathbf{U}^{{\star }}\mathbf{D}+%
\mathbf{D}^{{\star }}\mathbf{U}+\mathbf{D}^{{\star }}\mathbf{D})=\mathrm{d}%
\Lambda (\mathbf{G}^{{\star }}\mathbf{G}-\mathbf{U}^{{\star }}\mathbf{U})\ ,%
\eqno(2.9)
\end{equation*}%
where the QS derivative $\mathbf{G}^{{\star }}\mathbf{G}-\mathbf{U}^{{\star }%
}\mathbf{U}$ of the QS process $(U^{\ast }U)^{t}=U^{t{\ast }}U^{t}$ is
described in terms of the usual products $(\mathbf{U}^{{\star }}\mathbf{U}%
)(x)=\mathbf{U}^{{\star }}(x)\mathbf{U}(x)$ and the pseudo Hermitian
conjugation $\mathbf{U}^{{\star }}(x)=(\hat{I}\otimes \mathbf{g}(x))\mathbf{U%
}(x)^{\ast }(\hat{I}\otimes \mathbf{g}(x))$ of the triangular matrices 
\begin{equation*}
\mathbf{U}=\left[ 
\begin{matrix}
U & U_{0}^{-} & U_{+}^{-}\cr0 & U_{0}^{0} & U_{+}^{0}\cr0 & 0 & U%
\end{matrix}%
\right] ,\ \mathbf{D}=\left[ 
\begin{matrix}
0 & D_{0}^{-} & D_{+}^{-}\cr0 & D_{0}^{0} & D_{+}^{0}\cr0 & 0 & 0%
\end{matrix}%
\right] ,\ \mathbf{G}=\left[ 
\begin{matrix}
U, & G_{0}^{-}, & G_{+}^{-}\cr0 & G_{0}^{0} & G_{+}^{0}\cr0 & 0 & U%
\end{matrix}%
\right] ,
\end{equation*}%
with $U(x)=U^{t(x)},\mathbf{g}(x)=[g_{\nu }^{\mu
}(x)],g_{-}^{-}=1=g_{+}^{+},g_{0}^{0}(x)=I(x)$, otherwise $g_{\nu }^{\mu }=0$%
.
\end{theorem}

\begin{proof}
Let us firstly obtain an estimate for the representation $U=\epsilon (T)$ of
a relatively bounded operator-valued function $T$ in the sense (2.3). Due to
the inequalities (2.4) one obtains 
\begin{equation*}
\Vert U\Vert _{\xi ^{+}}^{\xi }\leq \exp \{\Vert \zeta _{+}^{-}\Vert
_{1}+(\Vert \zeta _{0}^{-}\Vert _{2}^{2}+\Vert \zeta _{+}^{0}\Vert
_{2}^{2})/2\varepsilon \}\Vert T\Vert (\pmb\zeta )\ ,
\end{equation*}%
if 
\begin{equation*}
\xi ^{+}\geq \varepsilon +\zeta ^{0},\quad \xi _{-}^{-1}\geq \zeta
_{0}^{-1}+\varepsilon \ ,\ \zeta ^{0}\zeta _{0}^{-1}\geq \Vert \zeta
_{0}^{0}\Vert _{\infty }^{2}\ .
\end{equation*}%
Hence for any $\xi ^{+}\xi _{-}^{-1}>\Vert \zeta _{0}^{0}\Vert _{\infty }^{2}
$ there exists an $\varepsilon >0$ such that this inequality holds, namely, 
\begin{equation*}
\varepsilon \leq \left( \xi ^{+}+\xi _{-}^{-1}-\sqrt{(\xi ^{+}-\xi
_{-}^{-1})^{2}+4\Vert \zeta _{0}^{0}\Vert _{\infty }^{2}}\right) \qquad
/2=\varepsilon (\xi ^{+},\xi _{-})\ ,
\end{equation*}%
where the upper bound $\varepsilon (\xi ^{+},\xi _{-})$ corresponds to the
solution $\varepsilon >0$ of the equation $\zeta ^{0}\zeta _{0}^{-1}=\Vert
\zeta _{0}^{0}\Vert _{\infty }^{2}$ with $\zeta ^{0}=\xi ^{+}-\varepsilon
>0,\;\zeta _{0}^{-1}=\xi _{-}^{-1}-\varepsilon >0$. Hence the operator $U$
is $\zeta $-bounded for any $\zeta >\zeta _{0}^{0}$ and also $U^{\ast }$ is $%
\zeta $-bounded due to $(\pmb\zeta ^{{\star }})_{0}^{0}=\zeta _{0}^{0}=(\pmb%
\zeta )_{0}^{0}$.

Now we show that the product formula (2.3) is valid for $T(\varkappa
)=X\otimes \mathbf{f}^{\otimes }(\pmb\varkappa )$, where $X\in \mathcal{B}(%
\mathcal{H})$, and 
\begin{equation*}
\mathbf{f}^{\otimes }(\pmb\varkappa )=\otimes _{\mu \leq \nu }f_{\nu }^{\mu
}(\varkappa _{\nu }^{\mu })
\end{equation*}%
with $f_{\nu }^{\mu }(\varkappa )=\otimes _{x\in \varkappa }f_{\nu }^{\mu
}(x)$ defined by the operator-valued elements $(f_{\nu }^{\mu })_{\nu
=0,+}^{\mu =-,0}$ of the matrix-function 
\begin{equation*}
\mathbf{f}(x)=\left[ 
\begin{matrix}
1 & f_{0}^{-} & f_{+}^{-}\cr0 & f_{0}^{0} & f_{+}^{0}\cr0 & 0 & 1%
\end{matrix}%
\right] (x),\quad 
\begin{array}{ll}
& f_{0}^{0}(x):\mathcal{E}(x)\rightarrow \mathcal{E}(x)\ ,\;f_{+}^{-}(x):%
\mathbb{C}\rightarrow \mathbb{C} \\ 
& f_{+}^{0}(x):\mathbb{C}\rightarrow \mathcal{E}(x)\ ,\;f_{0}^{-}(x):%
\mathcal{E}(x)\rightarrow \mathbb{C}%
\end{array}%
\end{equation*}%
with 
\begin{equation*}
\Vert f_{+}^{-}\Vert _{1}<\infty ,\quad \Vert f_{0}^{-}\Vert _{2}<\infty
,\quad \Vert f_{+}^{0}\Vert _{2}<\infty ,\quad \Vert f_{0}^{0}\Vert _{\infty
}<\infty \ .
\end{equation*}%
Let us find the action (2.1) of the operator $U=\epsilon (X\otimes \mathbf{f}%
^{\otimes })$ on the product vector $a=h\otimes k^{\otimes }\ ,\ h\in 
\mathcal{H}\ ,\ k\in \mathcal{K}$, where $k^{\otimes }(\varkappa )=\otimes
_{x\in \varkappa }k(x)$: 
\begin{align*}
& [Ua](\varkappa )=Xh\otimes \sum_{\varkappa _{0}^{0}\sqcup \varkappa
_{+}^{0}=\varkappa }\iint f_{+}^{-}(\varkappa _{+}^{-})f_{0}^{-}(\varkappa
_{0}^{-})k^{\otimes }(\varkappa _{0}^{-})f_{+}^{0}(\varkappa
_{+}^{0})\otimes f_{0}^{0}(\varkappa _{0}^{0})k^{\otimes }(\varkappa
_{0}^{0})\mathrm{d}\varkappa _{+}^{-}\mathrm{d}\varkappa _{0}^{-} \\
& =Xh\otimes \sum_{\varkappa _{0}^{0}\sqcup \varkappa _{+}^{0}=\varkappa
}\otimes _{x\in \varkappa _{+}^{0}}f_{+}^{0}(x)\otimes _{x\in \varkappa
_{0}^{0}}f_{0}^{0}(x)k(x)\int \dprod_{x\in \varkappa }f_{+}^{-}(x)\mathrm{d}%
\varkappa \int \dprod_{x\in \varkappa }f_{0}^{-}(x)k(x)\mathrm{d}\varkappa 
\\
& =Xh\otimes (f_{+}^{0}+f_{0}^{0}k)^{\otimes }(\varkappa )\exp \{\int
(f_{+}^{-}(x)+f_{0}^{-}(x)k(x))\mathrm{d}x\}
\end{align*}

In the same way, acting on the product vector $Xh\otimes
(f_{+}^{0}+f_{0}^{0}k)^{\otimes }$ by $U^{\ast }=\epsilon (X^{\ast }\otimes 
\mathbf{f}^{{\star }\otimes })$ with 
\begin{equation*}
\mathbf{f}^{{\star }}(x)_{0}^{0}=f_{0}^{0}(x)^{\ast },\quad \mathbf{f}^{{%
\star }}(x)_{+}^{-}=f_{+}^{-}(x)^{\ast },\quad \mathbf{f}^{{\star }%
}(x)_{+}^{0}=f_{0}^{-}(x)^{\ast },\quad \mathbf{f}^{{\star }%
}(x)_{0}^{-}=f_{+}^{0}(x)^{\ast }\ ,
\end{equation*}%
one obtains $[U^{\ast }Ua](\varkappa )=X^{\ast }Xh\otimes (f_{0}^{-{\ast }%
}+f_{0}^{0{\ast }}(f_{+}^{0}+f_{0}^{0}k))^{\otimes }(\varkappa )$. 
\begin{equation*}
\exp \{\int [f_{+}^{-}(x)^{\ast }+f_{+}^{0{\ast }%
}(x)(f_{+}^{0}(x)+f_{0}^{0}(x)k(x))+f_{+}^{-}(x)+f_{0}^{-}(x)k(x)]\mathrm{d}%
x\}=
\end{equation*}%
\begin{equation*}
=X^{\ast }Xh\otimes ((\mathbf{f}^{{\star }}\mathbf{f})_{+}^{0}+(\mathbf{f}^{{%
\star }}\mathbf{f})_{0}^{0}k)^{\otimes }(\varkappa )\exp \{\int ((\mathbf{f}%
^{{\star }}\mathbf{f})_{+}^{-}(x)+(\mathbf{f}^{{\star }}\mathbf{f}%
)_{0}^{-}k(x))\mathrm{d}x\}\ ,
\end{equation*}%
where the operator-valued functions 
\begin{equation*}
(\mathbf{f}^{{\star }}\mathbf{f})_{0}^{0}(x)=f_{0}^{0}(x)^{\ast
}f_{0}^{0}(x),(\mathbf{f}^{{\star }}\mathbf{\ f})_{+}^{-}(x)=f_{+}^{-}(x)^{%
\ast }+f_{+}^{0}(x)^{\ast }f_{+}^{0}(x)+f_{+}^{-}(x)
\end{equation*}%
\begin{equation*}
(\mathbf{f}^{{\star }}\mathbf{f})_{+}^{0}(x)=f_{0}^{-}(x)^{\ast
}+f_{0}^{0}(x)^{\ast }f_{+}^{0}(x),(\mathbf{f}^{\star }\mathbf{f}%
)_{0}^{-}(x)=f_{+}^{0}(x)^{\ast }f_{0}^{0}(x)+f_{0}^{-}(x)
\end{equation*}%
are defined as matrix elements of the product $(\mathbf{f}^{{\star }}\mathbf{%
f})(x)=\mathbf{f}(x)^{{\star }}\mathbf{f}(x)$ of triangular matrices $%
\mathbf{f}^{{\star }}$ and $\mathbf{f}$. Hence on the linear span of the
product vectors $a=h\otimes k^{\otimes }$ we have for $T=X\otimes \mathbf{f}%
^{{\star }}$ the $\star $-multiplicative property 
\begin{equation*}
\epsilon (T)^{\ast }\epsilon (T)=\epsilon (X^{\ast }X\otimes (\mathbf{f}^{{%
\star }}\mathbf{f})^{\otimes })=\epsilon (T^{{\star }}\cdot T)\ ,
\end{equation*}%
where the product $(T^{{\star }}\cdot T)(\pmb\varkappa )$ is defined as
(2.5) due to $(\mathbf{f}^{{\star }}\mathbf{f})^{\otimes }=\mathbf{f}^{{%
\otimes \star }}\cdot \mathbf{f}^{{\otimes }}$: 
\begin{equation*}
(\mathbf{f}^{{\star }}\mathbf{f})^{\otimes }(\pmb\varkappa )=\otimes _{x\in
\varkappa _{0}^{0}}(f_{0}^{0}(x)^{\ast }f_{0}^{0}(x))\otimes _{x\in
\varkappa _{+}^{0}}(f_{0}^{-}(x)^{\ast }+f_{0}^{0}(x)^{\ast
}f_{+}^{0}(x))\otimes 
\end{equation*}%
\begin{equation*}
\otimes _{x\in \varkappa _{0}^{-}}(f_{+}^{0}(x)^{\ast
}f_{0}^{0}(x)+f_{0}^{-}(x))\otimes _{x\in \varkappa
_{+}^{-}}(f_{+}^{-}(x)^{\ast }+f_{+}^{0}(x)^{\ast
}f_{+}^{0}(x)+f_{+}^{-}(x))=
\end{equation*}%
\begin{equation*}
=\sum_{\vartheta _{\nu }^{\mu }\subseteq \varkappa _{\nu }^{\mu }}^{\mu <\nu
}f_{0}^{0}(\varkappa _{0}^{0})^{\ast }f_{0}^{0}(\varkappa _{0}^{0})\otimes
f_{0}^{-}(\varkappa _{+}^{0}\backslash \vartheta _{+}^{0})^{\ast }\otimes
f_{0}^{0}(\vartheta _{+}^{0})^{\ast }f_{+}^{0}(\vartheta _{+}^{0})\otimes 
\end{equation*}%
\begin{equation*}
f_{+}^{0}(\vartheta _{0}^{-})^{\ast }f_{0}^{0}(\vartheta _{0}^{-}))\otimes
f_{0}^{-}(\varkappa _{0}^{-}\backslash \vartheta _{0}^{-})\sum_{\sigma
_{+}^{-}\cap \rho _{+}^{-}=\vartheta _{+}^{-}}^{\sigma _{+}^{-}\cup \rho
_{+}^{-}=\varkappa _{+}^{-}}f_{+}^{-}(\varkappa _{+}^{-}\backslash \sigma
_{+}^{-})^{\ast }f_{+}^{0}(\vartheta _{+}^{-})^{\ast }f_{+}^{0}(\vartheta
_{+}^{-})f_{+}^{-}(\varkappa _{+}^{-}\backslash \rho _{+}^{-})
\end{equation*}%
\begin{equation*}
=\sum_{\vartheta _{\nu }^{\mu }\subseteq \varkappa _{\nu }^{\mu }}^{\mu <\nu
}\sum_{\sigma _{+}^{-}\cap \rho _{+}^{-}=\vartheta _{+}^{-}}^{\sigma
_{+}^{-}\cup \rho _{+}^{-}=\varkappa _{+}^{-}}\mathbf{f}^{\otimes }\left( 
\begin{matrix}
\varkappa _{+}^{0}\backslash \vartheta _{+}^{0} & \varkappa
_{+}^{-}\backslash \sigma _{+}^{-}\cr\varkappa _{0}^{0}\sqcup \vartheta
_{+}^{0} & \vartheta _{0}^{-}\sqcup \vartheta _{+}^{-}%
\end{matrix}%
\right) ^{\ast }\mathbf{f}^{\otimes }\left( 
\begin{matrix}
\varkappa _{0}^{-}\backslash \vartheta _{0}^{-} & \varkappa
_{+}^{-}\backslash \rho _{+}^{-}\cr\varkappa _{0}^{0}\sqcup \vartheta
_{0}^{-} & \vartheta _{+}^{-}\sqcup \vartheta _{+}^{0}%
\end{matrix}%
\right) \ .
\end{equation*}%
As the operator-valued functions $X\otimes \mathbf{f}^{\otimes }(\pmb%
\varkappa )$ are relatively bounded $\Vert X\otimes \mathbf{f}^{\otimes }(%
\pmb\varkappa )\Vert (\pmb\zeta )=\Vert X\Vert $ with respect to $\zeta
_{\nu }^{\mu }(x)=\Vert f_{\nu }^{\mu }(x)\Vert $ and their linear span is
dense in inductive space $\mathcal{U}$, the product formula can be obtained
as a limit for any $T\in \mathcal{U}$, and 
\begin{equation*}
\Vert (T^{{\star }}\cdot T)(\pmb\varkappa )\Vert \leq \sum \Vert T^{{\star }%
}\left( 
\begin{matrix}
\vartheta _{0}^{-}\sqcup \vartheta _{+}^{-} & \varkappa _{+}^{-}\backslash
\sigma _{+}^{-}\cr\varkappa _{0}^{0}\sqcup \vartheta _{+}^{0} & \varkappa
_{+}^{0}\backslash \vartheta _{+}^{0}%
\end{matrix}%
\right) \Vert \;\Vert T\left( 
\begin{matrix}
\varkappa _{0}^{-}\backslash \vartheta _{0}^{-} & \varkappa
_{+}^{-}\backslash \rho _{+}^{-}\cr\varkappa _{0}^{0}\backslash \vartheta
_{0}^{-} & \vartheta _{+}^{-}\sqcup \vartheta _{+}^{0}%
\end{matrix}%
\right) \Vert \leq 
\end{equation*}%
\begin{equation*}
\Vert T\Vert ^{2}(\pmb\zeta )\sum \pmb\zeta ^{\otimes }\left( 
\begin{matrix}
\varkappa _{+}^{0}\backslash \vartheta _{+}^{0}, & \varkappa
_{+}^{-}\backslash \sigma _{+}^{-}\cr\varkappa _{0}^{0}\sqcup \vartheta
_{+}^{0}, & \vartheta _{0}^{-}\sqcup \vartheta _{+}^{-}%
\end{matrix}%
\right) \pmb\zeta ^{\otimes }\left( 
\begin{matrix}
\varkappa _{0}^{-}\backslash \vartheta _{0}^{-} & \varkappa
_{+}^{-}\backslash \rho _{+}^{-}\cr\varkappa _{0}^{0}\backslash \vartheta
_{0}^{-} & \vartheta _{+}^{-}\sqcup \vartheta _{+}^{0}%
\end{matrix}%
\right) =[\Vert T\Vert (\pmb\zeta )]^{2}(\pmb\zeta ^{2})^{\otimes }(\pmb%
\varkappa )
\end{equation*}%
this means $\Vert T^{{\star }}\cdot T\Vert (\pmb\zeta ^{{\star }}\pmb\zeta
)\leq \lbrack \Vert T\Vert (\pmb\zeta )]^{2}$. Due to the proven continuity
of the linear map $\epsilon $ on $\mathcal{U}$ into the $\ast $-algebra of
relatively bounded operators on the projective limit $\cap _{\xi >0}\mathcal{%
G}(\xi )$, the $\star $-multiplicative property of $\epsilon $ can be
extended on the whole $\star $-algebra $\mathcal{U}$ with the unity $I(\pmb%
\varkappa )=I\otimes \mathbf{1}^{\otimes }(\pmb\varkappa )$, $\mathbf{1}(x)$
is the identity matrix, having the representation $\epsilon (I)=\hat{I}$.

Now let us find the representation $U^{t}$ of the multiple quantum integral
(2.7), having the values $\epsilon \circ N_{[0,t)}(L)$ in $\mathcal{U}$ for
relatively bounded operator-valued functions $L(\pmb\vartheta ,\pmb\varkappa
)$ due to 
\begin{equation*}
\Vert T^{t}(\pmb\varkappa )\Vert \leq c\sum_{\vartheta _{0}^{0}\subseteq
\varkappa _{0}^{0}}^{t(\vartheta _{0}^{0})<t}\sum_{\vartheta
_{+}^{0}\subseteq \varkappa _{+}^{0}}^{t(\vartheta
_{+}^{0})<t}\sum_{\vartheta _{0}^{-}\subseteq \varkappa
_{0}^{-}}^{t(\vartheta _{0}^{-})<t}\sum_{\vartheta _{+}^{-}\subseteq
\varkappa _{+}^{-}}^{t(\vartheta _{+}^{-})<t}\Vert L(\pmb\vartheta ,\pmb%
\varkappa \backslash \pmb\vartheta )\Vert 
\end{equation*}%
\begin{equation*}
\leq c\dprod_{\nu =0,+}^{\mu =-,0}\sum_{\vartheta _{\nu }^{\mu }\leq
\varkappa _{\nu }^{\mu }}^{t(\vartheta _{\nu }^{\mu })<t}\lambda _{\nu
}^{\mu }(\vartheta _{\nu }^{\mu })\xi _{\nu }^{\mu }(\varkappa _{\nu }^{\mu
}\backslash \vartheta _{\nu }^{\mu })=c\dprod_{\nu =0,+}^{\mu =-,0}\zeta
_{\nu }^{\mu }(\varkappa _{\nu }^{\mu })\ ,
\end{equation*}%
where 
\begin{equation*}
\zeta (\varkappa )=\dprod_{x\in \varkappa }^{t(x)<t}[\lambda (x)+\xi
(x)]\dprod_{x\in \varkappa }^{t(x)\geq t}\xi (x)
\end{equation*}%
for 
\begin{equation*}
\lambda (\vartheta )=\dprod_{x\in \vartheta }\lambda (x),\quad \xi
(\varkappa )=\dprod_{x\in \varkappa }\xi (x)\ .
\end{equation*}

From the definitions (2.1) of $U^{t}=\epsilon (T^{t})$ we obtain the QS
integral (1.5): 
\begin{equation*}
\lbrack U^{t}a](\varkappa )=\sum_{\varkappa _{0}^{0}\sqcup \varkappa
_{+}^{0}=\varkappa }\iint \sum_{\pmb\vartheta \subseteq \pmb\varkappa ^{t}}L(%
\pmb\vartheta ,\pmb\varkappa \backslash \pmb\vartheta )a(\varkappa
_{0}^{0}\sqcup \varkappa _{0}^{-})\mathrm{d}\varkappa _{0}^{-}\mathrm{d}%
\varkappa _{+}^{-}
\end{equation*}%
\begin{equation*}
=\sum_{\vartheta \sqcup \vartheta _{+}\subseteq \varkappa ^{t}}\int_{%
\mathcal{X}^{t}}\int_{\mathcal{X}^{t}}\sum_{\varkappa _{0}^{0}\sqcup
\varkappa _{+}^{0}=\vartheta _{-}}\iint L(\pmb\vartheta ,\pmb\varkappa )\dot{%
a}(\vartheta \sqcup \vartheta ^{-},\varkappa _{0}^{0}\sqcup \varkappa
_{0}^{-})\mathrm{d}\varkappa _{0}^{-}\mathrm{d}\varkappa _{+}^{-}\mathrm{d}%
\vartheta ^{-}\mathrm{d}\vartheta _{+}^{-}\ ,
\end{equation*}%
where $\vartheta _{-}=\varkappa \backslash (\vartheta \sqcup \vartheta _{+}),%
\dot{a}(\vartheta ,\varkappa _{0}^{0})=a(\vartheta \sqcup \varkappa _{0}^{0})
$. Hence $\epsilon (T^{t})=\Lambda _{\lbrack 0,t)}(B)$ with 
\begin{equation*}
\lbrack B(\pmb\vartheta )\dot{a}(\vartheta \sqcup \vartheta ^{-})](\varkappa
)=\sum_{\varkappa _{0}^{0}\sqcup \varkappa _{+}^{0}=\varkappa }\iint L(\pmb%
\vartheta ,\pmb\varkappa )\dot{a}(\vartheta \sqcup \vartheta ^{-},\varkappa
_{0}^{0}\sqcup \varkappa _{0}^{-})\mathrm{d}\varkappa _{0}^{-}\mathrm{d}%
\varkappa _{+}^{-}\ ,
\end{equation*}%
that is $B(\pmb\vartheta )=\epsilon (L(\pmb\vartheta ))$. In particular, if $%
U^{t}=U^{0}+\Lambda ^{t}(\mathbf{D})$ with $U^{0}=\epsilon (T^{0})$ and $%
\mathbf{D}(x)=\epsilon (\mathbf{C}(x))$, then $U^{t}=\epsilon (T^{0}+N^{t}(%
\mathbf{C}))$, i.e. $\epsilon \circ N^{t}=\Lambda ^{t}\;\circ \epsilon $,
where 
\begin{equation*}
N^{t}(\mathbf{C})(\pmb\varkappa )=\sum_{x\in {\pmb\varkappa }^{t}}C(\mathbf{x%
},\pmb\varkappa \backslash \mathbf{x}),\quad C(\mathbf{x}_{\nu }^{\mu },\pmb%
\varkappa )=C_{\nu }^{\mu }(x,\pmb\varkappa )\ .
\end{equation*}

In the case $U^{t}=\Lambda _{\lbrack 0,t)}(B)$ with $B=\epsilon (L)$ the QS
derivatives%
\begin{equation*}
D_{\nu }^{\mu }(x)=\Lambda _{\lbrack 0,t(x))}(\dot{B}(\mathbf{x}_{\nu }^{\mu
}))=\epsilon (C_{\nu }^{\mu }(x))
\end{equation*}%
are defined by 
\begin{equation*}
C_{\nu }^{\mu }(x,\pmb\varkappa )=N_{[0,t(x))}(\dot{L}(\mathbf{x}_{\nu
}^{\mu })=\dot{T}^{t(x)]}(\mathbf{x}_{\nu }^{\mu },\pmb\varkappa )-\dot{T}%
^{t(x)}(\mathbf{x}_{\nu }^{\mu },\pmb\varkappa )\ ,
\end{equation*}%
where 
\begin{equation*}
N_{[0,t)}(\dot{L}(\mathbf{x}))=\sum_{\pmb\vartheta \subseteq \pmb\varkappa
^{t}}L(\pmb\vartheta \sqcup \mathbf{x},\pmb\varkappa \backslash \pmb%
\vartheta ),
\end{equation*}%
$\mathbf{x}$ is one of the four elementary tables $\mathbf{x}_{\nu }^{\mu }$
and 
\begin{equation*}
\dot{T}^{t(x)}(\mathbf{x},\pmb\varkappa )=\sum_{\pmb\vartheta \subseteq \pmb%
\varkappa ^{t(x)}}L(\pmb\vartheta ,\pmb\varkappa \sqcup \mathbf{x}%
)\backslash \pmb\vartheta )=T^{t(x)}\quad (\pmb\varkappa \sqcup \mathbf{x})\
,
\end{equation*}%
\begin{equation*}
\dot{T}^{t(x)]}(\mathbf{x},\pmb\varkappa )=\sum_{\pmb\vartheta \subseteq \pmb%
\varkappa ^{t(x)}\sqcup \mathbf{x}}L(\pmb\vartheta ,\pmb\varkappa \sqcup 
\mathbf{x})\backslash \pmb\vartheta )=T^{t(x)}\quad (\pmb\varkappa \sqcup 
\mathbf{x})
\end{equation*}%
\begin{equation*}
+\sum_{\pmb\vartheta \subseteq \varkappa ^{t(x)}}L(\pmb\vartheta \sqcup 
\mathbf{x},\pmb\varkappa \backslash \pmb\vartheta )=\dot{T}^{t(x)}(\mathbf{x}%
,\pmb\varkappa )+N_{[0,t(x))}(\dot{L}(\mathbf{x}))(\pmb\varkappa )
\end{equation*}%
due to 
\begin{equation*}
T^{s]}(\pmb\varkappa )=\sum_{\pmb\vartheta \subseteq \varkappa ^{s]}}L(\pmb%
\vartheta ,\pmb\varkappa \backslash \pmb\vartheta )=T^{s_{+}}(\pmb\varkappa )
\end{equation*}%
where 
\begin{equation*}
\varkappa ^{s]}=\{x\in \varkappa |t(x)\leq s\}\ ,\quad s_{+}=\min
\{t(x)>s|x\in \varkappa \}\ .
\end{equation*}%
Hence the derivatives $D_{\nu }^{\mu }(x),x\in X^{t}$, defining $%
U^{t}-U^{0}=\Lambda _{\lbrack 0,t)}(\mathbf{D})$ are represented as the
differences 
\begin{equation*}
D_{\nu }^{\mu }(x)=\epsilon \lbrack \dot{T}^{t(x)]}(\mathbf{x}_{\nu }^{\mu
})]-\epsilon \lbrack \dot{T}^{t(x)}(\mathbf{x}_{\nu }^{\mu })]
\end{equation*}%
of the operators (2.8), where $\dot{T}^{s]}(\mathbf{x},\pmb\varkappa )=T^{t}(%
\pmb\varkappa \sqcup \mathbf{x})=\dot{T}^{t}(\mathbf{x},\pmb\varkappa )$ for
any $t:s<t\leq s_{+}=\min \{t(x)>s|x\in \pmb\varkappa \}$.

Let us consider $\dot{T}^{t}(\mathbf{x})$ as elements $T_{\nu }^{\mu }(x)=%
\dot{T}(\mathbf{x}_{\nu }^{\mu })$ of the triangular operator-valued
matrix-function $\mathbf{T}^{t}(x)$ with $T_{\nu }^{\mu }=0$ for $\mu >\nu $%
, and $T_{-}^{-}(x)=T^{t}=T_{+}^{+}(x)$ independent of $x\in X$, defining
the triangular matrices $\mathbf{U}=[U_{\nu }^{\mu }]$ and $\mathbf{G}%
=[G_{\nu }^{\mu }]$ as $\mathbf{U}(x)=\epsilon (\mathbf{T}^{t(x)}(x))$ and $%
\mathbf{G}(x)=\epsilon (\mathbf{T}^{t(x)]}(x))$, where $\mathbf{T}^{s]}=%
\mathbf{T}^{t}$ for a $t\in (s,s_{+}]$. This helps to generalize the QS It%
\^{o} formula [1] for nonadapted processes as 
\begin{equation*}
U^{t{\ast }}U^{t}-U^{0{\ast }}U^{0}=\Lambda ^{t}(\mathbf{U}^{{\star }}%
\mathbf{D}+\mathbf{D}^{{\star }}\mathbf{U}+\mathbf{D}^{{\star }}\mathbf{D}),
\end{equation*}%
because 
\begin{equation*}
U^{t{\ast }}U^{t}=\epsilon (T^{t{\star }}\cdot T^{t}),(T^{{\star }}\cdot T)(%
\pmb\varkappa \sqcup \mathbf{x}_{\nu }^{\mu })=(\mathbf{T}^{{\star }%
}(x)\cdot \mathbf{T}(x))_{\nu }^{\mu }(\pmb\varkappa )\ ,
\end{equation*}%
as it follows directly from the formula (2.5), in terms of the usual product
of triangular matrices $T^{{\star }}$ and $T$, defined by the multiplication 
$\cdot $ of the matrix elements $T_{\nu }^{\mu }(x)$ and $\star $%
-multiplicative property 
\begin{equation*}
\epsilon (\mathbf{T}^{{\star }}(x)\cdot \mathbf{T}(x))\;=\;\epsilon (\mathbf{%
T}(x))^{{\star }}\epsilon (\mathbf{T}(x))\;.
\end{equation*}

Applying this for $t=t(x)$ and $t=t_{+}(x)=\min \{t\in t(\varkappa )|t>t(x)\}
$ to the representation 
\begin{equation*}
\epsilon \lbrack (\mathbf{T}^{t(x)]{\star }}\mathbf{T}^{t(x)]})(x)-(\mathbf{T%
}^{t(x){\star }}\mathbf{T}^{t(x)})(x)]
\end{equation*}%
of the QS derivative of the process $U^{t{\ast }}U^{t}$, we finally obtain
the formula 
\begin{equation*}
\mathrm{d}(U^{t{\ast }}U^{t})=\mathrm{d}\Lambda ^{t}[\epsilon (\mathbf{T}%
^{t]})^{{\star }}\epsilon (\mathbf{T}^{t]})-\epsilon (\mathbf{T}^{t})^{{%
\star }}\epsilon (\mathbf{T}^{t})]\ ,
\end{equation*}%
giving the multiplication table (2.9) in terms of the triangular matrices
(2.8) with $G_{-}^{-}(x)=U^{t(x)}=G_{+}^{+}(x)\
,U_{-}^{-}(x)=U^{t(x)}=U_{+}^{+}(x)$ and $D_{\nu }^{\mu }(x)=G_{\nu }^{\mu
}(x)-U_{\nu }^{\mu }(x)\ .$
\end{proof}

\begin{corollary}
The QS process $U^{t}=\epsilon (T^{t})$ is adapted, iff $T^{t}(\pmb\varkappa
)=T(\pmb\varkappa ^{t})\otimes 1(\pmb\varkappa _{\lbrack t})$ for almost all 
$\pmb\varkappa =(\varkappa _{\nu }^{\mu })$, where $\pmb\varkappa ^{t}=\pmb%
\varkappa \cap X^{t},\pmb\varkappa _{\lbrack t}=\pmb\varkappa \cap X_{[t}$,
and $1(\pmb\varkappa )=I(\varkappa _{0}^{0})$ for $\varkappa _{\nu }^{\mu
}=\emptyset ,\mu \not=\nu $, otherwise $1(\pmb\varkappa )=0$. The QS It\^{o}
formula for adapted processes $U^{t}$ can be written in the form 
\begin{equation*}
\mathrm{d}(U^{\ast }U)=\mathrm{d}\Lambda (\mathbf{G}^{{\star }}\mathbf{G}-%
\mathbf{U}^{{\ast }}\mathbf{U}\otimes \mathbf{1})=U^{\ast }\mathrm{d}U+%
\mathrm{d}U^{\ast }U+\mathrm{d}U^{\ast }\mathrm{d}U\ ,\eqno(2.10)
\end{equation*}%
where $\mathrm{d}U^{\ast }\mathrm{d}U=\mathrm{d}\Lambda (\mathbf{D}^{{\star }%
}\mathbf{D})$ is defined by the usual product of the triangular matrices $%
\mathbf{D}=[D_{\nu }^{\mu }]$, $\mathbf{D}^{{\star }}=(\hat{I}\otimes 
\mathbf{g})\mathbf{D}^{\ast }(\hat{I}\otimes \mathbf{g})$ with $D_{\nu
}^{\mu }=0$, if $\mu =+$ or $\nu =-\ .$
\end{corollary}

Indeed, if $T^{t}(\pmb\varkappa )=T(\pmb\varkappa ^{t})\otimes 1(\pmb%
\varkappa _{\lbrack t})$, then 
\begin{equation*}
\lbrack U^{t}(a^{t}\otimes c)](\varkappa )=\sum_{\varkappa _{0}^{0}\sqcup
\varkappa _{+}^{0}=\varkappa ^{t}}\iint T(\pmb\varkappa ^{t})a(\varkappa
_{0}^{0}\sqcup \varkappa _{0}^{-})\otimes c(\varkappa _{\lbrack t})\mathrm{d}%
\varkappa _{0}^{-}\mathrm{d}\varkappa _{+}^{-}\ ,
\end{equation*}%
for any $a^{t}\in \mathcal{G}^{t},c\in \mathcal{F}_{[t}$, where the integral
should be taken over $\varkappa _{0}^{-},\varkappa _{+}^{-}\in \mathcal{X}%
^{t}$, otherwise $T^{t}(\pmb\varkappa )=0$. Hence $U^{t}(a^{t}\otimes
c)=b^{t}\otimes c$ for a $b^{t}\in \mathcal{G}^{t}$. In this case 
\begin{equation*}
\dot{T}^{t(x)}(\mathbf{x}_{\nu }^{\mu },\pmb\varkappa )=T^{t(x)}(\pmb%
\varkappa \sqcup \mathbf{x}_{\nu }^{\mu })=T(\pmb\varkappa ^{t(x)})\otimes
1_{\nu }^{\mu }(x)\otimes 1(\pmb\varkappa _{t(x)})\ ,
\end{equation*}%
where $1_{\nu }^{\mu }(x)=0$, if $\mu \not=\nu
,1_{-}^{-}=1=1_{+}^{+},1_{0}^{0}(x)=I(x)$. This gives $U_{\nu }^{\mu
}(x)=\epsilon (\dot{T}^{t(x)}(\mathbf{x}_{\nu }^{\mu }))=U^{t}\otimes 1_{\nu
}^{\mu }(x),$ and 
\begin{equation*}
\mathrm{d}\Lambda (G^{{\star }}G-U^{\ast }U\otimes \mathbf{1})=\mathrm{d}%
\Lambda ((U^{\ast }\otimes \mathbf{1})\mathbf{D}+\mathbf{D}^{{\star }%
}(U\otimes \mathbf{1})+\mathbf{D}^{{\star }}\mathbf{D})=
\end{equation*}%
\begin{equation*}
U^{\ast }\mathrm{d}\Lambda (\mathbf{D})+\mathrm{d}\Lambda (\mathbf{D}^{{%
\star }})U+\mathrm{d}\Lambda (\mathbf{D}^{{\star }}\mathbf{D})=U^{\ast }%
\mathrm{d}U+\mathrm{d}U^{\ast }U+\mathrm{d}U^{\ast }\mathrm{d}U
\end{equation*}%
as it follows from corollary 1 for the adaptive $U^{t}$.

\section{A nonadapted QS evolution and chronological products}

The proved $\star $-homomorphism and continuity properties of the
representation $\epsilon $ of the unital inductive $\star $-algebra $%
\mathcal{U}$ of all operator-valued functions $T(\pmb\varkappa )$ of $%
\varkappa _{\nu }^{\mu }\in \mathcal{X}$, $(\mu ,\nu )\in \{-,0\}\times
\{0,+\}$, relatively bounded with respect to some $\pmb\zeta =[\zeta _{\nu
}^{\mu }(x)]$, into the $\star $-algebra $\mathcal{B}$ of all relatively
bounded operators on the projective limit $\mathcal{G}^{+}=\bigcap_{\xi >0}%
\mathcal{G}(\xi )$ enables us to construct a QS functional calculus.

Namely, if $T=f(Q_{1},\dots ,Q_{m})$ is an analytical function of $Q_{i}\in 
\mathcal{U}$ as a limit in $\mathcal{U}$ of polynomials $T_{n}$ with some
ordering of noncommuting $Q_{1},\dots ,Q_{m}$ in the sense $\Vert
T_{n}-T\Vert (\pmb\zeta )\rightarrow 0$ for a $\pmb\zeta $, then $U=\epsilon
(T)$ is the ordered function $f(X_{1},\dots ,X_{m})$ of $X_{i}=\epsilon
(Q_{i})$ as a limit on $\mathcal{G}^{+}$ of the corresponding polynomials $%
U_{n}=\epsilon (T_{n})$, that is $\Vert U_{n}-U\Vert _{\xi ^{+}}^{\xi
_{-}}\rightarrow 0$ for any $\xi _{-}>0$ and $\xi ^{+}>\xi _{-}\Vert \zeta
_{0}^{0}\Vert _{\infty }^{2}$. The function $U^{{\ast }}=f^{{\ast }}(X_{1}^{{%
\ast }},\dots ,X_{n}^{{\ast }})$ with the transposed ordering as $U^{{\ast }%
}=\epsilon (T^{{\star }})$ for $T^{{\star }}=f^{{\ast }}(Q_{1}^{{\star }%
},\dots ,Q_{n}^{{\star }})$ is also defined as $(\xi ^{+},\xi _{-})$-limit
due to $\Vert T_{n}^{{\star }}-T^{{\star }}\Vert (\pmb\zeta ^{{\star }%
})\rightarrow 0$ and $(\pmb\zeta ^{{\star }})_{0}^{0}=\zeta _{0}^{0}$.

The differential form of this calculus is given by the noncommutative and
nonadaptive generalization of the QS It\^{o} formula 
\begin{equation*}
\mathrm{d}X=\mathrm{d}\Lambda (\mathbf{A})\Rightarrow \mathrm{d}f(X)=\mathrm{%
d}\Lambda (f(\mathbf{X}+\mathbf{A})-f(\mathbf{X}))\eqno(3.1)
\end{equation*}%
defined for any analytical function $U^{t}=f(X^{t})$ of $\epsilon (Q^{t})$
as QS differential of $\epsilon (T^{t})$ for $T^{t}=f(Q^{t})$ with 
\begin{equation*}
U_{\nu }^{\mu }(x)=f(\mathbf{X})_{\nu }^{\mu }(x),\quad G_{\nu }^{\mu }(x)=f(%
\mathbf{X}+\mathbf{A})_{\nu }^{\mu }(x),
\end{equation*}%
where $f(\mathbf{Z})(x)=f(\mathbf{Z}(x))$ is the triangular matrix which is
the function of the matrix $\mathbf{Z}(x)$, representing $\mathbf{Q}%
^{t(x)}(x)$ and $\mathbf{Q}^{t(x)]}(x)$ correspondingly as $\mathbf{X}%
(x)=\epsilon (\mathbf{Q}^{t(x)}(x))$ and%
\begin{equation*}
\mathbf{X}(x)+\mathbf{A}(x)\ ,\ \;\;A_{\nu }^{\mu }(x)=\epsilon (\dot{Q}%
^{t(x)]}(\mathbf{x}_{\nu }^{\mu })-\dot{Q}^{t(x)}(\mathbf{x}_{\nu }^{\mu })).
\end{equation*}%
For the ordered functions $U^{t}=f(X_{1}^{t},\dots ,X_{n}^{t})$ this can be
written in terms of $\mathbf{X}_{i}$ with $\mathrm{d}\mathbf{X}_{i}=\mathrm{d%
}\Lambda (\mathbf{A}_{i})$ and $\mathbf{Z}_{i}=\mathbf{X}_{i}+\mathbf{A}_{i}$
as 
\begin{equation*}
\mathrm{d}U=\mathrm{d}\Lambda (f(\mathbf{Z}_{1},\dots ,\mathbf{Z}_{n})-f(%
\mathbf{X}_{1},\dots ,\mathbf{X}_{n}))\ .\eqno(3.2)
\end{equation*}%
In particular, if all the triangular matrices $\{\mathbf{X}_{i},\mathbf{Z}%
_{i}\}$ are commutative, then one can obtain the exponential function $%
U^{t}=\exp \{X^{t}\}$ for $X=\sum X_{i}$ as the solution of the following QS
differential equation 
\begin{equation*}
\mathrm{d}U=\mathrm{d}\Lambda \left[ (\exp \{\mathbf{A}\}-\hat{\mathbf{I}})\ 
\mathbf{U}\right] \ ,\eqno(3.3)
\end{equation*}%
with $\mathbf{A}=\sum \mathbf{A}_{i}$ and the initial condition $U^{0}=\hat{I%
}$. Now we shall study the problem of the solution of the general QS
evolution equation 
\begin{equation*}
\mathrm{d}U=\mathrm{d}\Lambda (\left( \mathbf{S}-\hat{\mathbf{I}}\right) \ 
\mathbf{U})=\mathrm{d}\Lambda \left( \mathbf{B}\ \mathbf{U}\right) ,\eqno%
(3.4)
\end{equation*}%
defined by a matrix-function $\mathbf{B}(x)=[B_{\nu }^{\mu }(x)]$ with
noncommutative operator-values 
\begin{align*}
B_{0}^{0}(x)& \colon \mathcal{G}\otimes \mathcal{E}(x)\rightarrow \mathcal{G}%
\otimes \mathcal{E}(x),\quad B_{+}^{-}(x)\colon \mathcal{G}\rightarrow 
\mathcal{G}\ , \\
B_{+}^{0}(x)& \colon \mathcal{G}\rightarrow \mathcal{G}\otimes \mathcal{E}%
(x),\quad B_{0}^{-}(x)\colon \mathcal{G}\otimes \mathcal{E}(x)\rightarrow 
\mathcal{G}\ ,
\end{align*}%
and $B_{\nu }^{\mu }(x)=0$, if $\mu =+\qquad \text{or}\qquad \nu =-$. In the
adapted case this equation can be written according to Corollary 1 as $%
\mathrm{d}U=\mathrm{d}\Lambda (\mathbf{B})U$, which shows that its solution
with $U^{0}=\hat{I}$ should be defined in some sense as a chronologically
ordered exponent $U^{t}=\Gamma _{\lbrack 0,t)}(\mathbf{B})$. In particular,
if $B_{\nu }^{\mu }(x)=\hat{I}\otimes l_{\nu }^{\mu }(x)$, where $\mathbf{l}%
=[l_{\nu }^{\mu }]$ is a triangular QS-integrable matrix-function with $%
l_{\nu }^{\mu }=0$, if $\mu =+$ or $\nu =-$, 
\begin{equation*}
l_{0}^{0}(x)\colon \mathcal{E}(x)\rightarrow \mathcal{E}(x),\;\Vert
l_{0}^{0}\Vert _{\infty }^{t}<\infty ;\;l_{+}^{0}(x)\in \mathcal{E}%
(x),\;l_{0}^{-}(x)\in \mathcal{E}^{{\ast }}(x),\;\Vert l\Vert
_{2}^{t}<\infty ;\Vert l_{+}^{-}\Vert _{1}^{t}<\infty ,
\end{equation*}%
then $U^{t}$ is defined as $I\otimes \Gamma _{\lbrack 0,t)}(\mathbf{l})$,
where $\Gamma _{\lbrack 0,t)}(\mathbf{l})=\epsilon \left( \mathbf{f}%
_{[0,t)}^{\otimes }\right) $ is the representation (2.1) of $\mathbf{f}%
_{[0,t)}^{\otimes }(\pmb\varkappa )=\otimes _{x\in \pmb\varkappa }\mathbf{f}%
^{t}(x)$ with $\mathbf{f}^{t}(x)=\mathbf{1}(x)+\mathbf{l}^{t}(x)$, $\mathbf{l%
}^{t}(x)=\mathbf{l}(x)$ if $t(x)<t$ and $\mathbf{l}^{t}(x)=0$, if $t(x)\geq
t,$ i.e. $\Gamma _{\lbrack 0,t)}(\mathbf{l})$ is the second quantization $%
\Gamma (\mathbf{l}^{t})$ of $\mathbf{l}^{t}$.

\begin{theorem}
The QS evolution equation (3.4) written in the integral form $%
U^{t}=U^{0}+\Lambda ^{t}(\mathbf{B}\mathbf{U})$ with $U^{0}=\epsilon (T^{0})$%
, $\mathbf{B}(x)=\epsilon (\mathbf{L}(x))$, is the representation $\epsilon $
of the recurrences 
\begin{equation*}
T^{t_{+}}(\pmb\varkappa )=\left[ F_{t(x)}\cdot T^{t}\right] (\pmb\varkappa
),x\in \varkappa =\sqcup \varkappa _{\nu }^{\mu }\ ,\eqno(3.5)
\end{equation*}%
defined for any partition $\pmb\varkappa =(\varkappa _{\nu }^{\mu })$ of a
chain $\varkappa \in \mathcal{X}$ with $t\in (t_{-}(x),t(x)]$,$%
\;t_{-}(x)=\max \;t\left( \varkappa ^{t(x)}\right) $ and $t_{+}\in
(t(x),t_{+}(x)]$, $t_{+}(x)=\min t(\varkappa _{t(x)})$, by 
\begin{equation*}
F_{t(x)}\left( \pmb\varkappa \sqcup \mathbf{x}_{\nu }^{\mu }\right) =L_{\nu
}^{\mu }(x,\pmb\varkappa )+I\left( \pmb\varkappa \sqcup \mathbf{x}_{\nu
}^{\mu }\right) \equiv F_{\nu }^{\mu }(x,\pmb\varkappa ),
\end{equation*}%
where $\mathbf{x}_{\nu }^{\mu }$ is one of the four single point tables $\pmb%
\vartheta =(\vartheta _{\lambda }^{\kappa })$ with $\vartheta _{\nu }^{\mu
}=x$.

The recurrency (3.5) with the initial condition $T^{t}(\pmb\varkappa )=T^{0}(%
\pmb\varkappa )$ for all $t\leq \min t(\varkappa )$ has the unique solution 
\begin{equation*}
T^{t}(\pmb\varkappa )=[F_{[0,t)}\cdot T^{0}](\pmb\varkappa
),\;\;\;\;F_{[0,t)}=\bullet _{0\leq s<t}^{\leftarrow }F_{s}\ ,\eqno(3.6)
\end{equation*}%
where $F_{s}(\pmb\varkappa )=I(\pmb\varkappa ),$ if $s\notin \sqcup
\varkappa _{\nu }^{\mu }=\varkappa $, defined for every chain $\varkappa
=(x_{1},\dots ,x_{m},\dots )\in \mathcal{X}$, $t\in (t_{m-1},t_{m}]$ as
product%
\begin{equation*}
\bullet _{x\in \varkappa t}^{\leftarrow }F_{t(x)}=F_{t_{m-1}}\cdot \cdot
\cdot F_{t_{1}}
\end{equation*}%
of $F_{t_{i}},t_{i}=t(x_{i})<t_{i+1}$ and $T^{0}$ in chronological order.
The solution $U^{t}=\epsilon (T^{t})$ of (3.4) is isometric $U^{{\ast }}U=%
\hat{I}$ (unitary: $U^{{\ast }}=U^{-1}$) up to a $t>0$, if $U^{0}$ is
isometric (unitary) and the triangular matrix-function $\mathbf{S}(x)=%
\mathbf{B}(x)+\hat{\mathbf{I}}(x)=\epsilon (\mathbf{F}(\varkappa ))$ is
pseudoisometric:%
\begin{equation*}
\mathbf{S}^{{\star }}(x)\mathbf{S}(x)=\hat{\mathbf{I}}(x)=\hat{I}\otimes 
\mathbf{1}(x)
\end{equation*}%
(pseudounitary $\mathbf{S}^{{\star }}(x)=\mathbf{S}(x)^{-1}$) for almost all 
$x\in X^{t}$, that is 
\begin{equation*}
S_{0}^{0}(x)^{{\ast }}S_{0}^{0}(x)=I(x),\quad S_{+}^{-}(x)^{{\ast }%
}+S_{+}^{0}(x)^{{\ast }}S_{+}^{0}(x)+S_{+}^{-}(x)=0\eqno(3.7)
\end{equation*}%
\begin{equation*}
S_{+}^{-}(x)^{{\ast }}+S_{0}^{0}(x)^{{\ast }}S_{+}^{0}(x)=0,\quad
S_{+}^{0}(x)^{{\ast }}S_{0}^{0}(x)+S_{0}^{-}(x)=0
\end{equation*}%
(and $S_{0}^{0}(x)$ is unitary $S_{0}^{0}(x)^{{\ast }}=S_{0}^{0}(x)^{-1}$
for almost all $x\in X^{t}$).
\end{theorem}

\begin{proof}
We are looking for the solution of the equation (3.4) as the representation $%
U=\epsilon (T)$ of some $T(\pmb\varkappa )$. If $\mathbf{B}(x)=\epsilon (%
\mathbf{L}(x))$, then $\mathbf{B}(x)\mathbf{U}(x)=\epsilon (\mathbf{L}(x)%
\mathbf{T}(x))$ and $\Lambda ^{t}(\mathbf{B}\mathbf{U})=\epsilon (N^{t}(%
\mathbf{LT})$ due to the property $\Lambda \circ \epsilon =\epsilon \circ N$%
, proved in theorem 2, and the multiplicative property $\epsilon (\mathbf{LT}%
)=\mathbf{B}\mathbf{U}$, where $\mathbf{U}(x)=\epsilon (\mathbf{T}(x))$, $%
\mathbf{T}(x)$ denotes the triangular matrix $[T_{\nu }^{\mu }(x)],T_{\nu
}^{\mu }=0$, if $\mu >\nu $, with $T_{-}^{-}(x)=T^{t(x)}=T_{+}^{+}(x)$ and $%
T_{\nu }^{\mu }(x)=\dot{T}^{t(x)}$ $(\mathbf{x}_{\nu }^{\mu })$ for $\mu
\not=+$, $\nu \not=-$. This gives a possibility to consider the equation
(3.4) in integral form as the representation $U^{t}=\epsilon (T^{0}+N^{t}(%
\mathbf{LT}))$ of the equation 
\begin{equation*}
T^{t}(\pmb\varkappa )=T^{0}(\pmb\varkappa )+N^{t}(\mathbf{LT})(\pmb\varkappa
),
\end{equation*}%
corresponding to $U^{0}=\epsilon (T^{0})$, where 
\begin{equation*}
N^{t}(\mathbf{C})(\pmb\varkappa )=\sum_{\nu =0,+}^{\mu =-,0}\ \sum_{x\in
\varkappa _{\nu }^{\mu }}^{t(x)<t}\left[ \mathbf{L}(x)\mathbf{T}(x)\right]
_{\nu }^{\mu }\left( \pmb\varkappa \backslash \mathbf{x}_{\nu }^{\mu
}\right) 
\end{equation*}%
depends on $T^{s}(\pmb\vartheta )$ with $s=t(x)<t$ for $x\in \sqcup
\varkappa _{\nu }^{\mu }\supseteq \sqcup \vartheta _{\nu }^{\mu }$. This
defines $T^{t}(\pmb\varkappa )$ for any partition $\pmb\varkappa =(\varkappa
_{\nu }^{\mu })$ of a chain $\varkappa =(x_{1},\dots ,x_{n})\in \mathcal{X}$
as the solution $T^{t}(\pmb\varkappa )=T_{m}(\pmb\varkappa ),\qquad
m=|\varkappa ^{t}|$ of the recurrency 
\begin{equation*}
T_{m}(\pmb\varkappa )=T_{0}(\pmb\varkappa )+\sum_{k=1}^{m}\left[ \left(
F_{t_{k}}-I\right) T_{k-1}\right] (\pmb\varkappa )=\left[ F_{t_{m}}T_{m-1}%
\right] (\pmb\varkappa ),
\end{equation*}%
where $T_{k-1}=T^{t_{k}}$ for $t_{k}=t(x_{k})$, and the product $\mathbf{LT}$
for $T_{\nu }^{\mu }(x,\pmb\varkappa )=T^{t(x)}(\pmb\varkappa \sqcup \mathbf{%
x}_{\nu }^{\mu })$ is written as%
\begin{equation*}
\lbrack \mathbf{L}(x)\mathbf{T}(x)]_{\nu }^{\mu }(\pmb\varkappa \backslash 
\mathbf{x}_{\nu }^{\mu })=\left[ \left( F_{t(x)}-I\right) T^{t(x)}\right] (%
\pmb\varkappa )
\end{equation*}%
for $x\in \varkappa _{\nu }^{\mu }$ in terms of $F_{t(x)}(\pmb\varkappa
\sqcup x_{\nu }^{\mu })=F_{\nu }^{\mu }(x,\pmb\varkappa )$ for $\mathbf{F}=%
\mathbf{L}+\mathbf{I}$. So, if the solution of (3.4) exists as $%
U^{t}=\epsilon (\mathbf{T}^{t})$, then it is uniquely defined by (3.6).

Let us suppose, that $(\mathbf{S}^{{\star }}\mathbf{S})(x)=\hat{I}\otimes 
\mathbf{1}(\mathbf{x})$ for almost all $x$ with $t(x)<t$, which is the
representation $\epsilon (\mathbf{F}^{{\star }}\mathbf{F})=\mathbf{S}^{{%
\star }}\mathbf{S}$ of $\left( F_{t(x)}^{{\star }}F_{t(x)}\right) (\pmb%
\varkappa \sqcup \mathbf{x})=I(\pmb\varkappa )\otimes 1(\mathbf{x})$ for
corresponding $\mathbf{F}(x)$, $t(x)<t$. By the recurrency $%
T_{m}=F_{t_{m}}T_{m-1}$ we obtain $(T_{k}^{{\star }}T_{k})(\pmb\varkappa )=I(%
\pmb\varkappa )$ for all $t_{k}<t$ from the initial condition $(T_{0}^{{%
\star }}T_{0})(\pmb\varkappa )=I(\pmb\varkappa )$. Hence, $(T^{t{\star }%
}T^{t})(\pmb\varkappa )=I(\pmb\varkappa )$ for almost all tables $\pmb%
\varkappa =(\varkappa _{\nu }^{\mu })$, namely those which are partitions of
the chains $\varkappa $. This gives $U^{t{\ast }}U^{t}=\hat{I}$ for $%
U^{t}=\epsilon (T^{t})$. In the same way one can obtain the condition $%
U^{t}U^{t{\ast }}=\hat{I}$ from $(\mathbf{SS}^{{\star }})(x)=\hat{I}\otimes 
\mathbf{1}(x)$ for almost all $x$ with $t(x)<t$. Writing the condition $%
\mathbf{S}^{{\star }}\mathbf{S}=\hat{I}\otimes \mathbf{1}$ in terms of
matrix elements, we obtain (3.7): 
\begin{equation*}
(\mathbf{S}^{{\star }}\mathbf{S})(x)=\left[ 
\begin{matrix}
1 & S_{+}^{0}(x)^{{\ast }} & S_{+}^{-}(x)^{{\ast }}\cr0 & S_{0}^{0}(x)^{{%
\ast }} & S_{0}^{-}(x)^{{\ast }}\cr0 & 0 & 1%
\end{matrix}%
\right] \left[ 
\begin{matrix}
1 & S_{0}^{-}(x) & S_{+}^{-}(x)\cr0 & S_{0}^{0}(x) & S_{+}^{0}(x)\cr0 & 0 & 1%
\end{matrix}%
\right] =\hat{I}\otimes \left[ 
\begin{matrix}
1 & 0 & 0\cr0 & I(x) & 0\cr0 & 0 & 1%
\end{matrix}%
\right] \ .
\end{equation*}

The unitary solution of the equation (3.4) under conditions (3.7) in terms
of $\mathbf{B}=\mathbf{S}-\hat{\mathbf{I}}$ was obtained in [1] in the
framework of It\^{o} (adapted) QS calculus for the stationary Markovian case 
$\mathbf{B}(t)=\mathbf{L}\otimes \hat{1}$, for nonstationary finite
dimensional Markovian case $\mathbf{B}(t)=\mathbf{L}(t)\otimes \hat{1}$ in
[11]; and for non Markovian adapted case $\mathbf{B}(t)=\mathbf{L}%
^{t}\otimes 1_{[t}$ in [5].
\end{proof}

\begin{corollary}
Let $U^{0}=T^{0}\otimes \hat{1},T^{0}\in \mathcal{B}(\mathcal{H})$ and $%
\mathbf{S}(x)=\mathbf{F}(x)\otimes \hat{1}$ be defined by the operators $%
F_{\nu }^{\mu }=L_{\nu }^{\mu },\ \mu <\nu ,\ F_{0}^{0}=L_{0}^{0}+I$, acting
as 
\begin{align*}
F_{0}^{0}(x)& \colon \mathcal{H}\otimes \mathcal{E}(x)\rightarrow \mathcal{H}%
\otimes \mathcal{E}(x),\quad F_{+}^{-}:\mathcal{H}\rightarrow \mathcal{H}, \\
F_{+}^{0}(x)& \colon \mathcal{H}\rightarrow \mathcal{H}\otimes \mathcal{E}%
(x),\quad F_{0}^{-}(x)\colon \mathcal{H}\otimes \mathcal{E}(x)\rightarrow 
\mathcal{H}
\end{align*}%
with 
\begin{equation*}
\Vert F_{0}^{0}\Vert _{\infty }^{(t)}<\infty ,\ \Vert F_{+}^{0}\Vert
_{2}^{(t)}<\infty ,\ \Vert F_{0}^{-}\Vert _{2}^{(t)}<\infty ,\ \Vert
F_{+}^{-}\Vert _{1}^{(t)}<\infty \ .\eqno(3.8)
\end{equation*}%
Then the solution $U^{t}=\epsilon \left( T^{t}\right) $ of the evolution
equation (3.4) is defined as $\zeta $-bounded operator for $\zeta >\Vert
F_{0}^{0}\Vert _{\infty }^{(t)}$ by adapted chronological product $T^{t}(\pmb%
\varkappa )=F_{[0,t)}(\pmb\varkappa )\cdot T^{0}$, satisfying the recurrency 
\begin{equation*}
T^{t_{+}}(\pmb\varkappa \sqcup \mathbf{x})=\ F(\mathbf{x})\cdot T^{t}(\pmb%
\varkappa ),\quad F(\mathbf{x}_{\nu }^{\mu })=F_{\nu }^{\mu }(x)\ ,\eqno(3.9)
\end{equation*}%
where $t_{-}(x)<t\leq t(x)<t_{+}\leq t_{+}(x),t_{-}(x)=\max \{t<t(x)|t\in
t(\varkappa )\},t_{+}(x)=\min \{t>t(x)|t\in t(\varkappa )\}$, and $\cdot $
means the semitensor product, defined in theorem 1.

The QS process $U^{t}=\epsilon (T^{t})$ is adapted, can be represented as
multiple QS integral (1.5) $U^{t}=\Lambda _{\lbrack 0,t)}(\mathbf{L}%
^{\triangleleft }\cdot T^{0}\otimes \hat{1})$ with semitensor chronological
products $\mathbf{L}^{\triangleleft }(\pmb\vartheta )=\bullet _{x\in \pmb%
\vartheta }^{\leftarrow }L(\mathbf{x})$ of $\mathbf{L}(x)=\mathbf{F}(x)-%
\mathbf{I}(x)$, and has the estimate 
\begin{equation*}
\Vert U^{t}\Vert _{\xi ^{+}}^{\xi _{-}}\leq \Vert T^{0}\Vert \exp \left\{
\int_{0}^{t}(\Vert L_{+}^{-}(x)\Vert +(\Vert L_{0}^{-}(x)\Vert ^{2}+\Vert
L_{+}^{0}(x)\Vert ^{2})/2\varepsilon \}\mathrm{d}x\right\} \eqno(3.10)
\end{equation*}%
for $\xi ^{+}/\xi _{-}>\mathrm{ess\sup }_{x\in X^{t}}\Vert F_{0}^{0}(x)\Vert 
$ and sufficiently small $\varepsilon >0$.
\end{corollary}

Indeed, if $\mathbf{S}=\mathbf{F}\otimes \hat{1}$ and $\mathbf{F}$ satisfies
the local integrability conditions (3.8), then 
\begin{equation*}
\Vert \bullet _{x\in \pmb\varkappa ^{t}}^{\leftarrow }F(\mathbf{x})\Vert
\leq \dprod_{x\in \pmb\varkappa ^{t}}\Vert F(\mathbf{x})\Vert =\dprod_{\nu
=0,+}^{\mu =-,0}\zeta ^{t}(\varkappa _{\nu }^{\mu })\ ,
\end{equation*}%
where $\zeta ^{t}(\varkappa _{\nu }^{\mu })=\dprod_{x\in \varkappa _{\nu
}^{\mu }}^{t(x)<t}\Vert F_{\nu }^{\mu }(x)\Vert $. Hence, $T^{t}(\pmb%
\varkappa )=F_{[0,t)}(\pmb\varkappa ^{t})\otimes 1(\pmb\varkappa _{\lbrack
t})$ is relatively bounded $\Vert T^{t}\Vert (\pmb\zeta ^{t})\leq \Vert
T^{0}\Vert $ with respect to $\pmb\zeta ^{t}(x)=\left( \Vert F_{\nu }^{\mu
}(x)\Vert \right) _{\nu =0,+}^{\mu =-,0}$, $x\in X^{t}$, and $\pmb\zeta
^{t}(x)=0$, $x\in X_{[t}$. Due to (2.4) and Theorem 2 this gives the $\zeta $%
-boundedness of the operator $U^{t}=\epsilon (t_{[0,t)})$ with respect to $%
\zeta >\Vert F_{0}^{0}\Vert _{\infty }^{(t)}$; and the estimate (3.10) for $%
\xi ^{+}>\xi _{-}\Vert F_{0}^{0}\Vert _{\infty }^{(t)},\varepsilon
<\varepsilon (\xi ^{+},\xi _{-})$ in terms of the norms (3.9) for $F_{\nu
}^{\mu }=L_{\nu }^{\mu }+I\otimes 1_{\nu }^{\mu }$. Taking into account that 
$\epsilon \circ N_{[0,t)}=\Lambda _{\lbrack 0,t)}\circ \epsilon $ and 
\begin{equation*}
N_{[0,t)}\left( \mathbf{L}^{\triangleleft }\cdot T^{0}\right) (\pmb\varkappa
)=\sum_{\pmb\vartheta \subset \pmb\varkappa ^{t}}\mathbf{L}^{\triangleleft }(%
\pmb\vartheta )\cdot T^{0}\otimes 1(\pmb\varkappa \backslash \pmb\vartheta )=%
\left[ \bullet _{x\in \pmb\varkappa ^{t}}^{\leftarrow }(I(\mathbf{x})+L(%
\mathbf{x}))\cdot T^{0}\right] \otimes 1(\pmb\varkappa _{t]})\ ,
\end{equation*}%
we obtain the QS integral representation $\Lambda _{\lbrack 0,t)}(\mathbf{L}%
^{\triangleleft }\cdot T^{0}\otimes \hat{1})$ of Wick chronological product $%
\epsilon (T^{t})$. This process is adapted and has the QS derivative 
\begin{equation*}
\mathbf{D}(x)=\Lambda _{\lbrack 0,t(x))}(\mathbf{L}(x)\cdot \mathbf{L}%
^{\triangleleft }\cdot T^{0}\otimes \hat{1})=\mathbf{L}(x)\odot U^{t(x)}\ ,
\end{equation*}%
where $(L\cdot U)_{\nu }^{\mu }=(L_{\nu }^{\mu }\otimes \hat{1})\cdot U$, $%
B_{0}^{\mu }(x)\cdot U=B_{0}^{\mu }(U\otimes I(x))$, $B_{+}^{\mu }\cdot
U=B_{+}^{\mu }U$. Hence multiple integral $U^{t}=\Lambda _{\lbrack 0,t)}(%
\mathbf{L}^{\triangleleft }\otimes \hat{1})$ satisfies the QS equation
(1.10) with $U^{0}=\hat{I}$ as the case $\mathrm{d}U=\mathrm{d}\Lambda (%
\mathbf{L}\otimes \hat{1})U$ of (3.4).

Finally let us define the solution of the unitary evolution equation (1.10)
with $\mathbf{L}=\mathrm{e}^{-\mathrm{i}\mathbf{H}}-\mathbf{I}$, $\mathbf{H}%
^{{\star }}=\mathbf{H}$, $U^{0}=\hat{I}$ as the representation (2.1) of
chronologically ordered products 
\begin{equation*}
T^{t}(\pmb\varkappa )=\bullet _{x\in \pmb\varkappa }^{\leftarrow }F^{t}(%
\mathbf{x})=\mathbf{F}^{\triangleleft }(\pmb\varkappa ^{t})\otimes 1(\pmb%
\varkappa _{\lbrack t})
\end{equation*}%
for $\mathbf{F}(x)=\exp \{-\mathrm{i}\mathbf{H}(x)\}$ 
\begin{equation*}
\mathbf{F}^{\triangleleft }(\pmb\varkappa )=F(\mathbf{x}_{n})\cdot \cdot
\cdot F(\mathbf{x}_{1})\qquad \text{for}\qquad \pmb\varkappa =\sqcup
_{i=1}^{n}\mathbf{x}_{i}\ .
\end{equation*}%
Here $F^{t}(\mathbf{x})=I(\mathbf{x})=I\otimes 1(\mathbf{x})$, if $t(x)\geq t
$, $F^{t}(\mathbf{x}_{\nu }^{\mu })=F_{\nu }^{\mu }(x)$, if $t<t(x)$, $\pmb%
\varkappa =(\varkappa _{\nu }^{\mu })$ is a partition $\varkappa =\sqcup
_{\nu =0,+}^{\mu =-,0}\varkappa _{\nu }^{\mu }$ of a chain $\varkappa
=(x_{1},\dots ,x_{n})\in \mathcal{X}$ ordered by $t(x_{i-1})<t(x_{i})$ with $%
x_{i}\in \varkappa _{\nu }^{\mu }$, corresponding to the single point table $%
\mathbf{x}_{i}=(\varkappa _{\nu }^{\mu })_{\nu =0,+}^{\mu =-,0}$ with $%
\varkappa _{\nu }^{\mu }=x_{i}$, and $F(\mathbf{x})\cdot T(\pmb\vartheta )$
is the semitensor product, 
\begin{equation*}
F(\pmb\varkappa )\cdot T(\pmb\vartheta )=(F(\pmb\varkappa )\otimes
I^{\otimes }(\vartheta _{0}^{0}\sqcup \vartheta _{+}^{0}))(T(\pmb\vartheta
)\otimes I^{\otimes }(\varkappa _{0}^{-}\sqcup \varkappa _{0}^{0})),
\end{equation*}%
which is the usual product $F(\mathbf{x})T(\pmb\vartheta )$, if $\dim 
\mathcal{E}=1$. As it follows from theorem 3, the solution $U^{t}=\epsilon
\left( F_{[0,t)}\right) $ is unitary, if the triangular matrix-function $%
\mathbf{F}(x)=\exp \{-\mathrm{i}\mathbf{H}(x)\}$ is pseudo unitary, that is
the Hamiltonian matrix-function $\mathbf{H}=[H_{\nu }^{\mu }]$ is pseudo
Hermitian $H^{{\star }}(x)=H(x)$ for almost all $x\in X$: 
\begin{equation*}
H_{0}^{0{\ast }}=H_{0}^{0},H_{+}^{0{\ast }}=H_{0}^{-},H_{0}^{-{\ast }%
}=H_{+}^{0},H_{+}^{-{\ast }}=H_{+}^{-},
\end{equation*}%
$(H_{\nu }^{\mu }=0$ for $\mu =+,$ or $\nu =-)$. One can easily find the
powers $\mathbf{H}^{n}$ of the triangular matrix $\mathbf{H}$: $\mathbf{H}%
^{0}=\mathbf{I}$, $\mathbf{H}^{1}=\mathbf{H}\ $, $\mathbf{H}^{2}$ is defined
by the table 
\begin{equation*}
\mathbf{H}^{2}=%
\begin{pmatrix}
H_{0}^{-}H_{0}^{0}, & H_{0}^{-}H_{+}^{0}\cr H_{0}^{0}H_{0}^{0}, & 
H_{0}^{0}H_{+}^{0}%
\end{pmatrix}%
,\ \mathbf{H}^{n+2}=%
\begin{pmatrix}
H_{0}^{-}H_{0}^{0n-1}, & H_{0}^{-}H_{0}^{0n}H_{+}^{0}\cr H_{0}^{0n+2}, & 
H_{0}^{0n-1}H_{+}^{0}%
\end{pmatrix}%
,n=1,2,\dots 
\end{equation*}%
and $\mathbf{F}=\sum_{n=0}^{\infty }(-\mathrm{i}\mathbf{H})^{n/_{n!}}$ as
the triangular matrix $F_{\nu }^{\mu }=0,\mu >\nu ,F_{-}^{-}=1=F_{+}^{+}$,%
\begin{eqnarray*}
F_{0}^{0} &=&e^{-\mathrm{i}H_{0}^{0}},\;\;\;\;\;F_{+}^{-}=H_{0}^{-}\left[
\left( e^{-\mathrm{i}H_{0}^{0}}-I_{0}^{0}+\mathrm{i}H_{0}^{0}\right)
/H_{0}^{0}H_{0}^{0}\right] H_{+}^{0}-\mathrm{i}H_{+}^{-} \\
F_{+}^{0} &=&\left[ \left( e^{-\mathrm{i}H_{0}^{0}}-I_{0}^{0}\right)
/H_{0}^{0}\right] H_{+}^{0},\;\;\;\;\;F_{0}^{-}=H_{0}^{-}\left[ \left( e^{-%
\mathrm{i}H_{0}^{0}}-I_{0}^{0}\right) /H_{0}^{0}\right] \ .
\end{eqnarray*}%
Representing the conjugated operators $F_{0}^{-},F_{+}^{0}$ in the form $%
H_{0}^{-}=F^{{\ast }}H_{0}^{0}+\mathrm{i}E^{{\ast }}$, $H_{+}^{0}=H_{0}^{0}F-%
\mathrm{i}E$, where $E(x),F(x)$ are uniquely defined by $F^{{\star }}E=0$,
one can obtain the following canonical decomposition of the table $(L_{\nu
}^{\mu })$ of the generating operators $L_{\nu }^{\mu }(x)=F_{\nu }^{\mu
}(x)-I\otimes 1_{\nu }^{\mu }(x)$ of the unitary QS evolution $U^{t}$: 
\begin{equation*}
\begin{pmatrix}
L_{0}^{-} & L_{+}^{-}\cr L_{0}^{0} & L_{+}^{0}%
\end{pmatrix}%
=%
\begin{pmatrix}
F^{{\ast }}L_{0}^{0} & F^{{\ast }}L_{0}^{0}F\cr L_{0}^{0} & L_{0}^{0}F%
\end{pmatrix}%
+%
\begin{pmatrix}
E^{{\ast }} & {\frac{1}{2}}E^{{\ast }}E\cr0 & -E%
\end{pmatrix}%
+%
\begin{pmatrix}
0 & -\mathrm{i}H\cr0 & 0%
\end{pmatrix}%
\ ,
\end{equation*}%
\begin{equation*}
H=H_{+}^{-}-FH_{0}^{0}F^{{\ast }},\;\;\;\;\ \ \ \;L_{0}^{0}=\exp \{-\mathrm{i%
}H_{0}^{0}\}-I_{0}^{0}\ .
\end{equation*}%
Each of these three tables $\mathbf{L}_{i},\ i=1,2,3$ corresponds to a
pseudounitary triangular matrix $\mathbf{F}_{i}=\mathbf{I}+\mathbf{L}_{i}$,
satisfying the condition $\dprod_{i=1}^{3}\mathbf{F}_{i}=\mathbf{I}%
+\sum_{i=1}^{3}\mathbf{L}_{i}=\mathbf{F}$ due to the orthogonality of $%
\mathbf{L}^{i}$. The first one can be diagonalized by the pseudounitary
transformation $\mathbf{F}_{0}^{{\star }}\mathbf{L}_{1}\mathbf{F}_{0}=%
\begin{pmatrix}
0 & 0\cr L_{0}^{0} & 0%
\end{pmatrix}%
$. This defines the QS unitary evolution as the composition of three
canonical types:

1) the Poissonian type evolution, given by the diagonal matrix-function $%
\mathbf{F}(x)$, corresponding to $H_{\nu }^{\mu }=0$ for all $(\mu ,\nu
)\not=0$, for which 
\begin{equation*}
U^{t}=\epsilon (F_{[0,t)})=F_{0_{[0,t)}}^{0},\qquad F_{0}^{0}(x)=\exp \{-%
\mathrm{i}H_{0}^{0}(x)\}
\end{equation*}%
that is $U^{t}=\int^{\oplus }U^{t}(\varkappa )\mathrm{d}\varkappa $, where $%
U^{t}(\varkappa )=F_{0}^{0}(x_{n})^{t}\cdot \cdot \cdot F_{0}^{0}(x_{1})^{t}$
for any chain $\varkappa =(x_{1},\dots ,x_{n})\in \mathcal{X}$, where $%
F_{0}^{0}(x)^{t}=F_{0}^{0}(x)$, if $t(x)<t$, otherwise $F_{0}^{0}(x)^{t}=I%
\otimes 1_{0}^{0}(x)$,

2) the quantum Brownian evolution, corresponding to $H_{0}^{0}=0=H_{+}^{-}$
with $\mathrm{i}H_{+}^{0}=E=\mathrm{i}H_{0}^{-{\ast }}$, and

3) Lebesgue type evolution, corresponding to $H_{\nu }^{\mu }=0$ for all $%
(\mu ,\nu )\not=(-,+)$ for which 
\begin{equation*}
U^{t}=\epsilon (F_{[0,t)})=\hat{1}\otimes \int_{\mathcal{X}^{t}}(-\mathrm{i}%
)^{|\varkappa |}\overleftarrow{\dprod }_{x\in \varkappa }H_{+}^{-}(x)\mathrm{%
d}\varkappa ={\overleftarrow{\exp }}\left\{ -\mathrm{i}%
\int_{X^{t}}H_{+}^{-}(x)\mathrm{d}x\right\} \otimes \hat{1}\ ,
\end{equation*}%
where $\overleftarrow{\dprod }_{x\in \varkappa }H(x)=H(x_{n})\cdots H(x_{1})$
is the usual chronological product of operators $H_{+}^{-}(x)$ in $\mathcal{H%
}$ defined for any chain $\varkappa =(x_{1},\dots ,x_{n})$ by $%
t(x_{i-1})<t(x_{i})$. The sufficient conditions for the existence of the
operators $U^{t}$ as the representations $\epsilon $ of chronological
products of the elements $F_{\nu }^{\mu }(x)$ of $\exp \{-\mathrm{i}\mathbf{H%
}(x)\}$ is the local QS integrability%
\begin{equation*}
\Vert H_{0}^{0}\Vert _{\infty }^{t}<\infty ,\Vert H_{+}^{0}\Vert
_{2}^{t}=\Vert H_{0}^{-}\Vert _{2}^{t}<\infty ,\Vert H_{+}^{-}\Vert
_{1}^{t}<\infty .
\end{equation*}%
These conditions define the QS integral $\sum \Lambda _{\nu }^{\mu }\left(
t,H_{\nu }^{\mu }\right) =\Lambda ^{t}(\mathbf{H})$ as a $(\xi ^{+},\xi _{-})
$-continuous operator for any $\xi ^{+}>1>\xi _{-}$ and the QS time ordered
exponential [14] $U^{t}={\overleftarrow{\exp }}\{-\mathrm{i}\Lambda ^{t}(%
\mathbf{H})\}$ as 
\begin{equation*}
U^{t}=\Gamma _{\lbrack 0,t)}(\mathbf{L})\equiv \epsilon \left(
F_{[0,t)}\right) 
\end{equation*}%
even if $\mathbf{H}(x)$ is not pseudo Hermitian.

\section{Non-Markovian QS processes and Langevin equations}

Let $\mathcal{A}\subseteq \mathcal{B}(\mathcal{H})$ be a unital $\ast $%
-algebra of operators $A\in \mathcal{A}$, acting on a Hilbert space $%
\mathcal{H}$, and $j^{t}:\mathcal{A}\rightarrow \mathcal{B}(\mathcal{G})$ be
a family of unital $\ast $-homomorphisms, representing $\mathcal{A}$ on $%
\mathcal{G}=\mathcal{H}\otimes \mathcal{F}$ as a QS process in the sense
[12,13]: 
\begin{equation*}
j^{t}(A^{{\ast }}A)=j^{t}(A)^{{\ast }}j^{t}(A)\ ,\ j^{t}(I)=\hat{I}\ .
\end{equation*}%
We shall assume that each process $A^{t}=j^{t}(A)$ has a QS differential $%
\mathrm{d}j(A)=\mathrm{d}\Lambda (\pmb\partial (A))$ in the sense 
\begin{equation*}
j^{t}(A)=j^{0}(A)+\sum \Lambda _{\mu }^{\nu }(t,\partial _{\nu }^{\mu }(A))\
,\eqno(4.1)
\end{equation*}%
where $\pmb\partial =(\partial _{\nu }^{\mu })_{\nu =0,+}^{\mu =-,0}$ is a
family of linear maps $\partial _{\nu }^{\mu }(x):\mathcal{A}\rightarrow 
\mathcal{B}(\mathcal{G})$, depending on $x\in X$ in such a way that the
table-function $\mathbf{D}(x)=\pmb\partial (x,A)$ is QS integrable for every 
$A\in \mathcal{A}$. The QS derivative $\pmb\partial $ of the process $j$ was
introduced by Hudson [14] in the Markovian case $\pmb\partial =j\circ \pmb%
\lambda $, corresponding to the assumption $\partial _{\nu }^{\mu }(x,%
\mathcal{A})\subseteq \;j^{t(x)}(\mathcal{A})$ for all $\mu ,\nu $ and $x$.
Using the adapted QS It\^{o} formula he obtained the cohomology conditions
for the maps $\lambda _{\nu }^{\mu }$: $\mathcal{A}\rightarrow \mathcal{A}$,
which are necessary and sufficient in the constant case $\pmb\lambda (x)=\pmb%
\lambda $ for homomorphism property of $j^{t}$. These conditions can be
written simply as unital $\star $-homomorphism property for $\pmb\varphi =%
\pmb\lambda +\mathbf{j}\ ,\ \mathbf{j}(A)=A\otimes \mathbf{1}$ 
\begin{equation*}
\pmb\varphi (x,A^{{\ast }}A)=\pmb\varphi (x,A)^{{\star }}\pmb\varphi (x,A)\
,\ \pmb\varphi (x,I)=\hat{I}
\end{equation*}%
in terms of the linear maps $\pmb\varphi (x):A\rightarrow \left[ \varphi
_{\nu }^{\mu }(x,A)\right] $ into the triangular block-matrices $\mathbf{A}%
=[A_{\nu }^{\mu }]=\pmb\varphi (A)$. In the scalar case $\mathcal{E}(x)=%
\mathbb{C}$ 
\begin{equation*}
\varphi _{\nu }^{\mu }(A)=\lambda _{\nu }^{\mu }(A),\mu <\nu ;\;\varphi
_{\nu }^{\mu }(A)=\lambda _{\nu }^{\mu }(A)+A,\mu =\nu ;\;\varphi _{\nu
}^{\mu }(A)=0,\mu >\nu ,
\end{equation*}%
is defined as the sum of the map $\pmb\lambda =[\lambda _{\nu }^{\mu }]$
into the triangular matrices: $\lambda _{\nu }^{\mu }(A)=0$ if $\mu =+$ or $%
\nu =-$ and the diagonal map $\mathbf{j}=[j_{\nu }^{\mu }]$, $j_{\nu }^{\mu
}(A)=0$, $\mu \not=\nu $, $j_{\nu }^{\mu }(A)=A$, if $\mu =\nu $. As an
example one can consider the spatial $\star $-homomorphism%
\begin{equation*}
\pmb\varphi (x,A)=\mathbf{F}^{{\star }}(x)(A\otimes \mathbf{1}(x))\mathbf{F}%
(x),
\end{equation*}%
where $\mathbf{F}(x)=\left[ F_{\nu }^{\mu }(x)\right] $ is a pseudounitary
triangular matrix $F_{\nu }^{\mu }(x)=0,\mu >\nu $ with $%
F_{-}^{-}=I=F_{+}^{+}$.

We shall prove as the consequence of the theorem 4 that the
pseudo-homomorphism property of a locally QS integrable function $\pmb%
\varphi =\{\pmb\varphi (x)\}$ is also sufficient for the uniqueness and
homomorphism property (4.1) of the solution $j^{t}$ of QS Langevin equation
(4.2) with a given initial $\star $-homomorphism $j^{0}$ in nonstationary
and non-Markovian, and even in the nonadapted case.

Before doing this let us describe a decomposable operator representation $%
\pmb{\mathcal A}=\int^{\oplus }\pmb{\mathcal A}(\varkappa )\mathrm{d}%
\varkappa $ of the unital $\star $-algebra $\pmb{\mathcal A}$ of relatively
bounded $\mathcal{A}$-valued operator-functions $T(\pmb\varkappa )$ in a
pseudo Hilbert space $\pmb{\mathcal G}=\mathcal{H}\otimes \pmb{\mathcal F}$.
Here $\pmb{\mathcal F}$ is the pseudo Fock space defined in [8,10] as usual
Fock space over the space $\pmb{\mathcal K}=L^{1}(X)\oplus \mathcal{K}\oplus
L^{\infty }(X)$ of $L^{p}$-integrable vector-function $\mathbf{k}(x)=[k^{\mu
}(x)|\;\mu =-,0,+]$ with the pseudoscalar product 
\begin{equation*}
(\mathbf{k}|\mathbf{k})=\langle k^{-}|k^{+}\rangle +\Vert k^{0}\Vert
^{2}+\langle k^{+}|k^{-}\rangle =\langle \mathbf{k}|\mathbf{gk}\rangle ,
\end{equation*}%
$k^{-}\in L^{1}(X),k^{+}\in L^{\infty }(X),k^{0}\in \mathcal{K}$. In general
case, when $\mathcal{K}$ is $L^{2}$-integral $\mathcal{K}=\int^{\oplus }%
\mathcal{E}(x)\mathrm{d}x$ of Euclidean (Hilbert) spaces $\mathcal{E}%
(x),x\in X$ with%
\begin{equation*}
k\in \mathcal{K}\Longleftrightarrow k(x)\in \mathcal{E}(x),\int \Vert
k(x)\Vert ^{2}\mathrm{d}x<\infty ,
\end{equation*}%
$\pmb{\mathcal F}$ consists of all integrable in the sense%
\begin{equation*}
\Vert \mathbf{k}\Vert =\sup\limits_{\varkappa ^{+}}[\int (\int \Vert
k(\varkappa ^{-},\varkappa ^{0},\varkappa ^{+})\Vert \mathrm{d}\varkappa
^{-})^{2}\mathrm{d}\varkappa ^{0}]^{1/2}<\infty 
\end{equation*}%
tensor-functions $k(\varkappa ^{-},\varkappa ^{0},\varkappa ^{+})\in 
\mathcal{E}^{\otimes }(\varkappa ^{0})=\otimes _{x\in \varkappa ^{0}}%
\mathcal{E}(x)$ of three chains $\varkappa ^{\mu }\in \mathcal{X}$, $\mu
=-,0,+$ with the pseudoscalar product $(\mathbf{k}|\mathbf{k})=\langle 
\mathbf{k}|\mathbf{g}^{\otimes }\mathbf{k}\rangle $, 
\begin{equation*}
(\mathbf{k}|\mathbf{k})=\iiint \langle k(\varkappa ^{-},\varkappa
^{0},\varkappa ^{+})|k(\varkappa ^{+},\varkappa ^{0},\varkappa ^{-})\rangle 
\mathrm{d}\varkappa ^{-}\mathrm{d}\varkappa ^{0}\mathrm{d}\varkappa ^{+},%
\eqno(4.2)
\end{equation*}%
Taking into account that $\bigcap \varkappa ^{\mu }=\emptyset $ almost
everywhere for the continuous measure $\mathrm{d}\varkappa $ on $\mathcal{X}$%
, and that for any $a\in \mathcal{H}\otimes \pmb{\mathcal F}$ 
\begin{equation*}
(a|a)=\int \sum_{\sqcup \varkappa ^{\mu }=\varkappa }\langle a(\varkappa
^{-},\varkappa ^{0},\varkappa ^{+})|a(\varkappa ^{+},\varkappa
^{0},\varkappa ^{-})\rangle \mathrm{d}\varkappa \equiv (\mathbf{a}|\mathbf{a}%
)
\end{equation*}%
one can consider the space $\mathbf{\mathcal{G}}$ as a pseudo Hilbert
integral $\int^{\oplus }\mathbf{\mathcal{G}}(\varkappa )\mathrm{d}\varkappa $
of tensor-functions $\mathbf{a}(\varkappa )\in \mathcal{H}\otimes %
\pmb{\mathcal E}^{\otimes }(\varkappa )=\mathbf{\mathcal{G}}(\varkappa ),%
\pmb{\mathcal E}(x)=\mathbb{C}\oplus \mathcal{E}(x)\oplus \mathbb{C}$ with
values in direct sums $\mathbf{a}(\varkappa )=\oplus _{\sqcup \varkappa
^{\mu }=\varkappa }a(\varkappa ^{-},\varkappa ^{0},\varkappa ^{+})$ over all
the partitions $\varkappa =\varkappa ^{-}\sqcup \varkappa ^{0}\sqcup
\varkappa ^{+}$ of $\varkappa \in \mathcal{X}$. The operator-valued
functions $T(\pmb\varkappa )$ of $\pmb\varkappa =(\varkappa _{\nu }^{\mu })$
are unequally defined by decomposable operator $\mathbf{T}=\int^{\oplus }%
\mathbf{T}(\varkappa )\mathrm{d}\varkappa $ in $\pmb{\mathcal G}$ acting as 
\begin{equation*}
\lbrack \mathbf{T}\mathbf{a}](\varkappa )=\sum_{\sqcup \varkappa ^{\mu
}=\varkappa }\oplus \lbrack Ta](\varkappa ^{-},\varkappa ^{0},\varkappa
^{+})=\mathbf{T}(\varkappa )\mathbf{a}(\varkappa ),
\end{equation*}%
\begin{equation*}
\lbrack Ta](\varkappa ^{-},\varkappa ^{0},\varkappa ^{+})=\sum_{\sqcup _{\nu
\geq \mu }\varkappa _{\nu }^{\mu }=\varkappa ^{\mu }}^{\mu =-,0,+}T%
\begin{pmatrix}
\varkappa _{0}^{-} & \varkappa _{+}^{-}\cr\varkappa _{0}^{0} & \varkappa
_{+}^{0}%
\end{pmatrix}%
a(\varkappa _{-}^{-},\varkappa _{0}^{-}\sqcup \varkappa _{0}^{0},\varkappa
_{+}^{-}\sqcup \varkappa _{+}^{0}\sqcup \varkappa _{+}^{+})\eqno(4.3)
\end{equation*}%
due to $\sqcup \varkappa ^{\mu }=\sqcup \varkappa _{\nu }$ for $\varkappa
^{\mu }=\sqcup _{\nu \geq \mu }\varkappa _{\nu }^{\mu }$, $\varkappa _{\nu
}=\sqcup _{\mu \leq \nu }\varkappa _{\nu }^{\mu }$. It is easy to check that
the pseudo conjugated operator $\mathbf{T}^{{\star }}=\int^{\oplus }\mathbf{T%
}(\varkappa )^{{\star }}\mathrm{d}\varkappa $ with respect to the pseudo
scalar product (4.2) is also decomposable: $\mathbf{T}^{{\star }}(\varkappa
)=\mathbf{T}(\varkappa )^{{\star }}$, and is defined as in (4.3), by $T^{{%
\star }}(\pmb\varkappa )$ and the product $(\mathbf{T}^{{\star }}\mathbf{T}%
)(\varkappa )=\mathbf{T}^{{\star }}(\varkappa )\mathbf{T}(\varkappa )$
corresponds to the product (2.5). Moreover, the Fock representation (2.1) of
the operator $\star $-algebra $\int^{\oplus }\pmb{\mathcal A}(\varkappa )%
\mathrm{d}\varkappa $ can be described as a spatial transformation $\epsilon
(T)=J^{{\star }}\mathbf{T}J$, where $J$ is a pseudoisometric operator $%
(Ja|Ja)=\Vert a\Vert ^{2}$, with $J^{{\star }}\colon \langle J^{{\star }}%
\mathbf{a}|a\rangle =(\mathbf{a}|Ja)$, acting as 
\begin{equation*}
\lbrack Ja](\varkappa ^{-},\varkappa ^{0},\varkappa ^{+})=\delta _{\emptyset
}(\varkappa ^{-})a(\varkappa ^{0}),\;\;\;\;\;[J^{{\star }}\mathbf{a}%
](\varkappa )=\int a(\varkappa ^{-},\varkappa ,\emptyset )\mathrm{d}%
\varkappa ^{-}\ ,
\end{equation*}%
($\delta _{\emptyset }$ means the vacuum function $\delta _{\emptyset
}(\varkappa )=0$, if $\varkappa \not=\emptyset ,\delta _{\emptyset
}(\emptyset )=1$). One can consider $J^{{\star }}\mathbf{T}J$ as a weak
limit $t\rightarrow \infty $ of the operators $J_{[0,t)}^{{\star }}\mathbf{T}%
J_{[0,t)}$, well defined on $\mathcal{G}$ as $J_{[0,t)}=\int_{\mathcal{X}%
^{t}}^{\oplus }J(\varkappa )\mathrm{d}\varkappa $ due to 
\begin{equation*}
\Vert J_{[0,t)}a\Vert ^{2}=\iiint_{\varkappa ^{\mu }\in \lbrack 0,t)}\delta
_{\emptyset }(\varkappa ^{-})\Vert a(\varkappa ^{0})\Vert ^{2}\mathrm{d}%
\varkappa ^{-}\mathrm{d}\varkappa ^{0}\mathrm{d}\varkappa ^{+}\leq \;\mathrm{%
e}^{t}\Vert a\Vert ^{2},
\end{equation*}%
and to prove directly the property 
\begin{equation*}
J^{{\star }}\mathbf{T}^{{\star }}JJ^{{\star }}\mathbf{T}J=J^{{\star }}%
\mathbf{T}^{{\star }}\mathbf{T}J
\end{equation*}%
corresponding to the multiplicativity property of $\epsilon $.

Let $U^{t}=J^{{\star }}\mathbf{T}^{t}J$ be the solution of the QS evolution
equation (3.4) with $\mathbf{S}=\mathbf{J}^{{\star }}\mathbf{FJ}$, where $%
\mathbf{J}=J\otimes \mathbf{1}(x),1_{\nu }^{\mu }=0$, if $\mu \not=\nu $, $%
1_{-}^{-}=1=1_{+}^{+}$, $1_{0}^{0}(x)=I(x)$ is the identity operator in $%
\mathcal{E}(x)$, and let $\pmb\tau ^{t}\colon \pmb{\mathcal A}\rightarrow %
\pmb{\mathcal A},\upsilon ^{t}\colon \mathcal{B}\rightarrow \mathcal{B}$ be
the corresponding transformations $\pmb\tau (\mathbf{A})=\mathbf{T}^{{\star }%
}\mathbf{AT}$, $\upsilon (B)=U^{{\star }}BU$ of the algebras of relatively
bounded operators respectively in $\pmb{\mathcal
G}$ and ${\mathcal{G}}$. Then one can obtain, denoting $\mathbf{E}=JJ^{{%
\star }}$,%
\begin{equation*}
J^{{\star }}(\pmb\tau ^{t}(A))J=J^{{\star }}\mathbf{T}^{t\star }\mathbf{AT}%
^{t}J=J^{{\star }}\mathbf{T}^{t\star }\mathbf{EAE}\mathbf{T}^{t}J=\upsilon
^{t}(J^{{\star }}\mathbf{A}J),
\end{equation*}%
that is the QS process $\epsilon ^{t}=\upsilon ^{t}\circ \epsilon $ over the 
$\star $-algebra $\pmb{\mathcal A}$ is the composition $\epsilon
^{t}=\epsilon \circ \pmb\tau ^{t}$ of the representation $\epsilon $ and $%
\pmb\tau ^{t}=\int^{\oplus }\pmb\tau ^{t}(\varkappa )\mathrm{d}\varkappa $,
where $\pmb\tau (\varkappa ,\mathbf{A})=\mathbf{T}(\varkappa )^{{\star }}%
\mathbf{AT}(\varkappa )$ is defined due to (3.6) as chronological
compositions 
\begin{equation*}
\phi _{\lbrack 0,t)}(\varkappa ,\mathbf{A})=\mathbf{F}_{t_{1}}^{{\star }%
}(\varkappa )\cdots \mathbf{F}_{t_{m}}^{{\star }}(\varkappa )\mathbf{A}%
\mathbf{F}_{t_{m}}(\varkappa )\cdots \mathbf{F}_{t_{1}}(\varkappa )=\left[
\circ _{x\in \varkappa ^{t}}^{\rightarrow }\pmb\phi _{t(x)}(\varkappa )%
\right] (\mathbf{A})
\end{equation*}%
of the maps $\pmb\phi _{t(x)}(\varkappa \sqcup x)=\pmb\phi (x,\varkappa )$,
and $\pmb\tau ^{0}(\varkappa ,A)=\mathbf{T}^{0}(\varkappa )^{{\star }}%
\mathbf{A}\mathbf{T}^{0}(\varkappa )$: $\pmb\tau (\varkappa )=\pmb\tau
^{0}(\varkappa )\circ \pmb\phi _{\lbrack 0,t)}(\varkappa )$, where 
\begin{equation*}
\pmb\phi (\dot{\mathbf{A}})=\int^{\oplus }\pmb\phi (\varkappa ,\dot{\mathbf{A%
}}(\varkappa ))\mathrm{d}\varkappa =\mathbf{F}^{{\star }}\dot{\mathbf{A}}%
\mathbf{F},\dot{\mathbf{A}}\in \pmb{\dot{\mathcal A}}=\int^{\oplus
}\int^{\oplus }\pmb{\mathcal A}(\varkappa \sqcup x)\mathrm{d}\varkappa 
\mathrm{d}x.
\end{equation*}

Moreover, if a $\mathcal{B}$-valued process $B^{t}=\epsilon (\mathbf{A}^{t})$
has a QS differential $\mathrm{d}B=\mathrm{d}\Lambda (\mathbf{D})$, then the
transformed process $\hat{B}^{t}=\upsilon ^{t}(B^{t})$ satisfies the QS
equation 
\begin{equation*}
\hat{B}^{t}=\hat{B}^{0}+\Lambda ^{t}(\hat{\pmb\sigma }(\mathbf{G})-\hat{%
\mathbf{B}}),\mathbf{G}=\mathbf{B}+\mathbf{D}\ ,\eqno(4.4)
\end{equation*}%
where $\hat{\mathbf{B}}(x)=\mathbf{U}^{{\star }}(x)\mathbf{B}(x)\mathbf{U}%
(x)\equiv \pmb\upsilon (x,\mathbf{B}(x))$, $\mathbf{B}(x)=\mathbf{J}^{{\star 
}}\dot{\mathbf{A}}^{t(x)}(x)\mathbf{J}$, $\pmb\sigma (\mathbf{G})=\mathbf{S}%
^{{\star }}\mathbf{GS}$ as it fallows directly from (3.4) and the main
formula (2.9): 
\begin{equation*}
\mathrm{d}\left( U^{{\star }}BU\right) =\mathrm{d}\Lambda \left( \mathbf{U}^{%
{\star }}\mathbf{S}^{{\star }}(\mathbf{B}+\mathbf{D})\mathbf{SU}-\mathbf{U}^{%
{\star }}\mathbf{BU}\right) \ .
\end{equation*}%
In particular case $\mathbf{D}=0$ this gives the QS Langevin (non adapted)
equation for the QS process $\epsilon ^{t}\colon \pmb{\mathcal A}\rightarrow 
\mathcal{B}$, written in the differential form as 
\begin{equation*}
\mathrm{d}\epsilon ^{t}(\mathbf{A})=\mathrm{d}\Lambda ^{t}(\pmb\epsilon
\circ \pmb\phi (\dot{\mathbf{A}})-\pmb\epsilon (\dot{\mathbf{A}}))=\mathrm{d}%
\Lambda ^{t}(\pmb\epsilon \circ \pmb\lambda (\dot{\mathbf{A}}))\ ,\eqno(4.5)
\end{equation*}%
where $\dot{\mathbf{A}}(x)=\int^{\oplus }\mathbf{A}(x\cup \varkappa )\mathrm{%
d}\varkappa \in \dot{\pmb{\mathcal A}}(x),$%
\begin{equation*}
\pmb\lambda (\dot{\mathbf{A}})=\pmb\phi (\dot{\mathbf{A}})-\dot{\mathbf{A}},%
\pmb\epsilon (x,\dot{\mathbf{A}})=\pmb\upsilon (x,\mathbf{J}^{{\star }}\dot{%
\mathbf{A}}\mathbf{J}),
\end{equation*}%
and $(\pmb\epsilon \circ \pmb\phi )(\dot{\mathbf{A}})=\pmb\upsilon \circ \pmb%
\sigma (\mathbf{J}^{{\star }}\dot{\mathbf{A}}\mathbf{J})$ due to 
\begin{equation*}
\mathbf{J}^{{\star }}\pmb\phi (\dot{\mathbf{A}})\mathbf{J}=\mathbf{J}^{{%
\star }}\mathbf{F}^{{\star }}\dot{\mathbf{A}}\mathbf{FJ}=\mathbf{J}^{{\star }%
}\mathbf{F}^{{\star }}\mathbf{E}\dot{\mathbf{A}}\mathbf{EFJ}=\mathbf{S}^{{%
\ast }}\mathbf{J}^{{\ast }}\dot{\mathbf{A}}^{{\ast }}\mathbf{J}\mathbf{S}=%
\pmb\sigma (\mathbf{J}^{{\star }}\dot{\mathbf{A}}\mathbf{J})\ .
\end{equation*}%
The restriction of the equation (4.7) on the $\star $-subalgebra $\mathcal{A}%
\otimes \mathbf{1}^{\otimes }\in \pmb{\mathcal A},\mathbf{1}^{\otimes
}(\varkappa )=\otimes _{x\in \varkappa }\mathbf{1}(x)$ gives the (non
Markovian) Langevin equation (4.2) for $j^{t}=\epsilon ^{t}\circ \mathbf{j}$%
, $\mathbf{j}(A)=A\otimes \mathbf{1}^{\otimes }$ with the QS derivative 
\begin{equation*}
\pmb\partial =\pmb\epsilon \circ \pmb\lambda \circ \mathbf{j},\;\;\;\;\;%
\mathbf{j}(x,A)=A\otimes \mathbf{1}(x)\otimes \mathbf{1}^{\otimes }
\end{equation*}%
over $\mathcal{A}\subset \mathcal{B}(\mathcal{H})$.

Let us find the solution of the general Langevin QS equation (4.5) with
nonspatial map $\pmb\varphi $. It is given by the following theorem, which
is an analog of the theorem 3 for the maps $\pmb\lambda ,\pmb\phi ,\pmb\tau $
instead of the corresponding operators $\mathbf{L},\mathbf{F},\mathbf{T},$.

\begin{theorem}
The QS equation (4.5), written as $\epsilon ^{t}=\epsilon ^{0}+\Lambda
^{t}\circ \pmb\delta $ for all $\mathbf{A}\in \pmb{\mathcal A}$ with 
\begin{equation*}
\delta _{\nu }^{\mu }(x,\mathbf{A})=\epsilon _{\nu }^{\mu }(x,\pmb\lambda (x,%
\dot{\mathbf{A}})),\;\;\;\;\epsilon ^{0}(\mathbf{A})=J^{{\star }}\pmb\tau
^{0}(\mathbf{A})J
\end{equation*}%
is defined by linear decomposable maps $\pmb\lambda (x)\colon \dot{%
\pmb{\mathcal A}}(x)\rightarrow \dot{\pmb{\mathcal A}}(x)$ with $\lambda
_{\nu }^{\mu }(x,\mathbf{A})=0$, if $\mu =+$ or $\nu =-$ and $\pmb\tau
^{0}\colon \pmb{\mathcal A}\rightarrow \pmb{\mathcal A}$ as the
representation 
\begin{equation*}
\epsilon ^{t}(\mathbf{A})=J^{{\star }}\pmb\tau ^{t}(\mathbf{A})J\ ,\ \pmb%
\epsilon (x,\dot{\mathbf{A}}(x))=\mathbf{J}^{{\star }}\dot{\pmb\tau }%
^{t(x)}(x,\dot{\mathbf{A}}(x))\mathbf{J}
\end{equation*}%
of the recurrences 
\begin{equation*}
\pmb\tau ^{t_{+}}(\varkappa )=\pmb\tau ^{t}(\varkappa )\circ \pmb\phi
_{t(x)}(\varkappa ),\quad x\in \varkappa \in \mathcal{X}\ ,\eqno(4.6)
\end{equation*}%
$t\in (t_{-}(\varkappa ),t(x)]$, $t_{+}\in (t(x),t_{+}(x)]$, where $t_{\pm
}=t_{m\pm 1}$ for a chain $\varkappa =(x_{1},\dots ,x_{m},\dots )$, $x=x_{m}$%
, $t_{m}=t(x_{m})$, and 
\begin{equation*}
\pmb\phi _{t(x)}(\varkappa \sqcup x,\mathbf{A})=\pmb\lambda (x,\varkappa ,%
\mathbf{A})+\mathbf{A}\equiv \pmb\phi (x,\varkappa ,\mathbf{A}),\mathbf{A}%
\in \pmb{\mathcal A}(x\sqcup \varkappa )\ .
\end{equation*}%
The recurrency (4.6) with initial condition $\pmb\tau ^{t}(\varkappa )=\pmb%
\tau ^{0}(\varkappa )$ for all $t\in \lbrack 0,t_{1}]$ has the unique
solution 
\begin{equation*}
\pmb\tau ^{t}(\varkappa )=\pmb\tau ^{0}(\varkappa )\circ \pmb\phi _{\lbrack
0,t)}(\varkappa )\ ,\ \pmb\phi _{\lbrack 0,t)}=\circ _{0\leq
s<t}^{\rightarrow }\pmb\phi _{s}\eqno(4.7)
\end{equation*}%
defined for every $\varkappa =(x_{1},\dots ,x_{m},\dots ),t\in
(t_{m-1},t_{m}]$ by the chronological composition $\circ _{x\in \varkappa
^{t}}^{\rightarrow }\pmb\phi _{t(x)}=\pmb\phi _{t_{1}}\circ \dots \circ \pmb%
\phi _{t_{m-1}}$ of $\pmb\phi _{t(x)}(\varkappa )=\pmb\phi (x,\varkappa
\backslash x)$ for $x\in \varkappa ^{t}=\{x\in \varkappa |t(x)<t\}$, $\pmb%
\phi _{s}(\varkappa )=\mathbf{i}(\varkappa )$, $i(\varkappa )$ is the
identity map $\pmb{\mathcal A}(\varkappa )\rightarrow \pmb{\mathcal A}%
(\varkappa )$, if $s\notin t(\varkappa )$. The solution $\{\epsilon ^{s}\}$
of (4.5) is Hermitian: $\epsilon ^{s}(\mathbf{A}^{{\star }})=\epsilon ^{s}(%
\mathbf{A})^{{\ast }}$ up to a $t>0$, if $\pmb\tau ^{0}$ and $\pmb\phi (x)$
are pseudo Hermitian: 
\begin{equation*}
\pmb\tau ^{0}(\mathbf{A}^{{\star }})=\pmb\tau ^{0}(\mathbf{A})^{{\star }},%
\pmb\phi (x,\dot{\mathbf{A}}^{{\star }}(x))=\pmb\phi (x,\dot{\mathbf{A}}%
(x))^{{\star }},x\in X^{t},
\end{equation*}%
and is (unital, faithful) QS-process, representing $\pmb{\mathcal A}$, if $%
\pmb\tau ^{0}\colon \pmb{\mathcal A}\rightarrow \pmb{\mathcal A}$ and $\pmb%
\phi (x)\colon \dot{\pmb{\mathcal A}}\rightarrow \dot{\pmb{\mathcal A}}(x)$
are (unital) $\star $-endomorphisms (automorphisms) of the algebras $%
\pmb{\mathcal A}$ and $\dot{\pmb{\mathcal A}}(x)$ for almost all $x\in X^{t}$%
.
\end{theorem}

\begin{proof}
We look for the solution of the equation (4.5) as for the representation $%
\epsilon ^{t}(\mathbf{A})=J^{{\star }}\mathbf{A}^{t}J$ of a process $\mathbf{%
A}^{t}=\pmb\tau ^{t}(\mathbf{A})$, transforming the $\star $-algebra $%
\pmb{\mathcal A}$. From definition of $\pmb\epsilon (x)$ and $\dot{\mathbf{A}%
}(x,\varkappa )=\mathbf{A}(\varkappa \sqcup x)$ we obtain 
\begin{equation*}
\pmb\epsilon (x,\dot{\mathbf{A}}(x))=\mathbf{J}^{{\star }}\hat{\mathbf{A}}(x)%
\mathbf{J},\hat{\mathbf{A}}(x)=\dot{\pmb\tau }^{t(x)}(x,\dot{\mathbf{A}}(x))=%
\dot{\mathbf{A}}^{t(x)}(x),
\end{equation*}%
and $\Lambda ^{t}(\pmb\delta (A))=\Lambda ^{t}(\mathbf{J}^{{\star }}\hat{%
\mathbf{L}}\mathbf{J})=J^{{\star }}N^{t}(\hat{\mathbf{L}})J$, where $\hat{%
\mathbf{L}}(x)=\dot{\pmb\tau }^{t(x)}(x,\dot{\mathbf{L}}(x)),\dot{\mathbf{L}}%
(x)=\pmb\lambda (x,\dot{\mathbf{A}}(x))$. This gives the equation (4.6) in
the integral form as the representation $J^{{\star }}(\dot{\mathbf{A}}%
^{0}+N^{t}(\hat{\mathbf{L}}))J=J^{{\star }}\mathbf{A}^{t}J$ of the equations 
\begin{align*}
\pmb\tau ^{0}(\varkappa ,\mathbf{A}(\varkappa ))& +\sum_{x\in \varkappa ^{t}}%
\dot{\pmb\tau }^{t(x)}(x,\varkappa \backslash x,\pmb\lambda (x,\varkappa
\backslash x,\dot{\mathbf{A}}(x,\varkappa \backslash x)))= \\
\pmb\tau ^{0}(\varkappa ,\mathbf{A}(\varkappa ))& +\sum_{x\in \varkappa ^{t}}%
\pmb\tau ^{t(x)}(\varkappa ,\pmb\lambda _{t(x)}(\varkappa ,\mathbf{A}%
(\varkappa )))=\pmb\tau ^{t}(\varkappa ,\mathbf{A}(\varkappa ))
\end{align*}%
where $\pmb\lambda _{t(x)}(\varkappa )=\pmb\lambda (x,\varkappa \backslash x)
$ for a $x\in \varkappa $. Denoting $\pmb\lambda _{t}(\mathbf{A})+\mathbf{A}$
as $\pmb\phi _{t}(\mathbf{A})$, we obtain the recurrency (4.6) for $\pmb\tau
^{t}(\varkappa )=\pmb\tau _{m}(\varkappa ),m=|\varkappa ^{t}|$ supposing the
linearity of the maps $\pmb\tau ^{t}$: 
\begin{equation*}
\pmb\tau _{m}(\varkappa )=\pmb\tau _{0}(\varkappa )+\sum_{k=1}^{m}(\pmb\tau
_{k-1}(\varkappa )\circ \pmb\phi _{t_{k}}(\varkappa )-\pmb\tau
_{k-1}(\varkappa ))=\pmb\tau _{m-1}(\varkappa )\circ \pmb\phi
_{t_{m}}(\varkappa ).
\end{equation*}%
This recurrency has the unique solution (4.7), which is linear as the
composition of the linear maps $\pmb\tau ^{0}$ and $\pmb\phi _{t}$, what
proves the uniqueness of the solution $\epsilon ^{t}=\epsilon \circ \pmb\tau
^{t}$ of the equation (4.5).

If the maps $\pmb\tau ^{0}(\varkappa )$ and $\pmb\phi (x,\varkappa )=\pmb%
\phi _{t(x)}(\varkappa \sqcup x)$ are pseudo Hermitian, then the composition 
$\pmb\tau ^{t}(\varkappa )$ is also pseudo Hermitian, and if they satisfy
the (unital) $\star $-endomorphism (automorphism) property 
\begin{align*}
\pmb\tau ^{0}(\mathbf{A}^{{\star }}\mathbf{A})& =\pmb\tau ^{0}(\mathbf{A})^{{%
\star }}\pmb\tau ^{0}(\mathbf{A}),\pmb\phi (x,\dot{\mathbf{A}}^{{\star }}%
\dot{\mathbf{A}})=\pmb\phi (x,\dot{\mathbf{A}})^{{\star }}\pmb\phi (x,%
\mathbf{A}) \\
(\pmb\tau ^{0}(\mathbf{A}^{-1})& =\pmb\tau ^{0}(\mathbf{A})^{-1},\ \pmb\phi
(x,\dot{A}^{-1})=\pmb\phi (x,\dot{A})^{-1})
\end{align*}%
for $x\in X^{t}$, then the compositions (4.7) have obviously the same
properties. This proves the Hermiticity and (unital) homomorphism
(isomorphism) property for the map $\epsilon ^{t}\colon \mathbf{A}\in %
\pmb{\mathcal A}\rightarrow J^{{\star }}\pmb\tau ^{t}(\mathbf{A})J$.
\end{proof}

Let us denote by $\pmb{\mathcal A}^{t}\subseteq \pmb{\mathcal A}$ the $\star 
$-subalgebra of relatively bounded operators $\mathbf{A}=\int^{\oplus }%
\mathbf{A}(\varkappa )\mathrm{d}\varkappa $ with $\mathbf{A}(\varkappa )=%
\mathbf{A}(\varkappa ^{t})\otimes \mathbf{1}^{\otimes }(\varkappa _{\lbrack
t})$, and by $\mathcal{B}^{t}\subseteq \mathcal{B}$ the corresponding
algebra of operators $B=B^{t}\otimes \hat{1}_{[t}$ with $B^{t}$, acting in $%
\mathcal{G}^{t}$. The adapted QS process $\epsilon ^{t}$ over $\pmb{\mathcal
A}$ is defined by the condition $\epsilon ^{t}(\pmb{\mathcal A}%
^{t})\subseteq \mathcal{B}^{t}$ for almost all $t$, and corresponds to the
adapted QS evolution $\upsilon ^{t}\colon \mathcal{B}\rightarrow \mathcal{B}%
,\upsilon ^{t}(\mathcal{B}^{t})\subseteq \mathcal{B}^{t}$, described as $%
\upsilon ^{t}(B)=\epsilon ^{t}(\mathbf{A})$ for $B=J^{{\star }}\mathbf{A}J$.

\begin{corollary}
The QS process $\epsilon ^{t}$, defined by the equation (4.6), is adapted,
if $\pmb\tau ^{0}(\pmb{\mathcal A}^{0})\subset \pmb{\mathcal A}^{0}$ and $%
\pmb\phi (x,\varkappa )=\pmb\phi (x,\varkappa ^{t(x)})\otimes \mathbf{i}%
(\varkappa _{\lbrack t(x)})$ for almost all $x\in X$. In that case the QS
evolution $\upsilon ^{t}$ is defined by adapted map $\pmb\sigma (x)\colon 
\mathbf{J}^{{\star }}\dot{\pmb{\mathcal A}}^{t(x)}(x)\mathbf{J}\rightarrow 
\mathbf{J}^{{\star }}\dot{\pmb{\mathcal A}}^{t(x)}(x)\mathbf{J}$; the
transformed adapted process $\hat{B}^{t}=\upsilon ^{t}(B^{t})$, with $B^{t}$%
, having the derivative $\mathbf{D}(x)\in \mathbf{J}^{{\star }}\dot{%
\pmb{\mathcal A}}^{t(x)}(x)\mathbf{J}$, satisfies the QS differential
equation 
\begin{equation*}
\mathrm{d}\upsilon ^{t}(B^{t})=\upsilon ^{t}[\mathrm{d}\Lambda ^{t}(\pmb%
\sigma (B\otimes \mathbf{1})+\pmb\sigma (\mathbf{D})-B\otimes \mathbf{1})]%
\eqno(4.8)
\end{equation*}%
In particular, if $B^{t}=A\otimes \hat{1}$ and $\pmb\sigma (x,A\otimes \hat{1%
}\otimes \mathbf{1})=\pmb\varphi (x,A)\otimes \hat{1}$, $\pmb\varphi (x,%
\mathcal{A})\in \pmb{\mathcal A}(x)$, where $\pmb{\mathcal A}(x)$ is the
algebra of $\mathcal{A}$-valued triangular matrices $\mathbf{A}=[A_{\nu
}^{\mu }],A_{\nu }^{\mu }=0$, if $\mu >\nu ,A_{-}^{-}=A=A_{+}^{+}$ then the
equation (4.8) has the form (4.1) in terms of $j^{t}(A)=\upsilon
^{t}(A\otimes \hat{1})$, $\partial _{\nu }^{\mu }(x)=j^{t(x)}\circ \lambda
_{\nu }^{\mu }(x)$, where 
\begin{equation*}
j^{0}(A)=\tau ^{0}(A)\otimes \hat{1},\pmb\lambda (x,A)=\pmb\varphi
(x,A)-A\otimes \mathbf{1}(x),A\in \mathcal{A}.
\end{equation*}%
If the maps $\varphi _{\nu }^{\mu }(x)$ are locally $L^{p}$-integrable in
the sense 
\begin{equation*}
\Vert \lambda _{0}^{0}\Vert _{\infty }^{t}<\infty ,\Vert \lambda
_{+}^{0}\Vert _{2}^{t}<\infty ,\Vert \lambda _{0}^{-}\Vert _{2}^{t}<\infty
,\Vert \lambda _{+}^{-}\Vert _{1}^{t}<\infty ,\eqno(4.9)
\end{equation*}%
where $\Vert \lambda \Vert _{p}^{t}=\left( \int_{X^{t}}\sup_{A\in \mathcal{A}%
}\{\Vert \lambda (x,A)\Vert /\Vert A\Vert \}^{p}\mathrm{d}x\right) ^{1/p}$,
then the solution $j^{t}(A)=J^{{\star }}\pmb\tau ^{0}(\pmb\phi _{\lbrack
0,t)}(A\otimes \hat{1}))J$ exists as relatively bounded QS integral 
\begin{equation*}
j^{t}(A)=\Lambda _{\lbrack 0,t)}(\pmb\tau ^{0}\circ \pmb\lambda
^{\triangleright }(A)\otimes \hat{1}),\pmb\lambda ^{\triangleright
}(\varkappa )=\circ _{x\in \varkappa }^{\rightarrow }\pmb\lambda
(x,\varkappa _{t(x)})
\end{equation*}%
where $\pmb\lambda (x,\varkappa )=\pmb\lambda (x)\otimes \mathbf{i}^{\otimes
}(\varkappa ),\mathbf{i}^{\otimes }(\varkappa )$ is the identity map for
operators in $\pmb{\mathcal E}^{\otimes }(\varkappa )$ and $\pmb\tau
^{0}(\varkappa )=\tau ^{0}\otimes \mathbf{i}^{\otimes }(\varkappa )$. It has
the estimate 
\begin{equation*}
\Vert j^{t}(A)\Vert _{\xi ^{+}}^{\xi _{-}}\leq \Vert \tau ^{\circ }\Vert
\exp \{\int_{X^{t}}(\Vert \lambda _{+}^{-}(x)\Vert +(\Vert \lambda
_{+}^{0}(x)\Vert ^{2}+\Vert \lambda _{0}^{-}(x)\Vert ^{2})/2\varepsilon )%
\mathrm{d}x\}\eqno(4.10)
\end{equation*}%
for $\xi ^{+}/\xi _{-}>\mathrm{ess\sup }_{x\in X^{t}}\Vert \varphi
_{0}^{0}(x)\Vert ,\Vert A\Vert \leq 1$, and sufficiently small $\varepsilon
>0$.
\end{corollary}

Indeed, the process $\epsilon ^{t}(A^{t})=J^{{\star }}\pmb\tau ^{t}(A^{t})J$
is adapted, iff $\pmb\tau ^{t}(\varkappa ,\mathbf{A})=\pmb\tau
^{t}(\varkappa ^{t},\mathbf{A}^{t})\otimes \mathbf{1}^{\otimes }(\varkappa
_{\lbrack t})$, as it was proven in the Corollary 2. But due to $\pmb\tau
(\varkappa ,\mathbf{A})=\pmb\tau (\varkappa ,\mathbf{A}(\varkappa ))$, it is
possible only in the case $\mathbf{A}^{t}(\varkappa )=\mathbf{A}%
^{t}(\varkappa ^{t})\otimes 1^{\otimes }(\varkappa _{\lbrack t})$ and $\pmb%
\tau ^{t}(\varkappa )=\pmb\tau ^{t}(\varkappa ^{t})\otimes \mathbf{1}%
^{\otimes }(\varkappa _{\lbrack t})$, what is equivalent to the
corresponding conditions for $\pmb\tau ^{0}$ and $\pmb\phi (x)$. If $%
B^{t}=J^{{\star }}\mathbf{A}^{t}J$ is an adapted process: $\mathbf{A}%
^{t}(\varkappa )=\mathbf{A}^{t}(\varkappa ^{t})\otimes \mathbf{1}^{\otimes
}(\varkappa _{\lbrack t})$, then 
\begin{align*}
\dot{\mathbf{A}}^{t}(x,\varkappa )& =\mathbf{A}^{t}(\varkappa \sqcup x)=%
\mathbf{A}^{t}(\varkappa )\otimes \mathbf{1}(x)\ ,\qquad \qquad \ \forall
t\leq t(x) \\
\mathbf{B}(x)& =\mathbf{J}^{{\star }}\dot{\mathbf{A}}^{t(x)}(x)\mathbf{J}%
=B^{t(x)}\otimes \mathbf{1}(x)\ ,\qquad \qquad \forall x\in X
\end{align*}%
and $\hat{\mathbf{B}}(x)=\upsilon (x,\mathbf{B}(x))=\hat{B}^{t(x)}\otimes 
\mathbf{1}(x)$, where $\hat{B}^{t}=\upsilon ^{t}(B^{t})$ for the transformed
process $\pmb\upsilon (x,B\otimes \mathbf{1}(x))=\upsilon ^{t(x)}(B)\otimes 
\mathbf{1}(x)$ evaluated in $B=B^{t(x)}$. This gives the equation (4.4) for
the adapted process $\hat{B}^{t}$ in the differential form (4.8): 
\begin{equation*}
\mathrm{d}\upsilon ^{t}(B^{t})=\mathrm{d}\Lambda ^{t}(\pmb\upsilon (\pmb%
\sigma (B\otimes \mathbf{1}+\mathbf{D})-B\otimes \mathbf{1}))=\upsilon ^{t}[%
\mathrm{d}\Lambda ^{t}(\pmb\beta (B)+\pmb\sigma (\mathbf{D}))],
\end{equation*}%
where $\pmb\beta (B)=\pmb\sigma (B\otimes \mathbf{1})-B\otimes \mathbf{1}$
and $\mathrm{d}\Lambda ^{t}\circ \pmb\upsilon =\upsilon ^{t}\circ \mathrm{d}%
\Lambda ^{t}$ for the adapted evolution $\upsilon ^{t}$ due to the same
arguments, as in Corollary 1. This equation, restricted on $B^{t}=A\otimes 
\hat{1}$ has the integral form (4.1) due to $\mathbf{D}=0$ where $%
j^{t}(A)=\upsilon ^{t}(A\otimes \hat{1}),\pmb\partial (x)=\upsilon
^{t(x)}\circ \pmb\beta (x)$ because $\pmb\upsilon (x)=\upsilon
^{t(x)}\otimes \mathbf{1}(x)$. If $\pmb\beta =\pmb\lambda \otimes \hat{1}$,
where $\pmb\lambda (x,A)\in \pmb{\mathcal A}(x)$, then it can be written as 
\begin{equation*}
\mathrm{d}j^{t}(A)=\mathrm{d}\Lambda ^{t}(j(\pmb\lambda (A))=j^{t}[\mathrm{d}%
\Lambda ^{t}(\pmb\lambda (A)\otimes \hat{1})]\ .
\end{equation*}%
The solution $j^{t}(A)=J^{{\star }}\pmb\tau ^{t}(A\otimes \mathbf{1}%
^{\otimes })J$ of this equation is defined by chronological composition
(4.7) as 
\begin{equation*}
\pmb\tau ^{t}(\varkappa ,A\otimes \mathbf{1}^{\otimes })=\pmb\tau
^{0}(\varkappa ,\pmb\varphi ^{\triangleright }(\varkappa ^{t},A)\otimes 
\mathbf{1}^{\otimes }(\varkappa _{\lbrack t})),\pmb\varphi ^{\triangleright
}(\varkappa )=\circ _{x\in \varkappa }^{\rightarrow }\pmb\varphi
(x,\varkappa _{t(x)})\ ,
\end{equation*}%
where $\pmb\tau ^{0}(\varkappa )=\tau ^{0}\otimes \mathbf{i}\otimes
(\varkappa )$ is an initial map, and $\pmb\varphi (x,\varkappa )=\pmb\varphi
(x)\otimes \mathbf{i}(\varkappa )$. It can be described as $%
j^{t}(A)=\epsilon (T^{t})$ by operator-valued function 
\begin{equation*}
T^{t}(\pmb\varkappa )=\tau ^{0}[\varphi (\mathbf{x}_{1},\varphi (\mathbf{x}%
_{2},\varphi (\mathbf{x}_{m},A)))]\otimes \mathbf{1}^{\otimes }(\pmb%
\varkappa _{\lbrack t})\ ,
\end{equation*}%
$t_{m}<t\leq t_{m+1}$, corresponding in (4.2) to the decomposable $\mathbf{T}%
^{t}=\pmb\tau ^{t}(A\otimes \mathbf{1}^{\otimes })$ for a partition $\pmb%
\varkappa =(\varkappa _{\nu }^{\mu })$ of the chain $\varkappa =(x_{1},\dots
,x_{m},\dots )\in \mathcal{X}$ with $\varphi (\mathbf{x}_{\nu }^{\mu
})=\varphi _{\nu }^{\mu }(x)$. Hence, the operator $\mathbf{T}^{t}$ is
relatively bounded for $A\in \mathcal{A}$: 
\begin{equation*}
\Vert T^{t}(\pmb\varkappa )\Vert \leq \Vert \tau ^{0}\Vert \ \Vert A\Vert
\dprod_{x\in \pmb\varkappa }\Vert \varphi ^{t}(\mathbf{x})\Vert =\Vert \tau
^{0}\Vert \ \Vert A\Vert \dprod_{\nu =0,+}^{\mu =-,0}\zeta ^{t}(\varkappa
_{\nu }^{\mu }),
\end{equation*}%
where $\pmb\varphi ^{t}(x)=\pmb\varphi (x)$, if $t(x)<t$, otherwise $\pmb%
\varphi ^{t}(x)=\mathbf{i}(x)$, and $\zeta ^{t}(\varkappa _{\nu }^{\mu
})=\dprod_{x\in \varkappa _{\nu }^{\mu }}^{t(x)<t}\left\Vert \varphi _{\nu
}^{\mu }(x)\right\Vert $. This proves like in the Corollary 3 the existence
of the solution $j^{t}(A)=J^{{\star }}\mathbf{T}^{t}J=\epsilon (T^{t})$ of
the equation (4.1) as relatively bounded operator $B^{t}=j^{t}(A)\in 
\mathcal{B}$ for any $A\in \mathcal{A}$, having the estimate (4.10) for $%
\Vert A\Vert \leq 1$ in terms of $\Vert \varphi _{\nu }^{\mu }(x)\Vert =\sup
\{\Vert \varphi _{\nu }^{\mu }(A)\Vert /\Vert A\Vert \}$, $\varphi _{\nu
}^{\mu }=\lambda _{\nu }^{\mu }$ for $\mu <\nu $. It can be written in the
form of the multiple QS integral (1.5) with respect to the operator function 
$B(\pmb\varkappa )=L(\pmb\varkappa )\otimes \hat{1}$ with 
\begin{equation*}
L(\pmb\varkappa ,A)=\tau ^{0}[\lambda (\mathbf{x}_{1},\lambda (\mathbf{x}%
_{2},\ldots ,\lambda (\mathbf{x}_{n},A)\dots ))]=\tau ^{0}[\lambda
^{\triangleright }(\pmb\varkappa ,A)]
\end{equation*}%
for any partition $\sqcup _{\nu =0,+}^{\mu =-,0}\varkappa _{\nu }^{\mu
}=(x_{1},\dots ,x_{n})\in \mathcal{X}$ where $\lambda (\mathbf{x}_{\nu
}^{\mu })=\lambda _{\nu }^{\mu }(x)$ on the single point table $\mathbf{x}%
=(\varkappa _{\nu }^{\mu }),\varkappa _{\nu }^{\mu }=x$.

The solution of the nonstationary Markov Langevin equation in the form of
multiple QS integral was obtained recently in the frame work of It\^{o} QS
calculus by Lindsay and Parthasarathy [15] for more restrictive conditions
then (4.9) (finite dimensional and local bounded $\lambda _{\nu }^{\mu }$).

\begin{acknowledgement}
I wish to thank Prof. L. Accardi and F. Guerra for stimulating discussions
of the results and for the hospitality of the University of Rome
\textquotedblleft La Sapienza\textquotedblright , and \textquotedblleft
Centro Matematico V. Volterra\textquotedblright , University of Rome II,
where this paper was written.
\end{acknowledgement}

\vfill\eject

\section{References}

\begin{enumerate}
\item Hudson R.L., Parthasarathy K.R. \textit{Quantum} \textit{It\^{o}'s
formula and stochastic evolution}, Comm. Math. Phys. \textbf{93}, (1984),
301-323.

\item Meyer P.A., \textit{Elements de Probabilit\'es Quantiques}, Exposes I
a IV, Institute de Mathematique, Universit\'e Louis Pasteur Strasbourg, 1985.

\item Maassen H., \textit{Quantum Markov processes in Fock space described
by integral kernels}, in: ``Quantum Probability and Applications II'', ed.
L. Accardi and W. von Waldenfels, Lecture Notes in Mathematics, 11--36,
Berlin Heidelberg New York: Springer, 1985.

\item Lindsay J.M., Maassen H., \textit{An integral kernel approach to noise}%
, in: \textquotedblleft Quantum probability and Applications
III\textquotedblright , ed. L. Accardi and W. von Waldenfels, Lecture Notes
in Mathematics, 192--208, Berlin Heidelberg New York: Springer 1988.

\item Accardi L., Quaegebeur J., \textit{The It\^{o} Algebra of Quantum
Gaussian Fields}, Journal of Functional Analysis, \textbf{85}, N. 2 (1989),
213--263.

\item Accardi L., Fagnola F., \textit{Stochastic integration}, in: ``Quantum
Probability and Applications III'', ed. L. Accardi and W. von Waldenfels,
Lecture notes in Mathematics, 6--19, Berlin Heidelberg New York: Springer
1988.

\item Belavkin V.P., \textit{Nondemolition measurements, nonlinear filtering
and dynamic programming of quantum stochastic processes}, in: ``Modeling and
Control of Systems'', ed. A. Blaquiere, Lecture Notes in Control and
Information Sciences, \textbf{121}, 245--265, Berlin Heidelberg, New York:
Springer 1988.

\item Belavkin V.P., \textit{Nondemolition stochastic calculus in Fock space
and nonlinear filtering and control in quantum systems}, in:
\textquotedblleft Stochastic methods in Mathematics and
Physics\textquotedblright , ed. R. Gielerak and W. Karwoski, 310--324,
Singapour, New Jersey, London, Hong Kong: World Scientific, 1988.

\item Belavkin V.P., \textit{Quantum stochastic calculus and quantum
nonlinear filtering}, Centro Matematico V. Volterra, Universit\`a degli
Studi di Roma II, Preprint N. 6, Marzo, 1989.

\item Belavkin V.P., \textit{A new form and $\star $-algebraic structure of
quantum stochastic integrals in Fock space}, Seminario matematico e
fisico\dots (1989),\dots

\item Holevo A.S., \textit{Time-ordered exponentials in quantum stochastic
calculus}, Preprint N. 517 Universit\"{a}t Heidelberg, June, 1989.

\item Accardi L., Frigerio A., Lewis J.T., \textit{Quantum stochastic
processes}, publications RIMS \textbf{48}, (1982), 97--133.

\item Belavkin V.P., \textit{Reconstruction theorem for quantum stochastic
processes}, Theor. Math. Phys., \textbf{62}, N. 3, (1985), 275--289.

\item Hudson R.L., \textit{Quantum diffusions and cohomology of algebras}.
Proceedings First World congress of Bernoulli society, vol. 1, Yu. Prohorov,
V.V. Sazonov (eds.), (1987), 479--485.

\item Lindsay J.M., Parthasarathy K.R., \textit{Cohomology of power sets
with applications in quantum probability}, Comm. Math. Phys. \textbf{124},
(1989), 337--364.

\item Malliavin P., \textit{Stochastic calculus of variations and
hypoelliptic operators} in: It\^{o} K., (ed) Proc. of Int. Symp. Stoch. D.
Eqs. Kyoto 1976 pp. 195--263. Tokyo: Kinokuniya--Wiley, 1978.
\end{enumerate}

\end{document}